\documentclass[10pt,a4paper,final,1p]{amsart}
\usepackage[all]{xy}
\usepackage{xspace}
\usepackage{amssymb}
\usepackage{tensor} % to define \FP
\usepackage[hidelinks]{hyperref}

\newcommand{\FP}[2]{\tensor[_{{#1}}]{\times}{_{{#2}}}} %requires tensor package
\newcommand{\gpoidarrows}{\ensuremath{\rightrightarrows}}
\newcommand{\FB}{\ensuremath{Fb}\xspace}
\newcommand{\FBS}{\ensuremath{Fbs}\xspace}

\newcommand{\LgpdC}{\ensuremath{LGpd}\xspace}
\newcommand{\lLgpdC}{\ensuremath{lLGpd}\xspace}

\newcommand{\VD}{\ensuremath{\mathcal{V}_d}\xspace}
\newcommand{\conjB}{\ensuremath{\ti{\SG}}\xspace}

\DeclareMathOperator{\ob}{ob}

\newcommand{\DC}{\ensuremath{{\mathcal{A}_d}}}
\newcommand{\SG}{\ensuremath{G}\xspace}
\newcommand{\jgg}{\ensuremath{\mathfrak{g}}\xspace}
\newcommand{\baseRing}[1]{\ensuremath{\mathbb{#1}}}
\newcommand{\Z}{\baseRing{Z}}
\newcommand{\R}{\baseRing{R}}
\newcommand{\CC}{\ensuremath{{\mathcal{A}}}}
\newcommand{\jdef}[1]{\index{#1}\emph{#1}}
\newcommand{\HLd}[1]{\ensuremath{{h_d^{{#1}}}}}

\newcommand{\ti}[1]{\widetilde{#1}}
\newcommand{\stext}[1]{\ensuremath{\quad\text{#1}\quad}}
\newcommand{\BC}{\ensuremath{{\mathcal{B}}\xspace}}
\newcommand{\BD}{\ensuremath{{\mathcal{B}_d}}}
\newcommand{\jgsg}{\ensuremath{\mathfrak{g}}\xspace}
\newcommand{\lie}[1]{\ensuremath{Lie(#1)}}
\newcommand{\conj}{\overline}

\newcommand{\VF}{\mathfrak{X}}

\DeclareMathOperator{\im}{Im}

\theoremstyle{plain}
\newtheorem{theorem}{Theorem}[section]
\newtheorem{corollary}[theorem]{Corollary}
\newtheorem{prop}[theorem]{Proposition}
\newtheorem{proposition}[theorem]{Proposition}
\newtheorem{lemma}[theorem]{Lemma}

\theoremstyle{definition}
\newtheorem{definition}[theorem]{Definition}
\newtheorem{remark}[theorem]{Remark}
\newtheorem{example}[theorem]{Example}

\numberwithin{equation}{section}

\newcommand{\PBmap}{\ensuremath{\pi}\xspace}
\newcommand{\PBtotal}{\ensuremath{Q}\xspace}
\newcommand{\PBbase}{\ensuremath{Q/\SG}\xspace}
\newcommand{\DCp}[1]{\ensuremath{{\mathcal{A}_d^{#1}}}}
\newcommand{\jcalU}{\ensuremath{\mathcal{U}}\xspace}

\begin{document}

\title[The discrete Atiyah sequence]{Discrete connections on principal
  bundles:\\ the Discrete Atiyah Sequence}

\author{Javier Fern\'andez}
\address{Instituto Balseiro, Universidad
  Nacional de Cuyo -- C.N.E.A.\\Av. Bustillo 9500, Bariloche, R{\'\i}o
  Negro, 8400, Rep\'ublica Argentina}
\email{jfernand@ib.edu.ar}

\author{Mariana Juchani}
\address{Departamento de Matem\'atica, Facultad de Ciencias Exactas,
    Universidad Nacional de La Plata\\Centro de Matem\'atica de La
    Plata (CMaLP)\\50 y 115, La Plata, Buenos Aires, 1900, Rep\'ublica
    Argentina}
\email{marianaevaj@gmail.com}

\author{Marcela Zuccalli}
\address{Departamento de Matem\'atica, Facultad de Ciencias Exactas,
    Universidad Nacional de La Plata\\Centro de Matem\'atica de La
    Plata (CMaLP)\\50 y 115, La Plata, Buenos Aires, 1900, Rep\'ublica
    Argentina}
\email{marce@mate.unlp.edu.ar}

  %%%%%%%%%%%%%%%%%%%%%%%%%%%%%%%%%%%%%%%%%%%%%%%%%%%%%%%%%%%%%%%%%% 

  \begin{abstract}
    In this work we study discrete analogues of an exact sequence of
    vector bundles introduced by M. Atiyah in 1957, associated to any
    smooth principal $G$-bundle $\PBmap:\PBtotal\rightarrow
    \PBbase$. In the original setting, the splittings of the exact
    sequence correspond to connections on the principal bundle
    $\PBmap$. The discrete analogues that we consider here can be
    studied in two different categories: the category of fiber bundles
    with a (chosen) section, \FBS, and the category of local Lie
    groupoids, \lLgpdC. In \FBS we find a correspondence between a)
    (semi-local) splittings of the discrete Atiyah sequence (DAS) of
    $\PBmap$, b) discrete connections on the same bundle $\PBmap$, and
    c) isomorphisms of the DAS with certain fiber product extensions
    in \FBS. We see that the right splittings of the DAS (in \FBS) are
    not necessarily right splittings in \lLgpdC: we use this
    obstruction to define the discrete curvature of a discrete
    connection.  Then, there is a correspondence between the right
    splittings of the DAS in \lLgpdC and discrete connections with
    trivial discrete curvature, that is, flat discrete connections. We
    also introduce a semidirect product between (some) local Lie
    groupoids and prove that there is a correspondence between
    semidirect product extensions and right splittings of the DAS in
    \lLgpdC.
  \end{abstract}

  
\thanks{This research was partially supported by grants from the
  Universidad Nacional de Cuyo (grants 06/C567 and 06/C496),
  Universidad Nacional de La Plata, and CONICET}

\bibliographystyle{amsalpha}

\maketitle

\tableofcontents

%%%%%%%%%%%%%%%%%%%%%%%%%%%%%%%%%%%%%%%%%%%%%%%%%%%%%%%%%%

\section{Introduction}
\label{sec:introduction}

Connections on bundles are useful tools in Geometry. They can be used
to probe the intrinsic structures of the underlying bundle space; they
can also be used to obtain convenient global expressions of local
computations. Even in Physics, connections on vector and principal
bundles are essential ingredients of modern Gauge Theory
(see~\cite{ar:daniel_viallet-the_geometrical_setting_of_gauge)theories_of_the_yang_milss_type}). A
connection on the principal $\SG$-bundle
$\PBmap:\PBtotal\rightarrow \PBbase$ can be defined by giving a
(smooth and $\SG$-equivariant) $\jgg$-valued $1$-form $\CC$ on
$T\PBtotal$, where $\jgg$ is the Lie algebra of $\SG$; $\CC$ is known
as the \jdef{connection form}. Alternatively, a connection on $\PBmap$
can be defined by a \jdef{horizontal lift}
$h:\PBmap^*(T(\PBbase))\rightarrow T\PBtotal$, satisfying appropriate
conditions (see~\cite{bo:kobayashi_nomizu-foundations-v1}).

Lagrangian mechanical systems on a configuration manifold $\PBtotal$
are dynamical systems on the tangent bundle $T\PBtotal$, whose
trajectories are critical points of a functional determined by a
Lagrangian function $L:T\PBtotal\rightarrow\R$. Approximate
trajectories of these systems ---useful to approximate the evolution
of realistic systems--- can be obtained by appropriate discretizations
of the problem via a discrete-time dynamical system on
$\PBtotal\times \PBtotal$, whose trajectories are the critical points
of a functional determined by a discrete Lagrangian
$L_d:\PBtotal\times \PBtotal\rightarrow \R$. It seems natural to view
$\PBtotal\times \PBtotal$ as a discrete version of $T\PBtotal$:
indeed, while elements $v_q\in T_q\PBtotal$ can be seen as velocity
vectors at the point $q$, a discrete (time) version of this picture
replaces $v_q$ with the pair $(q,q')\in \PBtotal\times \PBtotal$, seen
as a short trajectory from $q$ to $q'$ in the general direction of
$v_q$
(see~\cite{ar:marsden_west-discrete_mechanics_and_variational_integrators}). Following
A. Weinstein's~\cite{ar:weinstein-lagrangian_mechanics_and_groupoids},
this idea has been expanded to deal with more general continuous
dynamical systems, not necessarily defined on $T\PBtotal$ but, rather,
on Lie algebroids and with discrete approximations given by dynamical
systems defined on Lie groupoids
(see~\cite{ar:cortes_deLeon_marrero_martinDeDiego_martinez-a_survey_of_lagrangian_mechanics_and_control_on_lie_algebroids_and_groupoids},~\cite{ar:marrero_martinDeDiego_martinez-the_local_description_of_discrete_mechanics}
and~\cite{ar:marrero_martinDeDiego_stern-symplectic_groupoids_anddiscrete_constrained_lagrangian_mechanics}).

When a Lagrangian mechanical system $(\PBtotal,L)$ is symmetric, that
is, for example, when the Lie group $\SG$ acts on $\PBtotal$ making
$\PBmap:\PBtotal\rightarrow \PBbase$ a principal $\SG$-bundle and $L$
a $\SG$-invariant function, a connection $\CC$ on $\PBmap$ can be used
to reduce the system as well as to reconstruct the original dynamics
in terms of that of the reduced system
(see~\cite{ar:cendra_marsden_ratiu-geometric_mechanics_lagrangian_reduction_and_nonholonomic_systems}). Developing
a similar strategy for discrete-time mechanical systems led M. Leok,
J. Marsden and A. Weinstein to propose a notion of discrete connection
on a principal bundle
in~\cite{ar:leok_marsden_weinstein-a_discrete_theory_of_connections_on_principal_bundles}
and~\cite{bo:leok-foundations_of_computational_geometric_mechanics};
this idea was later refined by some of us
in~\cite{ar:fernandez_zuccalli-a_geometric_approach_to_discrete_connections_on_principal_bundles}. Briefly,
a \jdef{discrete connection} on $\PBmap$ is determined by a map
$\DC:\PBtotal\times \PBtotal\rightarrow \SG$ (subject to some
conditions) called a \jdef{discrete connection form} or,
alternatively, by a \jdef{discrete horizontal lift}
$\HLd{}:\PBtotal\times(\PBbase)\rightarrow \PBtotal\times \PBtotal$
(again, subject to some conditions)\footnote{There are other ways to
  define a discrete connection on a principal bundle but these two are
  the only ones that are relevant for this work.}. Things are, indeed,
a bit more subtle as there may be no discrete connections for
nontrivial principal bundles: usually, for instance, $\DC$ is defined
only in an open subset $\jcalU\subset \PBtotal\times \PBtotal$ called
its \jdef{domain}. Below, in Section~\ref{sec:review_DC}, we give
complete definitions and the reader should
consult~\cite{ar:fernandez_zuccalli-a_geometric_approach_to_discrete_connections_on_principal_bundles}
for further details.

Associated to any principal $\SG$-bundle
$\PBmap:\PBtotal\rightarrow \PBbase$, M. Atiyah introduced
in~\cite{ar:atiyah-complex_analytic_connections_in_fibre_bundles} the
following short exact sequence of vector bundles over $\PBbase$, known
nowadays as the \jdef{Atiyah sequence} of $\PBmap$,
\begin{equation}
  \label{eq:atiyah_sequence-def-intro}
  \xymatrix{
    {0} \ar[r] & {\ti{\jgg}} \ar[r]^-{\phi_1} &
    {(T\PBtotal)/\SG} \ar[r]^-{\phi_2}  &
    {T(\PBbase)} \ar[r] & {0}
  }
\end{equation}
where $\ti{\jgg}:=(\PBtotal\times \jgg)/\SG$ is the \jdef{adjoint
  bundle}\footnote{When the Lie group $G$ acts on the manifold $X$ (on
  the left), we denote the action by $l^X:\SG\times X\rightarrow X$
  and the quotient map by $\pi^{X,\SG}:X\rightarrow X/\SG$.}  (for
$\SG$ acting on $\jgg$ by the Adjoint action),
$\phi_1(\pi^{\PBtotal\times\jgg,\SG}(q,\xi)) :=
\pi^{T\PBtotal,\SG}(\xi_\PBtotal(q))$ and
$\phi_2(\pi^{T\PBtotal,\SG}(v_q)) := T_q\PBmap(v_q)$.  He goes on to
define connections on $\PBmap$ as splittings
of~\eqref{eq:atiyah_sequence-def-intro} but, a well known argument
proves that this definition coincides with the standard definition of
principal connection on $\PBmap$. It is also possible to
view~\eqref{eq:atiyah_sequence-def-intro} as an exact sequence in the
category of Lie algebroids, as presented
in~\cite{bo:mackenzie-lie_groupoids_and_algebroids_in_differential_geometry}. In
this case, a right splitting $s_R$
of~\eqref{eq:atiyah_sequence-def-intro} (in the category of vector
bundles) may not be a morphism of Lie algebroids, with the obstruction
being the curvature of the connection associated to $s_R$.

%%%

A discrete analogue of~\eqref{eq:atiyah_sequence-def-intro} is
sketched
in~\cite{ar:leok_marsden_weinstein-a_discrete_theory_of_connections_on_principal_bundles},
although the precise algebraic or categorical context is not made
explicit, which is a nontrivial point because the categories possibly
involved are not abelian, so that the basic notions, like exactness,
require some careful interpretation. The main objective of this work
is to provide a precise definition of a \jdef{discrete Atiyah
  sequence} and prove that results analogous to the ones
described in the previous paragraph for the continuous setting are
valid in the discrete setting. For any principal $\SG$-bundle
$\PBmap:\PBtotal\rightarrow \PBbase$ we consider the sequence
\begin{equation}\label{eq:discrete_atiyah_sequence-def-intro}
  \xymatrix{
    {(\PBtotal\times \SG)/\SG} \ar[r]^-{F_1} &
    {(\PBtotal\times \PBtotal)/\SG} \ar[r]^-{F_2} &
    {(\PBbase)\times (\PBbase)}
  }
\end{equation}
for
\begin{equation}\label{eq:F_1_and_F_2-def}
  \begin{gathered}  
    F_1(\pi^{\PBtotal\times\SG,\SG}(q,g)) :=
    \pi^{\PBtotal\times \PBtotal,\SG}(q,l^\PBtotal_g(q)) \stext{ and }\\
    F_2(\pi^{\PBtotal\times \PBtotal,\SG}(q_0,q_1)) := (\PBmap(q_0),\PBmap(q_1))
  \end{gathered}
\end{equation}
as a sequence of smooth fiber bundles over $\PBbase$. The space
$\conjB := (\PBtotal\times \SG)/\SG$, obtained when $\SG$ acts on
$\PBtotal$ by $l^\PBtotal$ and on $\SG$ by conjugation, will be called
the \jdef{conjugate bundle} over $\PBbase$. We
call~\eqref{eq:discrete_atiyah_sequence-def-intro} the \jdef{discrete
  Atiyah sequence} (DAS) of $\PBmap$. In fact, for technical reasons
we have to work in the category $\FBS_{\PBbase}$ of fiber bundles with
a section (over $\PBbase$). As a first step we prove that the DAS is
an \jdef{extension} (Definition~\ref{def:extensionM_in_FBS}) of
$(\PBbase)\times (\PBbase)$ by $\conjB$ in the $\FBS_{\PBbase}$
category. Just as discrete connection forms may not be globally
defined, splittings of the DAS may not be defined globally but,
rather, semi-locally (Definition~\ref{def:semi-local_morphism}). Then
we prove that the (semi-local) left splittings of the DAS satisfying a
certain equivariance condition are in a bijective correspondence with
discrete connections on $\PBmap$, while (semi-local) right splittings
of the DAS are in a bijective correspondence with discrete horizontal
lifts on $\PBmap$. As discrete horizontal lifts and discrete
connections on $\PBmap$ are known to be equivalent notions, we have
proved the equivalence between some left splittings of the DAS, right
splittings of the DAS, and discrete connections on $\PBmap$. Notice
that in the non-abelian category $\FBS_{\PBtotal/\SG}$, in almost all
cases, the correspondence is between right splittings and only some
left splittings of the DAS. Another standard result, for instance in
the category of vector bundles, is that split exact sequences are
equivalent to direct sum sequences. We prove that the DAS is
semi-locally right (or left) split if and only if the sequence is
semi-locally equivalent to a certain fiber product extension of
$(\PBbase)\times (\PBbase)$ by $\conjB$. A complete ``map'' of the
relationship between the different objects considered so far is
provided by Theorem~\ref{thm:bijections_DAS_in_FBS}.

In the continuous case additional aspects of connections on $\PBmap$
are revealed when they are considered in the context of Lie
algebroids. Inspired by this fact and that a discrete connection form
$\DC:\PBtotal\times \PBtotal\rightarrow \SG$ on $\PBmap$ is a map
between Lie groupoids we ask if $\DC$ is a morphism of Lie
groupoids. It turns out that the interesting question is if $\DC$ is a
morphism from the pair groupoid $\PBtotal\times \PBtotal$ into the Lie
groupoid associated to the \jdef{opposite Lie group} $\SG^{op}$. In
general $\DC$ need not be such a morphism and we use the obstruction
to define the \jdef{curvature} $\BD$ of $\DC$
(Definition~\ref{def:discrete_curvature}). Then, naturally,
$\DC:\PBtotal\times \PBtotal\rightarrow \SG^{op}$ is a morphism of Lie
groupoids if and only if $\BD$ is trivial. Discrete connections with
trivial curvature are called \jdef{flat}.

As we mentioned before, in the continuous case, the
sequence~\eqref{eq:atiyah_sequence-def-intro} can be seen as a
sequence in the category of Lie algebroids. It was noted
in~\cite{ar:leok_marsden_weinstein-a_discrete_theory_of_connections_on_principal_bundles}
that~\eqref{eq:discrete_atiyah_sequence-def-intro} can be seen as a
sequence in the category of Lie groupoids $\LgpdC$, a fact that we
review here, and that has also been considered elsewhere in the
literature
(\cite{ar:androulidakis-classification_of_extensions_of_principal_bundles_and_transitive_lie_groupoids_with_prescribed_kernel_and_cokernel}
and~\cite{ar:casimiro_rodrigo-reduction_of_forward_operators_in_principal_G_bundles}). Our
interest is to explore the relationship between the (right) splittings
of~\eqref{eq:discrete_atiyah_sequence-def-intro} in the $\LgpdC$
category and discrete connections on $\PBmap$.  Still, a little twist
is required because discrete connections are defined in open subsets
$\jcalU\subset \PBtotal\times \PBtotal$ that are usually proper, so
that the right splittings that they induce are defined on open subsets
$(\PBmap\times \PBmap)(\jcalU)\subset (\PBbase)\times (\PBbase)$,
hence, they cannot be morphisms in $\LgpdC$ as these are globally
defined. So we move from the $\LgpdC$ category to the category of
\jdef{local Lie groupoids}, \lLgpdC.  We prove
that~\eqref{eq:discrete_atiyah_sequence-def-intro} is an extension in
the $\lLgpdC$ category
(Definition~\ref{def:extensionM-lLGpd}). Discrete connections on
$\PBmap$ induce (semi-local) right splittings $s_R$
of~\eqref{eq:discrete_atiyah_sequence-def-intro} in the
$\FBS_{\PBbase}$ category, but not necessarily in the
$\lLgpdC_{\PBbase}$ category; those $s_R$ that are morphisms in
$\lLgpdC_{\PBbase}$ and, so, split the DAS in the $\lLgpdC_{\PBbase}$
category are precisely those coming from \emph{flat} discrete
connections on $\PBmap$. In fact, we prove in
Proposition~\ref{prop:bijection_right_splittings_DAS_flat_DC-lLgpd}
that there is a bijective correspondence between the right splittings
of the DAS in the $\lLgpdC_{\PBbase}$ category and the flat discrete
connections on $\PBmap$. We briefly comment on how these
correspondences could be used to relate flat discrete connections to
flat connections on $\PBmap$ (Remark~\ref{rem:CC_and_DC}).

Inspired by the structure introduced
in~\cite{ar:metere_montoli-semidirect_products_of_internal_groupoids},
we consider the \jdef{semidirect product} of a (totally intransitive,
see Definition~\ref{def:local_lie_groupoid}) Lie groupoid by a local
Lie groupoid, that is a local Lie groupoid and can be used to
construct a certain semidirect product extension in
$\lLgpdC$. Theorem~\ref{prop:equiv_split_extensions_and_semi_direct_products}
proves that, for a given extension in the $\lLgpdC$ category, there is
a bijective correspondence between its right splittings and the
isomorphisms with those semidirect product extensions. As a
consequence, the DAS is right split in the $\lLgpdC$ category if and only
if it is isomorphic to a semidirect product extension, thus giving yet
another characterization of the flat discrete connections on $\PBmap$
(Corollary~\ref{cor:flat_DC_right)splittings_and_semidirect_product}).

The plan for the paper is as follows. In
Section~\ref{sec:some_background} we review some background material
regarding the continuous and the discrete Atiyah sequence, we also
revisit some of the basic notions about discrete connections on a
principal bundle. In Section~\ref{sec:DAS_in_FBS} we study the DAS of
a principal bundle in the category $\FBS$ of fiber bundles with a
section, a category that we introduce in
Section~\ref{sec:category_FBS}. The various bijective correspondences
proved in this Section are summarized in
Theorem~\ref{thm:bijections_DAS_in_FBS}. Lie groupoids and their
morphisms are reviewed in Section~\ref{sec:curvature_of_DC}, where we
define the \jdef{curvature} of a discrete connection on a principal
bundle. In Section~\ref{sec:DAS_in_lLgpd} we study the DAS in the
category of local Lie groupoids. The semidirect product of (some)
local Lie groupoids is introduced in
Section~\ref{sec:semidirect_product-lLGpd}, where its relationship to
split extensions in the $\lLgpdC$ category is established.

\emph{Notation:} throughout the paper the left action of a group $\SG$
on the set $X$ is denoted by $l^X$ and the corresponding quotient map
by $\pi^{X,\SG}:X\rightarrow X/\SG$. In general, we use
$\PBmap:\PBtotal\rightarrow \PBbase$ for a smooth principal
$\SG$-bundle induced by the (left) $\SG$-action $l^\PBtotal$ on
$\PBtotal$; we usually denote this bundle by $\PBmap$ or
$\PBtotal$. In addition to $l^{\PBtotal}$ we will consider some other
(left) $\SG$-actions: the lifted action on $T\PBtotal$, the diagonal
action on $\PBtotal\times\PBtotal$, the action on the second variable
of $\PBtotal\times\PBtotal$ and a ``diagonal'' action on
$\PBtotal\times\SG$:
\begin{equation}\label{eq:basic_actions-def}
  \begin{gathered}
    l^{T\PBtotal}_g(v_q) := T_ql^\PBtotal_g(v_q),\quad
    l^{\PBtotal\times\PBtotal}_g(q,q'):=(l^\PBtotal_g(q),l^\PBtotal_g(q')),\\
    l^{\PBtotal\times\PBtotal_2}_g(q,q') := (q,l^\PBtotal_g(q')),
    \stext{ and }
    l^{\PBtotal\times\SG}_g(q,h):=(l^{\PBtotal}_g(q),ghg^{-1}).
  \end{gathered}
\end{equation}
Given $f:X\rightarrow Y$ and subsets $X'\subset X$ and $Y'\subset Y$
such that $f(X)\subset Y'$, the maps $f|_{X'}:X'\rightarrow Y$ and
$f|^{Y'}:X\rightarrow Y'$ are the restriction and co-restriction of
$f$. Last, if $\phi_j:X_j\rightarrow Y$ are maps ($j=1,2$) then
$X_1\FP{\phi_1}{\phi_2}X_2:=\{(x_1,x_2)\in X_1\times X_2: \phi_1(x_1)
= \phi_2(x_2)\}$ is the \jdef{fiber product} of $X_1$ and $X_2$.

%%%%%%%%%%%%%%%%%%%%%%%%%%%%%%%%%%%%%%%%%%%%%%%%%%%%%%%%%%

\section{Some Background}
\label{sec:some_background}

In this section we review some of the main objects we work with in the
paper: the Atiyah sequence associated to a principal bundle and the
discrete connections of those spaces. We also introduce the discrete
version of the Atiyah sequence.

%%%%%%%%%%%%%%%%%%%%%%%%%%%%%%

\subsection{The Atiyah sequence}
\label{sec:the_atiyah_sequece}

The Atiyah sequence was introduced by M. Atiyah
in~\cite{ar:atiyah-complex_analytic_connections_in_fibre_bundles} to
study holomorphic connections on (holomorphic) principal bundles. The
construction has been successfully applied to smooth objects too, as
we review next. Later, in
Section~\ref{sec:the_atiyah_sequece-discrete_case}, we present the
discrete analogue (in the smooth setting).

\subsubsection{Continuous case}
\label{sec:the_atiyah_sequece-continuous_case}

Given the principal $\SG$-bundle $\PBmap:\PBtotal\rightarrow \PBbase$,
we have the following sequence of vector bundles over $\PBtotal$:
\begin{equation}\label{eq:pre_atiyah_sequence-def}
  \xymatrix{
    {0} \ar[r] & {\PBtotal\times\jgsg}  \ar[r]^-{\widehat{\phi_1}} &
    {T\PBtotal} \ar[r]^-{\widehat{\phi_2}} & {\PBmap^*T(\PBbase)} \ar[r] & {0}  
  }
\end{equation}
where $\jgsg := \lie{\SG}$,
$\widehat{\phi_1}(q,\xi) := \xi_\PBtotal(q) =
\frac{d}{dt}|_{t=0}l^\PBtotal_{\exp(t\xi)}(q)$ and
$\widehat{\phi_2}(v_q) := (q,T_q\PBmap(v_q))$. It is easy to see
that~\eqref{eq:pre_atiyah_sequence-def} is, indeed, an exact sequence
of vector bundles over $\PBtotal$.  If we now consider the
$\SG$-actions
$l^{\PBtotal\times \jgsg}_g(q,\xi) := (l^\PBtotal_g(q),Ad_g(\xi))$,
$l^{T\PBtotal}_g(v_q) = T_ql^\PBtotal_g(v_q)$ and
$l^{\PBmap^*T(\PBbase)}_g(q,r_{\PBmap(q)}) :=
(l^\PBtotal_g(q),r_{\PBmap(q)}) =
(l^\PBtotal_g(q),r_{\PBmap(l^\PBtotal_g(q))})$, both
$\widehat{\phi_1}$ and $\widehat{\phi_2}$ are $\SG$-equivariant, thus
we get a new sequence of vector bundles over $\PBbase$:
\begin{equation*}
  \xymatrix{
    {0} \ar[r] & {(\PBtotal\times\jgsg)/\SG} \ar[r] &{(T\PBtotal)/\SG} \ar[r] &
    {(\PBmap^*T(\PBbase))/\SG} \ar[r] & {0}
  }
\end{equation*}
The vector bundle $\ti{\jgsg} := (\PBtotal\times \jgsg)/\SG$ is the
\jdef{adjoint bundle}. Also, $(\PBmap^*T(\PBbase))/\SG$ and
$T(\PBbase)$ are isomorphic as vector bundles over $\PBbase$ via
$\pi^{\PBmap^*T(\PBbase),\SG}(q,r_{\PBmap(q)}) \mapsto r_{\PBmap(q)}$. Thus, the
previous sequence is isomorphic to the sequence
\begin{equation}
  \label{eq:atiyah_sequence-def}
  \xymatrix{
    {0} \ar[r] & {\ti{\jgg}} \ar[r]^-{\phi_1} &{(T\PBtotal)/\SG} \ar[r]^-{\phi_2}  &
  {T(\PBbase)} \ar[r] & {0}
}
\end{equation}
of vector bundles over $\PBbase$ that is exact and is called the
\jdef{Atiyah sequence} (AS) of $\PBmap$. Explicitly,
$\phi_1(\pi^{\PBtotal\times \jgsg,\SG}(q,\xi)) :=
\pi^{T\PBtotal,\SG}(\xi_\PBtotal(q))$ and
$\phi_2(\pi^{T\PBtotal,\SG}(v_q)) := T_q\PBmap(v_q)$.

The AS is closely related to connections on $\PBmap$. Indeed, Atiyah
defines a connection as a splitting of~\eqref{eq:atiyah_sequence-def}
(Definition on page 188
of~\cite{ar:atiyah-complex_analytic_connections_in_fibre_bundles}) and
so does Mackenzie (Definition 4.1 in Appendix A
of~\cite{bo:mackenzie-lie_groupoids_and_algebroids_in_differential_geometry}).
A lift of the image of a right splitting
of~\eqref{eq:atiyah_sequence-def} to $T\PBtotal$ corresponds precisely
to the horizontal distribution of such a connection, while its
connection form is equivalent to a lift of a left splitting
of~\eqref{eq:atiyah_sequence-def} (see pp. 292-3
in~\cite{bo:mackenzie-lie_groupoids_and_algebroids_in_differential_geometry}).

As stated in the Introduction, the AS~\eqref{eq:atiyah_sequence-def}
is an exact sequence in the category of Lie algebroids. Not all the
splittings described in the previous paragraph are morphisms in this
category. Only those that arise from flat connections correspond to
right splittings of~\eqref{eq:atiyah_sequence-def} in the category of
Lie algebroids.

For more details on these subjects good references are the Appendix A
of~\cite{bo:mackenzie-lie_groupoids_and_algebroids_in_differential_geometry},
and~\cite{ar:grabowski_kotov_poncin-geometric_structures_encoded_in_the_lie_structure_of_an_atiyah_algebroid}.

%%%%%%%%%%%%%%%

\subsubsection{Discrete case}
\label{sec:the_atiyah_sequece-discrete_case}

As we mentioned in the Introduction, the discrete analogue of
$T\PBtotal$ that we consider is the space $\PBtotal\times
\PBtotal$. Following ideas
of~\cite{ar:leok_marsden_weinstein-a_discrete_theory_of_connections_on_principal_bundles},
next we introduce a discrete analogue
of~\eqref{eq:atiyah_sequence-def}. A difference to notice is that
while the objects that appear in~\eqref{eq:atiyah_sequence-def} are
vector bundles, their discrete analogues are fiber bundles, so that no
linear structure is canonically available in the fibers.

\begin{definition}
  A \jdef{fiber bundle} $(E,\phi,M,S)$ consists of manifolds $E, M, S$
  and a smooth map $\phi:E\rightarrow M$ such that, for each $m\in M$,
  there is an open subset $U\subset M$ with $m\in U$ and such that
  $E|_U:=\phi^{-1}(U)$ is diffeomorphic to $U\times S$ by a fiber
  preserving diffeomorphism so that the following diagram is
  commutative.
  \begin{equation*}
    \xymatrix{
      {E|_U} \ar[r]^-{\Phi} \ar[d]_{\phi} & {U\times S}\ar[dl]^{p_1} \\
      {U} & {}
    }
  \end{equation*}
  If $(E_j,\phi_j,M_j,S_j)$ (for $j=1,2$) are fiber bundles, a smooth
  map $F:E_1\rightarrow E_2$ is a \jdef{bundle map} over the smooth
  function $f:M_1\rightarrow M_2$ if the following diagram commutes.
  \begin{equation*}
    \xymatrix{{E_1} \ar[r]^F \ar[d]_{\phi_1} & {E_2} \ar[d]^{\phi_2} \\
      {M_1} \ar[r]_{f} & {M_2}
    }
  \end{equation*}
  A \jdef{section} of the fiber bundle $(E,\phi,M,S)$ is a smooth
  function $\sigma: M\rightarrow E$ such that $\phi\circ \sigma =
  id_M$. The set of all sections of $E$ is denoted by $\Gamma(E)$.
\end{definition}

We notice that, by the local product structure, the projection map
$\phi$ of a fiber bundle is a submersion. It is customary to refer to
a fiber bundle $(E,\phi,M,S)$ just as $E$ or $\phi$, something that we
do in what follows unless it is necessary to provide all the data.

Inspired by~\eqref{eq:pre_atiyah_sequence-def}, given a principal
$\SG$-bundle $\PBmap:\PBtotal\rightarrow \PBbase$ we consider the following
sequence of fiber bundles over $\PBtotal$ (with the projection onto the first
variable as bundle projection in the three cases)
\begin{equation}
  \label{eq:discrete_pre_atiyah_sequence-def}
  \begin{gathered}
    \xymatrix{ {\PBtotal\times\SG} \ar[r]^-{\widehat{F_1}} & {\PBtotal\times \PBtotal}
      \ar[r]^-{\widehat{F_2}} &
      {\PBtotal\FP{\PBmap}{p_1}((\PBbase)\times (\PBbase))}}\\
    \text{where }\quad \widehat{F_1}(q,g):=(q,l^\PBtotal_g(q)) \stext{ and }
    \widehat{F_2}(q_0,q_1) := (q_0,(\PBmap(q_0),\PBmap(q_1))).
  \end{gathered}
\end{equation}
Clearly $\widehat{F_1}$ is injective and $\widehat{F_2}$ is
onto;~\eqref{eq:discrete_pre_atiyah_sequence-def} is a sequence of
fiber bundles and maps over $\PBtotal$. If we consider the
$\SG$-actions $l^{\PBtotal\times\PBtotal}$ and $l^{\PBtotal\times\SG}$
defined in~\eqref{eq:basic_actions-def}, and
\begin{equation*}
  l^{\PBtotal\FP{\PBmap}{p_1}((\PBbase)\times
    (\PBbase))}_g(q_0,(\PBmap(q_0),\PBmap(q_1))) :=
  (l^\PBtotal_g(q_0),(\PBmap(q_0),\PBmap(q_1))),
\end{equation*}
both $\widehat{F_1}$ and $\widehat{F_2}$ are $\SG$-equivariant. Thus,
from~\eqref{eq:discrete_pre_atiyah_sequence-def} we obtain the
sequence
\begin{equation*}
  \xymatrix{
    {(\PBtotal\times\SG)/\SG} \ar[r] &
    {(\PBtotal\times \PBtotal)/\SG} \ar[r] &
    {(\PBtotal\FP{\PBmap}{p_1}((\PBbase)\times (\PBbase)))/\SG}
  }.
\end{equation*}
Notice that
$(\PBmap(q_0),\PBmap(q_1))\mapsto
\pi^{\PBtotal\FP{\PBmap}{p_1}((\PBbase)\times(\PBbase)),\SG}(q_0,(\PBmap(q_0),
\PBmap(q_1)))$ is an isomorphism of fiber bundles over $\PBbase$, from
$(\PBbase)\times (\PBbase)$ onto
$(\PBtotal\FP{\PBmap}{p_1}((\PBbase)\times (\PBbase)))/\SG$. Thus, the previous
sequence produces the following sequence of fiber bundles over $\PBbase$
\begin{equation} \label{eq:discrete_atiyah_sequence-def}
  \begin{gathered}
    \xymatrix{ {\conjB} \ar[r]^(.3){F_1} & {(\PBtotal\times \PBtotal)/\SG}
      \ar[r]^(.43){F_2} & {(\PBbase)\times (\PBbase)}
    } \quad \text{ for }\\
    F_1(\pi^{\PBtotal\times\SG,\SG}(q,g)) := \pi^{\PBtotal\times \PBtotal,\SG}(q,l^\PBtotal_g(q))
    \stext{ and } F_2(\pi^{\PBtotal\times \PBtotal,\SG}(q_0,q_1)) :=
    (\PBmap(q_0),\PBmap(q_1)),
  \end{gathered}
\end{equation}
where $\conjB:=(\PBtotal\times \SG)/\SG$ is the \jdef{conjugate
  bundle} ---also known as the \jdef{gauge bundle}--- associated to
$\PBtotal$ over $\PBbase$.  By analogy
with~\eqref{eq:atiyah_sequence-def} we call
sequence~\eqref{eq:discrete_atiyah_sequence-def} the \jdef{discrete
  Atiyah sequence} (DAS) of $\PBmap$.

While the Atiyah sequence~\eqref{eq:atiyah_sequence-def} is a sequence
of vector bundles or, even, of Lie algebroids over $\PBbase$, so that
notions as exactness or being split have a well defined meaning, the
DAS~\eqref{eq:discrete_atiyah_sequence-def} is a sequence of fiber
bundles or, as we will discuss later, of Lie groupoids: in both cases,
the categorical notions mentioned above are not available. The study
of the DAS in these two contexts and its relationship to discrete
connections and their curvature are the subject of this paper.

%%%%%%%%%%%%%%%%%%%%%%%%%%%%%%

\subsection{Discrete connections on principal fiber bundles}
\label{sec:review_DC}

As we mentioned in the Introduction, discrete connections on principal
bundles were introduced
in~\cite{ar:leok_marsden_weinstein-a_discrete_theory_of_connections_on_principal_bundles}
and refined
in~\cite{ar:fernandez_zuccalli-a_geometric_approach_to_discrete_connections_on_principal_bundles}. We
refer the reader to the latter for further details on the subject.
Let $l^\PBtotal$ be a $\SG$-action on $\PBtotal$ so that
$\PBmap:\PBtotal\rightarrow \PBbase$ is a smooth principal
$\SG$-bundle. On $\PBtotal\times \PBtotal$ we consider the two
$\SG$-actions $l^{\PBtotal\times \PBtotal}$ and
$l^{\PBtotal\times \PBtotal_2}$ defined
in~\eqref{eq:basic_actions-def}. As usual, for any manifold $X$,
$\Delta_X\subset X\times X$ is the diagonal submanifold and we define
the \jdef{discrete vertical submanifold}
$\VD:= l^{\PBtotal\times \PBtotal_2}_\SG(\Delta_\PBtotal) =
\{(q,l^{\PBtotal}_g(q)) \in \PBtotal\times \PBtotal: q\in \PBtotal
\text{ and } g\in \SG\}$.

\begin{definition}\label{def:D_type_subset}
  An open subset $\jcalU\subset \PBtotal\times \PBtotal$ is said to be
  of \jdef{$pD$-type} or \jdef{pre-$D$-type} if it is
  $l^{\PBtotal\times \PBtotal}$-invariant and $\VD\subset \jcalU$ (in
  particular, $\Delta_\PBtotal\subset \jcalU$). A $pD$-type subset
  $\jcalU$ is of \jdef{$D$-type} if it is $\SG\times
  \SG$-invariant. If $\jcalU\subset \PBtotal\times \PBtotal$ is of
  $D$-type (or $pD$-type) and is $\Z_2$-invariant for the standard
  action that interchanges components it is said to be of
  \jdef{symmetric $D$-type} (or \jdef{symmetric $pD$-type}).
\end{definition}

Since $\PBmap\times \PBmap$ and $id_\PBtotal\times \PBmap$ are open maps, if
$\jcalU\subset \PBtotal\times \PBtotal$ is open (in particular, when
$\jcalU$ is a $pD$-type subset),
\begin{equation}
  \label{eq:prime_and_2prime_subsets-def}
  \jcalU':=(id_\PBtotal\times \PBmap)(\jcalU)\subset \PBtotal\times (\PBbase)
  \stext{ and }
  \jcalU'':=(\PBmap\times \PBmap)(\jcalU)\subset (\PBbase)\times (\PBbase)
\end{equation}
are open subsets. We will use this $'$ and $''$ notation in what
follows.

\begin{definition}\label{def:pre-DC}
  Let $\jcalU\subset \PBtotal\times \PBtotal$ be of $pD$-type. A map
  $\DC:\jcalU\rightarrow \SG$ is a \jdef{pre-discrete connection form}
  on $\PBmap$ if, for all $g\in\SG$,
  $\DC\circ l^{\PBtotal\times \PBtotal}_g = l^\SG_g\circ \DC$ (for
  $l^\SG_g(g'):=gg'g^{-1}$) over $\jcalU$ and
  $\DC(q,l^\PBtotal_g(q))=g$ for all $q\in \PBtotal$. The set of all
  discrete connection forms defined on $\jcalU$ is denoted by
  $\Sigma_C'(\jcalU)$.
\end{definition}

\begin{definition}\label{def:DC}
  Let $\jcalU\subset \PBtotal\times \PBtotal$ be of $D$-type. A
  pre-discrete connection form on $\PBmap$, $\DC:\jcalU\rightarrow G$,
  is a \jdef{discrete connection form} if it satisfies
  \begin{equation}
    \label{eq:Ad_GxG}
    \DC(l^\PBtotal_{g_0}(q_0),l^\PBtotal_{g_1}(q_1)) = g_1 \DC(q_0,q_1) g_0^{-1},
  \end{equation}
  for all $(q_0,q_1)\in\jcalU$ and $g_0,g_1\in \SG$. The set of all
  discrete connection forms defined on $\jcalU$ is denoted by
  $\Sigma_C(\jcalU)$. Clearly,
  $\Sigma_C(\jcalU)\subset \Sigma_C'(\jcalU)$.
\end{definition}

Horizontal lifts are useful tools for working with connections. The
following definition introduces their discrete analogues.

\begin{definition}\label{def:HLd}
  Let $\jcalU\subset \PBtotal\times \PBtotal$ be of $D$-type. A smooth
  function $\HLd{}:\jcalU'\rightarrow \PBtotal\times \PBtotal$ is a
  \jdef{discrete horizontal lift} on $\PBmap$ if the following
  conditions hold.
  \begin{enumerate}
  \item \label{it:HL_properties-HL_is_smooth}
    $\HLd{}:\jcalU'\rightarrow \PBtotal\times \PBtotal$ is
    $G$-equivariant for the $G$-actions
    $l^{\PBtotal\times \PBtotal}$ and
    $l^{\PBtotal\times (\PBbase)}_g(q,r):=(l^{\PBtotal}_g(q),r)$.
  \item \label{it:HL_properties-HL_is_a_section} $\HLd{}$ is a section
    of
    $(id_\PBtotal\times \PBmap):\PBtotal\times \PBtotal\rightarrow
    \PBtotal\times (\PBbase)$ over $\jcalU'$, that is,
    $(id_\PBtotal\times \PBmap) \circ \HLd{} = id_{\jcalU'}$.
  \item \label{it:HL_properties-HL_normalization} For every $q\in \PBtotal$,
    $\HLd{}(q,\PBmap(q)) = (q,q)$.
  \end{enumerate}
  We denote the set of all discrete horizontal lifts defined on the
  set $\jcalU'$ by $\Sigma_H(\jcalU)$.
\end{definition}

Discrete connection forms and discrete horizontal lifts are related as
follows. First, recall the following construction: consider the fiber
product $\PBtotal\FP{\PBmap}{\PBmap}\PBtotal$ of $\PBmap$ with
itself. Let
$\kappa:\PBtotal\FP{\PBmap}{\PBmap}\PBtotal\rightarrow \SG$ be defined
by
\begin{equation}
  \label{eq:kappa-def}
  \kappa(q_0,q_1):=g \stext{ if and only if } l^\PBtotal_g(q_0)=q_1.
\end{equation}
It is easy to check that $\kappa$ is a smooth function. Let
$\jcalU\subset \PBtotal\times \PBtotal$ be a $D$-type subset. Given a
discrete connection form on $\PBmap$, $\DC:\jcalU\rightarrow G$, we
define
\begin{equation}\label{eq:HLd_from_DC}
  h_\DC:\jcalU'\rightarrow \PBtotal\times \PBtotal \stext{ by }
  h_\DC(q,r) := (q,l^\PBtotal_{\DC(q,q')^{-1}}(q')),
\end{equation}
for any $q'\in \PBmap^{-1}(\{r\})$. Conversely, given a discrete
horizontal lift on $\PBmap$,
$\HLd{}:\jcalU'\rightarrow \PBtotal\times \PBtotal$, we define
\begin{equation}\label{eq:DC_from_HLd}
  \DCp{\HLd{}}:\jcalU\rightarrow G \stext{ by }
  \DCp{\HLd{}}(q_0,q_1):=\kappa(p_2(\HLd{}(q_0,\PBmap(q_1))),q_1).
\end{equation}

\begin{theorem}\label{th:DC_forms_and_HLd_are_equivalent}
  The maps $h_\DC$ and $\DCp{\HLd{}}$ defined
  by~\eqref{eq:HLd_from_DC} and~\eqref{eq:DC_from_HLd} are a discrete
  horizontal lift and a discrete connection form on $\PBmap$
  respectively. Furthermore, the two constructions are inverses of
  each other.
\end{theorem}

\begin{proof}
  It follows from Theorems 3.6 and 4.6
  of~\cite{ar:fernandez_zuccalli-a_geometric_approach_to_discrete_connections_on_principal_bundles}. See
  Remark~\ref{rem:DC_from_manifold} below.
\end{proof}

As a consequence of Theorem~\ref{th:DC_forms_and_HLd_are_equivalent}
discrete connection forms and discrete horizontal lifts are
alternative descriptions of a single object, the \jdef{discrete
  connections} with domain $\jcalU$. In this spirit we speak of a
``discrete connection'' with domain $\jcalU$ as the object defined by
either one of these maps. We also derive from
Theorem~\ref{th:DC_forms_and_HLd_are_equivalent} that the map
$F_{CH}:\Sigma_C(\jcalU)\rightarrow \Sigma_H(\jcalU)$ that assigns
$\DC\mapsto h_\DC$ (given by~\eqref{eq:HLd_from_DC}) is a bijection;
we denote its inverse by
$F_{HC}:\Sigma_H(\jcalU)\rightarrow \Sigma_C(\jcalU)$.

\begin{remark}\label{rem:DC_from_manifold}
  Just as connections on principal bundles can be defined by a
  horizontal distribution, discrete principal connections can be
  defined by a \jdef{horizontal submanifold}
  $Hor_\DC\subset \PBtotal\times \PBtotal$. Indeed, this is the
  approach
  of~\cite{ar:leok_marsden_weinstein-a_discrete_theory_of_connections_on_principal_bundles}
  and~\cite{ar:fernandez_zuccalli-a_geometric_approach_to_discrete_connections_on_principal_bundles}. But,
  as shown in Theorem 3.6
  of~\cite{ar:fernandez_zuccalli-a_geometric_approach_to_discrete_connections_on_principal_bundles},
  giving a horizontal submanifold is equivalent to giving a discrete
  connection form $\DC$. The horizontal submanifold associated to
  $\DC$ is $Hor_\DC:=\DC^{-1}(\{e\})$.
\end{remark}

\begin{remark}
  Connections on principal bundles
  $\PBmap:\PBtotal\rightarrow \PBbase$ always exist and can be
  constructed, for instance, using a $\SG$-invariant Riemannian metric
  on $Q$. Theorem 5.2
  in~\cite{ar:fernandez_zuccalli-a_geometric_approach_to_discrete_connections_on_principal_bundles}
  proves that the same thing happens for discrete connections on
  $\PBmap$.
\end{remark}

\begin{remark}\label{rem:HL_is_horizontal}
  If $\DC:\jcalU\rightarrow G$ is a discrete connection form on
  $\PBmap$ and $h_\DC$ is the associated discrete horizontal
  lift~\eqref{eq:HLd_from_DC}, then $\DC(h_{\DC}(q,r)) = e$ for all
  $(q,r)\in \jcalU'$.
\end{remark}

%%%%%%%%%%%%%%%%%%%%%%%%%%%%%%%%%%%%%%%%%%%%%%%%%%%%%%%%%%

\section{The discrete Atiyah sequence in the category \FBS}
\label{sec:DAS_in_FBS}

%%%%%%%%%%%%%%%%%%%%%%%%

\subsection{The category of fiber bundles with a section}
\label{sec:category_FBS}

It is well known that the smooth fiber bundles and the smooth bundle
maps (together with the standard identity and composition of maps)
form a category that we denote by $\FB$.  Unfortunately $\FB$ is not
an abelian or an exact category, so there is no general notion of
exactness or, even, of kernel. In order to recover some of these
notions, we work in a slightly enriched version of the category $\FB$.
The category of \jdef{fiber bundles with a section}, \FBS, has objects
\begin{equation*}
  \ob_\FBS :=\{((E,\phi,M,S),\sigma):\text{ such that }
  (E,\phi,M,S) \in \ob_\FB \text{ and } \sigma \in \Gamma(E)\},
\end{equation*}
and morphisms
\begin{equation*}
  \begin{split}
    \hom_\FBS((E_1,\sigma_1),(E_2,\sigma_2)) :=& \{
    (F,f)\in\hom_\FB(E_1,E_2) \text{ such that } F\circ \sigma_1 =
    \sigma_2 \circ f\}.
  \end{split}
\end{equation*}
With the standard identity and composition of bundle maps, $\FBS$ is
indeed a category. For any smooth manifold $M$ we denote by
$\FBS_M\subset \FBS$ the subcategory consisting of fiber bundles over
$M$ and morphisms over $id_M$; since all such morphisms are of the
form $(F,id_M)$ we refer to them just by $F$ in what follows.

Next we introduce some examples of fiber bundles with a section.

\begin{example}\label{ex:elements_of_FBS_Q/G}
  We saw in Section~\ref{sec:the_atiyah_sequece-discrete_case} that
  the discrete Atiyah sequence~\eqref{eq:discrete_atiyah_sequence-def}
  is a sequence in the $\FB_{\PBbase}$ category. A trivial check shows
  that, with the maps
    \begin{equation*}
    \begin{split}
      \sigma_{\conjB}:\PBbase\rightarrow \conjB \stext{ where }&
      \sigma_{\conjB}(\PBmap(q)):=\pi^{\PBtotal\times \SG,\SG}(q,e),\\
      \sigma_{(\PBtotal\times \PBtotal)/\SG}:\PBbase\rightarrow
      (\PBtotal\times \PBtotal)/\SG \stext{ where }&
      \sigma_{(\PBtotal\times \PBtotal)/\SG}(\PBmap(q)) :=
      \pi^{\PBtotal\times \PBtotal,\SG}(q,q),\\
      \sigma_{\PBbase\times \PBbase}:\PBbase\rightarrow
      (\PBbase)\times (\PBbase) \stext{ where }& \sigma_{\PBbase\times
        \PBbase)}(\PBmap(q)):=(\PBmap(q),\PBmap(q)),
    \end{split}
  \end{equation*}
  the DAS becomes a sequence in the $\FBS_{\PBbase}$ category.
\end{example}

\begin{example}\label{ex:initial_and_final_object_in_FBS_M}
  For any smooth manifold $M$, $(M,id_M,M,\{0\})$ is a (trivial) fiber
  bundle over $M$. Clearly, $id_M$ is a section of that bundle. Hence,
  \begin{equation*}
    M^\dagger:=((M,id_M,M,\{0\}),id_M) \in \ob_{\FBS_M}.
  \end{equation*}
  For $((E,\phi,M,S),\sigma) \in \ob_{\FBS_M}$, define
  $0^E:M\rightarrow E$ by $0^E := \sigma$ and $0_E:E\rightarrow M$ by
  $0_E := \phi$.  It is easy to verify that
  $\hom_{\FBS_M}(M^\dagger,(E,\sigma)) =\{0^E\}$ and
  $\hom_{\FBS_M}((E,\sigma),M^\dagger) =\{0_E\}$. An object in a
  category that has exactly one morphism to all objects in the
  category is called \jdef{initial}\footnote{We use a few simple
    categorical
    notions;~\cite{bo:adamek_herrlich_strecker-abstract_and_concrete_categories}
    could be used for further context.}, while an object so that every
  object has exactly one morphism into it is called \jdef{terminal};
  also, an object that is both initial and terminal in a category is
  called a \jdef{zero} object. Thus, $M^\dagger$ is a zero object in
  $\FBS_M$.
\end{example}

\begin{lemma}\label{le:F_2_is_a_surjective_submersion}
  The map $F_2$ defined by~\eqref{eq:discrete_atiyah_sequence-def} is
  a a surjective submersion.
\end{lemma}

\begin{proof}
  It follows from
  $F_2\circ \pi^{\PBtotal\times \PBtotal,\SG} = \PBmap\times \PBmap$
  and that $\PBmap\times \PBmap$ is a surjective submersion.
\end{proof}

As $\FBS_M$ is not an abelian category, the standard categorical
notions of exact sequence and extension and their consequences are not
naturally available. Still, in $\FBS_M$, we consider the following
naive notion, inspired by one that is in common use in the category of
Lie groupoids (see Definition~\ref{def:extensionM-LGpd}).

\begin{definition}\label{def:extensionM_in_FBS}
  A sequence
  $(E_1,\sigma_1)\xrightarrow{\eta_1} (E_2,\sigma_2)
  \xrightarrow{\eta_2} (E_3,\sigma_3)$ in $\FBS_M$ is said to be an
  \jdef{extension} of $(E_3,\sigma_3)$ by $(E_1,\sigma_1)$ if $\eta_1$
  is an embedding, $\eta_2$ is onto, and the subsets
  $\im(\eta_1):=\eta_1(E_1)$ and
  $\ker(\eta_2):=\eta_2^{-1}(\sigma_3(M))$ are equal. An extension as
  above is said to be \jdef{right split} if $\eta_2$ has a right
  inverse and \jdef{left split} if $\eta_1$ has a left inverse.
\end{definition}

\begin{remark}\label{rem:im_and_ker_in_FBS_for_extensionM_in_FBS}
  If a sequence
  $(E_1,\sigma_1)\xrightarrow{\eta_1} (E_2,\sigma_2)
  \xrightarrow{\eta_2} (E_3,\sigma_3)$ in $\FBS_M$ is an extension, it
  is easy to see that $\im(\eta_1)$ is a fiber bundle over $M$ with
  the projection $\phi_2|_{\im(\eta_1)}$. Furthermore, as
  $\eta_1\in \hom_{\FBS_M}((E_1,\sigma_1),(E_2,\sigma_2))$,
  $((\im(\eta_1),\sigma_2|^{\im(\eta_1)}) \in\ob_{\FBS_M}$. As
  $\ker(\eta_2) = \im(\eta_1)$, it also follows that
  $(\ker(\eta_2),\sigma_2|^{\ker(\eta_2)}) \in\ob_{\FBS_M}$.
\end{remark}

\begin{example}\label{ex:DAS_is_extensionM_in_FBS}
  By Example~\ref{ex:elements_of_FBS_Q/G} we know
  that~\eqref{eq:discrete_atiyah_sequence-def} is a sequence in
  $\FBS_{\PBbase}$. By Lemma~\ref{le:F_2_is_a_surjective_submersion},
  $F_2$ is onto. $F_1$ is an embedding: as
  $\hat{F}_1:\PBtotal\times\SG\rightarrow \PBtotal\times \PBtotal$ defined
  in~\eqref{eq:discrete_pre_atiyah_sequence-def} is a closed
  immersion, it is an embedding onto $\VD$ that is, also,
  $\SG$-equivariant and it follows easily that the induced map in the
  quotient, $F_1$, is an embedding onto $\VD/\SG$.  Last, as
  \begin{equation*}
    \ker(F_2) = F_2^{-1}(\sigma_{\PBbase\times \PBbase}(\PBbase\times
    \PBbase)) = l^{\PBtotal\times \PBtotal_2}_\SG(\Delta_\PBtotal)/\SG = \VD/\SG = \im(F_1),
  \end{equation*}
  the discrete Atiyah sequence~\eqref{eq:discrete_atiyah_sequence-def}
  is an extension in $\FBS_{\PBbase}$.
\end{example}

\begin{remark}\label{rem:categorical_stuff_in_FBS}
  It may seem that the notions of kernel and extension are quite
  ad-hoc. Still, there are categorical notions that support these
  ideas. For example, in a category ${\mathcal C}$ with a zero object
  $0$, given a morphism $f\in\hom_{\mathcal C}(A_1,A_2)$ a pair
  $(K,j)$ where $j\in\hom_{\mathcal C}(K,A_1)$ is said to be a
  \jdef{categorical kernel} of $f$ if the sequence
  $0\xleftarrow{0_K} K \xrightarrow{j} A_1$ is a pullback of
  $0\xrightarrow{0^{A_2}} A_2 \xleftarrow{f} A_1$, where $0_K$ and
  $0^{A_2}$ are the corresponding (unique) initial and final
  morphisms. As we saw in
  Example~\ref{ex:initial_and_final_object_in_FBS_M}, $M^\dagger$ is a
  zero object in the $\FBS_M$ category. Given a submersion
  $F\in\hom_{\FBS_M}((E_1,\sigma_1),(E_2,\sigma_2))$ we can define
  $K:=F^{-1}(\sigma_2(M))\subset E_1$ and $j:K\rightarrow E_1$. It is
  not hard to see that when $\phi_1|_K:K\rightarrow M$ is a fiber
  bundle, $(K,j)$ is a categorical kernel in the previous sense. In
  addition, back in a category ${\mathcal C}$ as above, a sequence
  $A_1\xrightarrow{f_1} A_2 \xrightarrow{f_2} A_3$ in ${\mathcal C}$
  is said to be a \jdef{categorical extension} if $f_2$ is an
  epimorphism and $(A_1,f_1)$ is a categorical kernel of $f_2$. In the
  $\FBS_M$ category it can be proved that an extension
  (Definition~\ref{def:extensionM_in_FBS}) is always a categorical
  extension. In particular, the
  DAS~\eqref{eq:discrete_atiyah_sequence-def} is a categorical
  extension in $\FBS_M$. Finally, we notice that a categorical
  extension
  $(E_1,\sigma_1)\xrightarrow{\eta_1} (E_2,\sigma_2)
  \xrightarrow{\eta_2} (E_3,\sigma_3)$ in $\FBS_M$ where
  $\phi_2|_{\eta_2^{-1}(\sigma_3(M))}:\eta_2^{-1}(\sigma_3(M))\rightarrow
  M$ is a fiber bundle is always an extension
  (Definition~\ref{def:extensionM_in_FBS}); one case where this
  condition is true is when $\eta_2$ is a proper submersion.
\end{remark}

%%%%%%%%%%%%%%%%%%%%%%%%%%%%%%%

Sometimes we need to work with maps between objects of $\FBS$ that are
not globally defined. The following definition introduces a notion
that is adequate in those cases.

\begin{definition}\label{def:semi-local_morphism}
  Let $(E_j,\sigma_j)\in \ob_{\FBS_M}$ for $j=1,2$,
  $\jcalU\subset E_1$ be open subset and
  $\Phi:\jcalU\rightarrow E_2$ be a smooth map. $\Phi$ is called
  a \jdef{semi-local morphism}\footnote{The reason we call these maps
    semi-local rather than just local is that they are defined in an
    open subset containing ``the diagonal'' $\sigma_1(M)$ rather than
    on any open subset.} from $(E_1,\sigma_1)$ to $(E_2,\sigma_2)$ if
  \begin{enumerate}
  \item \label{it:semi-local_morphism-phi}
    $\phi_2\circ \Phi = \phi_1|_\jcalU$, where
    $\phi_j:E_j\rightarrow M$ are the fiber bundle projections, and
  \item \label{it:semi-local_morphism-sigma}
    $\sigma_1(M)\subset \jcalU$ and
    $\Phi\circ \sigma_1 = \sigma_2$.
  \end{enumerate}
  We say that the semi-local morphism $\Phi$ is a \jdef{semi-local
    isomorphism} if there is a semi-local morphism
  $\Psi:\mathcal{W}\rightarrow E_1$ with $\mathcal{W}\subset E_2$ such
  that $\Psi\circ \Phi = id_\jcalU$ and
  $\Phi\circ \Psi=id_\mathcal{W}$.
\end{definition}

\begin{lemma}\label{le:semi-local_morph_if_section}
  Let $(E_j,\sigma_j)\in \ob_{\FBS_M}$ for $j=1,2$,
  $\jcalU\subset E_1$ be an open subset,
  $\Phi:\jcalU\rightarrow E_2$ be a smooth map and
  $\eta\in\hom_{\FBS_M}((E_2,\sigma_2),(E_1,\sigma_1))$.
  \begin{enumerate}
  \item \label{it:semi-local_morph_if_section-right} If
    $\eta\circ \Phi = id_{\jcalU}$ and
    condition~\ref{it:semi-local_morphism-sigma} in
    Definition~\ref{def:semi-local_morphism} hold or
  \item \label{it:semi-local_morph_if_section-left} if
    $\Phi\circ \eta = id_{\eta^{-1}(\jcalU)}$,
    $\sigma_1(M)\subset \jcalU$ and
    condition~\ref{it:semi-local_morphism-phi} in
    Definition~\ref{def:semi-local_morphism} hold,
  \end{enumerate}
  then $\Phi$ is a semi-local morphism from $(E_1,\sigma_1)$ to
  $(E_2,\sigma_2)$.
\end{lemma}

\begin{proof}
  It is easily deduced from the definitions.
\end{proof}

\begin{definition}
  Let $(E_j,\sigma_j)\in \ob_{\FBS_M}$ for $j=1,2$,
  $\eta\in\hom_{\FBS_M}((E_2,\sigma_2),(E_1,\sigma_1))$ and
  $\Phi:\jcalU\rightarrow E_2$ be a semi-local morphism from
  $(E_1,\sigma_1)$ to $(E_2,\sigma_2)$. We say that $\Phi$ is a
  \jdef{semi-local left inverse} of $\eta$ if
  $\Phi\circ \eta|_{\eta^{-1}(\jcalU)} = id_{\eta^{-1}({\mathcal
      U})}$. Analogously, we say that $\Phi$ is a \jdef{semi-local
    right inverse} of $\eta$ if $\eta\circ \Phi = id_{\jcalU}$.
\end{definition}

\begin{definition}
  A \jdef{semi-local left splitting} of the sequence
  $(E_1,\sigma_1) \xrightarrow{\eta_1} (E,\sigma) \xrightarrow{\eta_2}
  (E_2,\sigma_2)$ in $\FBS_M$ with domain $\jcalU$ is a
  semi-local left inverse of $\eta_1$ with domain
  $\jcalU$. Similarly a \jdef{semi-local right splitting} of the
  same sequence with domain $\jcalU$ is a semi-local right
  inverse of $\eta_2$ with domain $\jcalU$.
\end{definition}

%%%%%%%%%%%%%%%%%

\subsection{Left splittings and discrete connections}
\label{sec:left_splittings_and_connections}

In~\cite{ar:leok_marsden_weinstein-a_discrete_theory_of_connections_on_principal_bundles},
Section 4.6, it is discussed how discrete connections of
$\PBmap:\PBtotal\rightarrow \PBbase$ and left splittings of the
discrete Atiyah sequence~\eqref{eq:discrete_atiyah_sequence-def} are
related. Here we want to revisit that analysis making explicit the
categorical context.

Let $\jcalU\subset \PBtotal\times \PBtotal$ be of $pD$-type. Define
$\Sigma_L'(\jcalU)$ as the set of all semi-local left splittings of
the DAS~\eqref{eq:discrete_atiyah_sequence-def} defined over
$\jcalU/\SG$ (for the $\SG$-action $l^{\PBtotal\times \PBtotal}$).

Given $A\in\Sigma_C'(\jcalU)$ we define
\begin{equation*}
  \label{eq:tis_L_from_A-def}
  \ti{s_L}:\jcalU\rightarrow \PBtotal\times \SG \stext{ by }
  \ti{s_L}(q_0,q_1) := (q_0,A(q_0,q_1))
\end{equation*}
that is clearly smooth and $\SG$-equivariant (for the actions
$l^{\PBtotal\times \PBtotal}$ an $l^{\PBtotal\times\SG}$ defined
in~\eqref{eq:basic_actions-def}), so that it defines a smooth map
\begin{equation}
  \label{eq:s_L_from_A-def}
  s_L:\jcalU/\SG\rightarrow \conjB \stext{ such that }
  \pi^{\PBtotal\times \SG,\SG} \circ \ti{s_L} = s_L\circ \pi^{\PBtotal\times \PBtotal,\SG}. 
\end{equation}
It is easy to check that $s_L\in \Sigma_L'(\jcalU)$. Then,
\begin{equation}\label{eq:F_C'L'-def}
  F_{C'L'}:\Sigma_C'(\jcalU)\rightarrow \Sigma_L'(\jcalU)
  \stext{ such that } F_{C'L'}(A) := s_L
\end{equation}
with $s_L$ as in~\eqref{eq:s_L_from_A-def} is well defined.

Next we introduce the auxiliary function
\begin{equation*}
  \kappa_2:\PBtotal\FP{\PBmap}{\check{p}_1} \conjB\rightarrow \SG
  \stext{ such that } \kappa_2(q,\pi^{\PBtotal\times \SG,\SG}(q',g')) :=
  l^\SG_{\kappa(q,q')^{-1}}(g'),
\end{equation*}
where $\kappa$ is defined in~\eqref{eq:kappa-def}.  It is easy to
check that
\begin{equation*}\label{eq:kappa_2-equivariance}
  \kappa_2(l^\PBtotal_g(q),\pi^{\PBtotal\times \SG,\SG}(q',g')) =
  g\kappa_2(q,\pi^{\PBtotal\times \SG,\SG}(q',g')) g^{-1}.
\end{equation*}

Given $s_L\in \Sigma_L'(\jcalU)$ we define
\begin{equation*}\label{eq:hat_s_L-def}
  \hat{s_L}:\jcalU\rightarrow \PBtotal\FP{\PBmap}{\check{p}_1} \conjB
  \stext{ by }
  \hat{s_L}(q_0,q_1):=(q_0,s_L(\pi^{\PBtotal\times \PBtotal,\SG}(q_0,q_1)))
\end{equation*}
and 
\begin{equation}\label{eq:DC_from_s_L}
  A_{s_L}:\jcalU\rightarrow \SG \stext{ by } A_{s_L}
  := \kappa_2\circ \hat{s_L}.
\end{equation}
Being compositions of smooth functions, $\hat{s_L}$ and $A_{s_L}$ are
smooth; it is easy to verify that $A_{s_L}\in
\Sigma_C'(\jcalU)$. Then,
\begin{equation}\label{eq:F_L'C'-def}
  F_{L'C'}:\Sigma_L'(\jcalU) \rightarrow \Sigma_C'(\jcalU)
  \stext{ such that } F_{L'C'}(s_L) := A_{s_L}
\end{equation}
with $A_{s_L}$ as in by~\eqref{eq:DC_from_s_L} is well defined.

\begin{proposition}\label{prop:F_C'L'_and_FL'C'_are_bijective}
  The functions $F_{C'L'}$ and $F_{L'C'}$ defined
  by~\eqref{eq:F_C'L'-def} and~\eqref{eq:F_L'C'-def} are mutually
  inverses.
\end{proposition}

\begin{proof}
  It is routine checking that both compositions equal the identities.
\end{proof}

So far, we have assumed that $\jcalU\subset \PBtotal\times \PBtotal$
is of $pD$-type. Suppose that $\jcalU$ is of $D$-type. The elements of
$\Sigma_C(\jcalU)$ are precisely the discrete connection forms with
domain $\jcalU$ (on the fixed principal $\SG$-bundle
$\PBmap:\PBtotal\rightarrow \PBbase$). We define
$\Sigma_L(\jcalU) := F_{C'L'}(\Sigma_C(\jcalU))$ and
$F_{CL}:\Sigma_C(\jcalU)\rightarrow \Sigma_L(\jcalU)$ as the
appropriate restriction and co-restriction of $F_{C'L'}$; as
$F_{C'L'}$ is a bijection, so is $F_{CL}$. Notice that if
$F_{LC}:\Sigma_L(\jcalU)\rightarrow \Sigma_C(\jcalU)$ is the
corresponding restriction and co-restriction of $F_{L'C'}$, then
$F_{LC}^{-1} = F_{CL}$.

\begin{remark}\label{rem:intrinsic_description_of_Sigma_C}
  For any $s_L\in \Sigma_L'(\jcalU)$ we define
  \begin{equation}\label{eq:tis_L-def}
    \ti{s_L}:\jcalU\rightarrow \PBtotal\times \SG \stext{ by }
    \ti{s_L}(q_0,q_1) :=
    (q_0,\kappa_2(q_0,s_L(\pi^{\PBtotal\times \PBtotal,\SG}(q_0,q_1))))
  \end{equation}
  that, being a composition of smooth functions, is smooth.  It is
  easy to verify that the diagram of manifolds and smooth maps
  \begin{equation*}
    \xymatrix{
      {\PBtotal\times \SG} \ar[d]_{\pi^{\PBtotal\times\SG,\SG}}&
      {\jcalU} \ar[d]^{\pi^{\PBtotal\times \PBtotal,\SG}} \ar[l]_-{\ti{s_L}}   \\
      {\conjB} & {\jcalU/\SG} \ar[l]^{s_L} 
    }
  \end{equation*}
  is commutative, so that $\ti{s_L}$ is a lift of $s_L$. When
  $\jcalU$ is of $D$-type, we can consider the restriction of the
  $\SG$-action $l^{\PBtotal\times \PBtotal_2}$ to $\jcalU$ and, also, the
  $\SG$-action $l^{\PBtotal\times\SG_{lm}}_g(q,h):= (q,gh)$. It is then easy
  to check that $s_L \in \Sigma_L(\jcalU)$ if and only if
  $\ti{s_L}$ is $\SG$-equivariant for those actions.
\end{remark}

\begin{remark}
  Even for trivial principal $\SG$-bundles $\PBmap:\PBtotal\rightarrow \PBbase$
  with abelian $\SG$ it can be seen that
  $\Sigma_L(\jcalU)\subsetneq \Sigma_L'(\jcalU)$.
\end{remark}

%%%%%%%%%%%%%%%%%%%%%%%%%

\subsection{Right splittings and horizontal liftings}
\label{sec:right_splittings_and_horizontal_liftings}

Let $\PBmap:\PBtotal\rightarrow \PBbase$ be a principal $\SG$-bundle,
$\jcalU\subset \PBtotal\times \PBtotal$ be of $D$-type, and $\jcalU'$,
$\jcalU''$ as in~\eqref{eq:prime_and_2prime_subsets-def}. In this
section we study the relationship between discrete horizontal lifts of
$\PBmap$ and semi-local right splittings of the discrete Atiyah
sequence~\eqref{eq:discrete_atiyah_sequence-def}.

Before we tackle the real issues, we introduce two auxiliary
functions. Define
\begin{equation*}
  \ti{\lambda}:\PBtotal\FP{\PBmap}{(\PBmap\circ p_1)}(\PBtotal\times \PBtotal)
  \rightarrow \PBtotal\times \PBtotal
  \stext{ by } 
  \ti{\lambda}(q,(q_0,q_1)) := l^{\PBtotal\times
    \PBtotal}_{\kappa(q_0,q)}(q_0,q_1).
\end{equation*}
Let $l^{\PBtotal^3}$ be the $\SG$-action on
$\PBtotal\times (\PBtotal\times \PBtotal)$ defined by
$l^{\PBtotal^3}_g(q,(q_0,q_1)):= (q,l^{\PBtotal\times
  \PBtotal}_g(q_0,q_1))$; then
$\PBtotal\FP{\PBmap}{(\PBmap\circ p_1)}(\PBtotal\times \PBtotal)$ is a
$\SG$-invariant closed submanifold, so $l^{\PBtotal^3}$ restricts to a
$\SG$-action on
$\PBtotal\FP{\PBmap}{(\PBmap\circ p_1)}(\PBtotal\times \PBtotal)$ that
we still denote by $l^{\PBtotal^3}$; it is easy to see that this
action is free and proper.  A straightforward computation shows that
$\ti{\lambda}$ is $\SG$-invariant and, then, it induces a smooth map
\begin{equation}\label{eq:lambda-def}
  \lambda:\PBtotal\FP{\PBmap}{\check{p}_1}(\PBtotal\times \PBtotal)/\SG \rightarrow \PBtotal\times
  \PBtotal \stext{ such that }
  \lambda(q,\pi^{\PBtotal\times \PBtotal,\SG}(q_0,q_1)) =
  \ti{\lambda}(q,(q_0,q_1)).
\end{equation}
In addition, as, for $(q,(q_0,q_1))\in \PBtotal^3$ and $g\in\SG$,
\begin{equation*}
  \ti{\lambda}(l^\PBtotal_g(q),(q_0,q_1)) = l^{\PBtotal\times
    \PBtotal}_g(\ti{\lambda}(q,(q_0,q_1))),
\end{equation*}
it is easy to see that $\lambda$ is $\SG$-equivariant for the actions
$l^{\PBtotal\times \PBtotal}$ and
\begin{equation*}
  l^{\PBtotal\times (\PBtotal\times \PBtotal)/\SG}_g(q,\pi^{\PBtotal\times \PBtotal,\SG}(q_0,q_1)) :=
  (l^\PBtotal_g(q),\pi^{\PBtotal\times \PBtotal,\SG}(q_0,q_1)).
\end{equation*}

Let $\Sigma_R(\jcalU)$ be the set of all semi-local right
splittings of~\eqref{eq:discrete_atiyah_sequence-def} defined over
$\jcalU''$. Then, for $s_R\in \Sigma_R(\jcalU)$, we define
\begin{equation}
  \label{eq:h_from_sR-def}
  h:\jcalU'\rightarrow \PBtotal\times \PBtotal \stext{ by } 
  h(q,r):=\lambda(q,s_R(\PBmap(q),r)).
\end{equation}
It is straightforward that $h\in \Sigma_H(\jcalU)$, so we define
the map 
\begin{equation}\label{eq:F_RH-def}
  F_{RH}:\Sigma_R(\jcalU)\rightarrow \Sigma_H(\jcalU) \stext{ by }
  F_{RH}(s_R):=h
\end{equation}
where $h$ is as in~\eqref{eq:h_from_sR-def}.

Given $h\in \Sigma_H(\jcalU)$ we define
\begin{equation*}
  \hat{s_R}:\jcalU'\rightarrow (\PBtotal\times \PBtotal)/\SG \stext{ by }
  \hat{s_R}(q,r):=\pi^{\PBtotal\times \PBtotal,\SG}(h(q,r)).
\end{equation*}
Being a composition of smooth functions, $\hat{s_R}$ is smooth and as,
for any $g\in\SG$,
\begin{equation*}
  \begin{split}
    \hat{s_R}(l^\PBtotal_g(q),r) =& \pi^{\PBtotal\times
      \PBtotal,\SG}(h(l^\PBtotal_g(q_0),\PBmap(q_1))) = \pi^{\PBtotal\times \PBtotal,\SG}(l^{\PBtotal\times
      \PBtotal}_g(h(q_0,\PBmap(q_1)))) \\=& \pi^{\PBtotal\times \PBtotal,\SG}(h(q_0,\PBmap(q_1))) =
    \hat{s_R}(q,r),
  \end{split}
\end{equation*}
there is a unique smooth function
\begin{equation}
  \label{eq:s_R_from_h-def}
  s_R:\jcalU''\rightarrow (\PBtotal\times \PBtotal)/\SG \stext{ such that }
  s_R(\PBmap(q),r) = \hat{s_R}(q,r).
\end{equation}
It is easy to check that $s_R\in \Sigma_R(\jcalU)$, so that the map
\begin{equation}
  \label{eq:F_HR-def}
  F_{HR}:\Sigma_H(\jcalU)\rightarrow \Sigma_R(\jcalU)
  \stext{ such that } F_{HR}(h):=s_R
\end{equation}
with $s_R$ as in~\eqref{eq:s_R_from_h-def} is well defined.

\begin{proposition}\label{prop:F_RH_and_F_HR_are_bijective}
  The functions $F_{RH}$ and $F_{HR}$ defined by~\eqref{eq:F_RH-def}
  and~\eqref{eq:F_HR-def} are mutually inverses.
\end{proposition}

\begin{proof}
  It is a routine matter of checking that both compositions equal the
  identities.
\end{proof}

\begin{remark}
  By Proposition~\ref{prop:F_RH_and_F_HR_are_bijective}, when $\jcalU$
  is of $D$-type, there is a bijective correspondence between
  semi-local right splittings
  of~\eqref{eq:discrete_atiyah_sequence-def} and discrete horizontal
  liftings on $\PBmap$ which, in turn, are equivalent to discrete
  connections on $\PBmap$. This equivalence is the discrete analogue
  of the equivalence between right splittings of the Atiyah
  sequence~\eqref{eq:atiyah_sequence-def} and connections on $\PBmap$
  as
  in~\cite{bo:mackenzie-lie_groupoids_and_algebroids_in_differential_geometry}. Notice
  that while right and left splittings
  of~\eqref{eq:atiyah_sequence-def} are equivalent ---and both are
  equivalent to connections--- in the discrete case an additional
  equivariance condition is required of a left splitting
  of~\eqref{eq:discrete_atiyah_sequence-def} to correspond to a right
  splitting ---or a discrete connection---, contrary to what is
  asserted
  in~\cite{ar:leok_marsden_weinstein-a_discrete_theory_of_connections_on_principal_bundles}.
\end{remark}

%%%%%%%%%%%%%%%%%%%%%%%%%

\subsection{Fiber product isomorphism associated with a discrete connection}
\label{sec:fiber_product_associated_with_DC-FBS}

Let $((E_j,\phi_j,M,S_j),\sigma_j)\in \ob_{\FBS_M}$ for $j=1,2$. As
$\phi_j$ are submersions, it is easy to see that,
$E_\times := E_1\FP{\phi_1}{\phi_2} E_2$ is a manifold and, in fact,
$E_\times$ is a pullback of
$E_1\xrightarrow{\phi_1} M \xleftarrow{\phi_2} E_2$ in the category of
smooth manifolds. If we let $\phi_\times :E_\times\rightarrow M$ be
such that $\phi_\times := \phi_1\circ p_1^r = \phi_2\circ p_2^r$, it
is easy to check that $(E_\times,\phi_\times,M,S_1\times S_2)$ is a
fiber bundle and that, furthermore, if
$\sigma_\times(m):=(\sigma_1(m),\sigma_2(m))$, then
$(E_\times,\sigma_\times) \in \ob_{\FBS_M}$. As
$E_\times\subset E_1\times E_2$ is an embedded submanifold, we can
define the following four smooth maps:
\begin{gather*}
  F^\times_1:E_1\rightarrow E_\times,\quad F^\times_2:E_\times
  \rightarrow E_2,\quad s^\times_1: E_\times \rightarrow E_1,\quad
  s^\times_2:E_2\rightarrow E_\times \quad \text{ by }\\
  F^\times_1(e_1) := (e_1,\sigma_2(\phi_1(e_1))),\quad
  F^\times_2(e_1,e_2):=e_2,\\ s^\times_1(e_1,e_2):=e_1,\quad
  s^\times_2(e_2) := (\sigma_1(\phi_2(e_2)),e_2).
\end{gather*}
It is easy to check that $F^\times_1,F^\times_2,s^\times_1,s^\times_2$
are morphisms in $\FBS_M$ and that the relations
\begin{equation}
  \label{eq:relations_in_product_bundle}
  s^\times_1\circ F^\times_1 = id_{E_1},\quad F^\times_2\circ s^\times_2 = id_{E_2},
  \quad F^\times_2\circ F^\times_1 = \sigma_2\circ \phi_1,\quad
  s^\times_1\circ s^\times_2 = \sigma_1\circ \phi_2
\end{equation}
are satisfied. The sequence
\begin{equation}\label{eq:fiber_product_sequence_in_Fbs_M}
  \xymatrix{{(E_1,\sigma_1)} \ar[r]^{F^\times_1} & {(E_\times,\sigma_\times)}
    \ar[r]^{F^\times_2} & {(E_2,\sigma_2)}
  }
\end{equation}
will be called the \jdef{fiber product sequence} of $(E_1,\sigma_1)$
and $(E_2,\sigma_2)$.  It can be verified that
$(E_1,\sigma_1) \xleftarrow{s_1^\times} (E_\times,\sigma_\times)
\xrightarrow{F_2^\times} (E_2,\sigma_2)$ is a pullback of
$(E_1,\sigma_1) \xrightarrow{\phi_1} (M^\dagger,id_M)
\xleftarrow{\phi_2} (E_2,\sigma_2)$ in $\FBS_M$.

Also, the fiber product sequence is an extension in $\FBS_M$. Indeed,
$F_2^\times$ is a surjective submersion and, as $F_1^\times$ is a
right inverse of the smooth map $s_1^\times$,
Lemma~\ref{le:section_of_smooth_map_is_embedding} proves that
$F_1^\times$ is an embedding. Last, a direct verification shows that
$\im(F_1^\times) = \ker(F_2^\times)$.

\begin{lemma}\label{le:section_of_smooth_map_is_embedding}
  Let $\phi:X\rightarrow Y$ be a smooth map, $V\subset Y$ an open
  subset and $\sigma:V\rightarrow X$ smooth such that
  $\phi\circ \sigma = id_V$. Then, $\sigma$ is an embedding.
\end{lemma}

\begin{proof}
  It follows immediately from $\phi\circ \sigma = id_V$ that $\sigma$
  is an injective immersion. It is easy to check that, for any
  $U\subset V$, we have $\sigma(U) = \phi^{-1}(U)\cap \sigma(V)$. In
  particular, when $U$ is open, $\sigma(U)$ is open in $\sigma(V)$, so
  that $\sigma|^{\sigma(V)}:V\rightarrow \sigma(V)$ is open and, in
  turn, a homeomorphism. All together, $\sigma$ is an embedding.
\end{proof}

Given an extension in $\FBS_M$
\begin{equation}\label{eq:extensionM-FBS}
  \xymatrix{ {(E_1,\sigma_1)} \ar[r]^{F_1} & {(E,\sigma)} \ar[r]^{F_2} &
    {(E_2,\sigma_2)}
  }
\end{equation}
we want to study different relationships between semi-local morphisms
from~\eqref{eq:extensionM-FBS} into the fiber product
sequence~\eqref{eq:fiber_product_sequence_in_Fbs_M} and semi-local
splittings of~\eqref{eq:extensionM-FBS}. Fix an open subset
$\jcalU \subset E$ containing $\sigma(M)$.  Let
$\Phi:\jcalU\rightarrow E_\times$ be a semi-local morphism that
makes 
\begin{equation}\label{eq:E_diagram_with_Phi_and_s_1^times-FBS}
  \xymatrix{
    {E_1} \ar@{<->}[d]_{id_{E_1}} \ar[r]^{F_1^\times} &
    {E_\times} \ar[r]^{F_2^\times} & {E_2} \ar@{<->}[d]^{id_{E_2}}\\
    {E_1} \ar[r]_{F_1} & {E} \ar[u]_{\Phi} \ar[r]_{F_2} & {E_2}
  }
\end{equation}
a commutative diagram in $\FBS_M$. We define
\begin{equation}
  \label{eq:extensionM_s_1_from_Phi-FBS-def}
  s_1:\jcalU\rightarrow E_1 \stext{ by }
  s_1:=s_1^\times \circ \Phi
\end{equation}
It is easy to check that $s_1$ is a semi-local left splitting
of~\eqref{eq:extensionM-FBS}.

Conversely, given a semi-local left splitting
$s_1:\jcalU\rightarrow E_1$ of~\eqref{eq:extensionM-FBS}, we
define
\begin{equation}
  \label{eq:extensionM_Phi_from_s_1-FBS-def}
  \Phi:\jcalU\rightarrow E_\times \stext{ by }
  \Phi:=(s_1, F_2|_\jcalU).
\end{equation}
A simple verification shows that $\Phi$ is a well defined semi-local
morphism from $E$ into $E_\times$ for whom
diagram~\eqref{eq:E_diagram_with_Phi_and_s_1^times-FBS} is
commutative.

\begin{proposition}\label{prop:extensionM_morphisms_and_splittings-FBS}
  For the extension~\eqref{eq:extensionM-FBS} in the category
  $\FBS_M$, the construction of the semi-local left splitting $s_1$
  from a semi-local morphism $\Phi$
  using~\eqref{eq:extensionM_s_1_from_Phi-FBS-def} and the
  construction of the semi-local morphism $\Phi$ from the semi-local
  left splitting $s_1$
  using~\eqref{eq:extensionM_Phi_from_s_1-FBS-def} are mutually
  inverse operations.
\end{proposition}

\begin{proof}
  Straightforward check.
\end{proof}

We apply the previous constructions to $M=\PBbase$,
$(E_1,\sigma_1)=(\conjB,\sigma_{\conjB})$ and
$(E_2,\sigma_2)=((\PBbase)\times(\PBbase),
\sigma_{(\PBbase)\times(\PBbase)})$. Thus, we have the sequence in
$\FBS_{\PBbase}$
\begin{equation*}
  \xymatrix{
    {\conjB} \ar@<0.5ex>[r]^-{F_1^\times} &
    {\conjB\FP{\check{p_1}}{p_1}((\PBbase)\times (\PBbase))}
    \ar@<0.5ex>[r]^-{F_2^\times} \ar@<0.5ex>[l]^-{s_1^\times}&
    {(\PBbase)\times (\PBbase) \ar@<0.5ex>[l]^-{s_2^\times}}
  }.
\end{equation*}

The map
$\Theta:\conjB \times (\PBbase)\rightarrow
\conjB\FP{\check{p_1}}{p_1}((\PBbase)\times (\PBbase))$ defined by
\begin{equation*}
  \Theta(\pi^{\PBtotal\times\SG,\SG}(q_0,g),\PBmap(q_1)) :=
  (\pi^{\PBtotal\times\SG,\SG}(q_0,g),(\PBmap(q_0),\PBmap(q_1)))
\end{equation*}
is, actually, an isomorphism of fiber bundles over $\PBbase$. If we
define
\begin{equation*}
  \sigma_{\conjB \times (\PBbase)}(\PBmap(q)) := (\pi^{\PBtotal\times
  \SG,\SG}(q,e),\PBmap(q)),
\end{equation*}
it is easy to verify that $\Theta$ is an
isomorphism in $\FBS_{\PBbase}$. The previous sequence in $\FBS_{\PBbase}$
becomes
\begin{equation}\label{eq:product_sequence_in_FBS-def}
  \xymatrix{{\conjB} \ar@<0.5ex>[r]^-{\hat{F}_1} &
    {\conjB\times (\PBbase)}
    \ar@<0.5ex>[r]^-{\hat{F}_2} \ar@<0.5ex>[l]^-{\hat{s}_1} &
    {(\PBbase)\times (\PBbase) \ar@<0.5ex>[l]^-{\hat{s}_2}}
  }
\end{equation}
where
\begin{gather*}
  \hat{F}_1(\pi^{\PBtotal\times \SG,\SG}(q,g)):=(\pi^{\PBtotal\times
    \SG,\SG}(q,g),\PBmap(q)),\quad\quad
  \hat{F}_2(\pi^{\PBtotal\times\SG,\SG}(q,g),r):= (\PBmap(q),r),\\
  \hat{s}_1(\pi^{\PBtotal\times\SG,\SG}(q,g),r):=\pi^{\PBtotal\times\SG,\SG}(q,g),
  \quad \quad \hat{s}_2(\PBmap(q),r):=(\pi^{\PBtotal\times
    \SG,\SG}(q,e),r).
\end{gather*}
Sequence~\eqref{eq:product_sequence_in_FBS-def} is the \jdef{fiber
  product sequence} of $(\conjB, \sigma_{\conjB})$ with
$((\PBbase)\times(\PBbase),\\\sigma_{(\PBbase)\times(\PBbase)})$.

Going back to the existence of a discrete connection on a principal
bundle, we have the following result.

\begin{proposition}\label{prop:Phi_DC_and_Psi_DC_are_isomorphism}
  Let $\jcalU\subset Q\times Q$ be of $D$-type and
  $\DC\in \Sigma_C(\jcalU)$ (on $\PBmap$).  Let
  $\mathcal{W}:= \{(q_0,g,r_1)\in \PBtotal\times\SG\times(\PBbase) :
  (q_0,r_1)\in\jcalU'\}$ and consider the $\SG$-actions
  $l^{\PBtotal\times \PBtotal}$ and
  $l^{\PBtotal\times \SG\times (\PBbase)}_g(q_0,w_0,r_1) :=
  (l^\PBtotal_g(q_0), l^\SG_g(w_0), r_1)$. Then, the maps
  $\Phi_\DC:\jcalU/\SG\rightarrow \mathcal{W}/\SG$ and
  $\Psi_\DC:\mathcal{W}/\SG\rightarrow \jcalU/\SG$ defined by
  \begin{equation}\label{eq:explicit_Phi_DC_and_Psi_DC}
    \begin{split}
      \Phi_\DC(\pi^{\PBtotal\times \PBtotal,\SG}(q_0,q_1)) :=&
      (\pi^{\PBtotal\times\SG,\SG}(q_0,\DC(q_0,q_1)),\PBmap(q_1)),\\
      \Psi_\DC(\pi^{\PBtotal\times\SG,\SG}(q_0,w_0),r_1) :=&
      \pi^{\PBtotal\times \PBtotal,\SG}(l^{\PBtotal\times
        \PBtotal_2}_{w_0}(h_\DC(q_0,r_1)))
    \end{split}
  \end{equation}
  are mutually inverse semi-local isomorphisms in $\FBS_{\PBbase}$. In
  addition, ignoring the smaller domains of $\Phi_\DC$ and $\Psi_\DC$,
  the following diagrams in $\FBS_{\PBbase}$ are commutative.
  \begin{equation*}
    \xymatrix{
      {\conjB} \ar[r]^-{\hat{F}_1} \ar@{<->}[d]_{id_{\conjB}} &
      {\conjB\times(\PBbase)} \ar[r]^-{\hat{F}_2} &
      {(\PBbase)\times (\PBbase)} \ar@{<->}[d]^{id_{(\PBbase)\times (\PBbase)}}\\
      {\conjB} \ar[r]_-{F_1} & {(\PBtotal\times \PBtotal)/\SG} \ar[r]_-{F_2}
      \ar[u]_{\Phi_\DC}&
      {(\PBbase)\times (\PBbase)}
    }
  \end{equation*}
  and
  \begin{equation*}
    \xymatrix{
      {\conjB} \ar[r]^-{\hat{F}_1} \ar@{<->}[d]_{id_{\conjB}} &
      {\conjB\times(\PBbase)} \ar[r]^-{\hat{F}_2} \ar[d]_{\Psi_\DC} &
      {(\PBbase)\times (\PBbase)} \ar@{<->}[d]^{id_{(\PBbase)\times (\PBbase)}}\\
      {\conjB} \ar[r]_-{F_1} & {(\PBtotal\times \PBtotal)/\SG} \ar[r]_-{F_2} &
      {(\PBbase)\times (\PBbase)}
    }
  \end{equation*}
\end{proposition}

\begin{proof}
  It follows from Proposition 4.19
  in~\cite{ar:fernandez_tori_zuccalli-lagrangian_reduction_of_discrete_mechanical_systems}
  and chasing the diagrams.
\end{proof}

\begin{corollary}
  Given a discrete connection with domain $\jcalU$ on the principal
  $\SG$-bundle $\PBmap:\PBtotal\rightarrow \PBbase$, the discrete
  Atiyah sequence~\eqref{eq:discrete_atiyah_sequence-def} is
  semi-locally equivalent to the fiber product
  sequence~\eqref{eq:product_sequence_in_FBS-def} in $\FBS_{\PBbase}$,
  in a canonical way.
\end{corollary}

In the next two sections we explore how, semi-local equivalences in
$\FBS_{\PBbase}$ between the
DAS~\eqref{eq:discrete_atiyah_sequence-def} and the fiber product
sequence~\eqref{eq:product_sequence_in_FBS-def} may, in turn, define
discrete connections on $\PBmap$.

%%%%%%%%%%%%%%%%%%%%%%%%

\subsection{Left splittings and fiber product decomposition}
\label{sec:left_splittings_and_fiber_product_decomposition}

In this section we explore the relation between having a semi-local
morphism $\Phi$ from the DAS into the fiber product sequence and
having a semi-local left splitting of the DAS.

Let $\jcalU\subset \PBtotal\times \PBtotal$ be of $pD$-type and let
$\Sigma_U'(\jcalU)$ be the set of all semi-local morphisms
$\Phi:\jcalU/\SG\rightarrow \conjB\times (\PBbase)$ that make the
following diagram in $\FBS_{\PBbase}$ commutative (wherever it makes
sense, because $\Phi$ is semi-local).
\begin{equation}\label{eq:diagram_vert_morph_phi}
  \xymatrix{
    {\conjB} \ar[r]^-{\hat{F}_1} \ar@{<->}[d]_{id_{\conjB}} &
    {\conjB\times(\PBbase)} \ar[r]^-{\hat{F}_2}&
    {(\PBbase)\times (\PBbase)} \ar@{<->}[d]^{id_{(\PBbase)\times (\PBbase)}}\\
    {\conjB} \ar[r]_-{F_1} & {(\PBtotal\times \PBtotal)/\SG} \ar[r]_-{F_2}
    \ar[u]_{\Phi}&
    {(\PBbase)\times (\PBbase)}
  }
\end{equation}

The constructions introduced in
Section~\ref{sec:fiber_product_associated_with_DC-FBS} can now be
specialized to the current context to provide bijections between
$\Sigma_U'(\jcalU)$ and $\Sigma_L'(\jcalU)$. Given
$\Phi\in \Sigma_U'(\jcalU)$, we define
\begin{equation*}
  \label{eq:DAS_sequence_s_1_from_Phi-FBS-def}
  s_L:\jcalU/\SG\rightarrow \conjB \stext{ by }
  s_L:=\hat{s}_1 \circ \Phi,
\end{equation*}
with $\hat{s}_1$ from~\eqref{eq:product_sequence_in_FBS-def}.  It was
proved in Section~\ref{sec:fiber_product_associated_with_DC-FBS} that
$s_L\in \Sigma_L'(\jcalU)$. Thus, we have a map
\begin{equation}\label{eq:F_U'L'-def}
  F_{U'L'}:\Sigma_U'(\jcalU)\rightarrow \Sigma_L'(\jcalU)
  \stext{ defined by } F_{U'L'}(\Phi):=s_L.
\end{equation}

In the same spirit, given $s_L\in \Sigma_L'(\jcalU)$ we define 
\begin{equation}
  \label{eq:DAS_sequence_Phi_from_s_1-FBS-def}
  \Phi:\jcalU/\SG\rightarrow \conjB\times(\PBbase)
  \stext{ by } \Phi:=(s_L, p_2\circ F_2|_{\jcalU/\SG}).
\end{equation}
It was proved in Section~\ref{sec:fiber_product_associated_with_DC-FBS} that
$\Phi\in \Sigma_U'(\jcalU)$, hence we can define a map
\begin{equation}\label{eq:F_L'U'-def}
  F_{L'U'}:\Sigma_L'(\jcalU)\rightarrow \Sigma_U'(\jcalU)
  \stext{ by } F_{L'U'}(s_L):=\Phi.
\end{equation}

The following result is a direct specialization of
Proposition~\ref{prop:extensionM_morphisms_and_splittings-FBS}.
\begin{corollary}\label{cor:F_L'U'_and_F_U'L'_are_bijective}
  The functions $F_{U'L'}$ and $F_{L'U'}$ defined
  by~\eqref{eq:F_U'L'-def} and~\eqref{eq:F_L'U'-def} are mutually
  inverses.
\end{corollary}

When $\jcalU$ is of $D$-type, the subset
$\Sigma_L(\jcalU)\subset \Sigma_L'(\jcalU)$ corresponds to
those semi-local left splittings that arise from discrete connections
on $\PBmap$ with domain $\jcalU$. Thus,
$\Sigma_U(\jcalU) := F_{L'U'}(\Sigma_L(\jcalU))$ is the set
of those semi-local morphisms in $\Sigma_U'(\jcalU)$ that arise
from a discrete connection on $\PBmap$ with domain $\jcalU$. We
define $F_{LU}:\Sigma_L(\jcalU)\rightarrow \Sigma_U(\jcalU)$
as the restriction and co-restriction of $F_{L'U'}$ to the
corresponding subsets. It is clear that $F_{LU}$ is still a bijection
and that, if
$F_{UL}:\Sigma_U(\jcalU)\rightarrow \Sigma_L(\jcalU)$ is the
corresponding restriction and co-restriction of $F_{U'L'}$, then
$F_{LU}^{-1} = F_{UL}$.

\begin{remark}\label{rem:intrinsic_description_of_Sigma_U}
  Assume that $\jcalU$ is of $D$-type. For any
  $\Phi\in \Sigma_U'(\jcalU)$, let $s_L:=F_{U'L'}(\Phi)$ and
  $\ti{s_L}$ be the lift defined by~\eqref{eq:tis_L-def}. Then,
  $s_L\in \Sigma_L(\jcalU)$ if and only if $\ti{s_L}$ is
  $\SG$-equivariant in the sense of
  Remark~\ref{rem:intrinsic_description_of_Sigma_C}. Define
  \begin{equation}\label{tiPhi-def}
    \ti{\Phi}:\jcalU \rightarrow (\PBtotal\times \SG)\times (\PBbase)
    \stext{ by } \ti{\Phi}(q_0,q_1):=(\ti{s_L}(q_0,q_1),\PBmap(q_1)).
  \end{equation}
  Clearly $\ti{\Phi}$ is smooth and it is easy to check that the diagram
  \begin{equation*}
    \xymatrix{
      {\PBtotal\times \PBtotal} \ar[r]^-{\ti{\Phi}} \ar[d]_{\pi^{\PBtotal\times \PBtotal,\SG}} &
      {(\PBtotal\times\SG)\times (\PBbase)} \ar[d]^{\pi^{\PBtotal\times\SG,\SG} \times id_{\PBbase}}\\
      {(\PBtotal\times \PBtotal)/\SG} \ar[r]_-{\Phi} & {\conjB\times (\PBbase)}
    }
  \end{equation*}
  in the category of manifolds and smooth maps is commutative. Hence,
  $\ti{\Phi}$ is a lift of $\Phi$. It can be checked that the
  condition of $\ti{s_L}$ being $\SG$-equivariant is equivalent to the
  $\SG$-equivariance of $\ti{\Phi}$ for the $\SG$-actions
  $l^{\PBtotal\times \PBtotal_2}$ and
  $l^{\PBtotal\times \SG_{lm}\times(\PBbase)}_g(q,h,r):=(q,gh,r)$. Thus,
  $\Phi\in \Sigma_U(\jcalU)$ if and only if its lift $\ti{\Phi}$
  is $\SG$-equivariant for these actions.
\end{remark}

\begin{proposition}\label{prop:Phi_from_DC_is_isomorphism}
  If $\jcalU$ is of $D$-type and $\Phi\in \Sigma_U(\jcalU)$,
  then $\Phi$ is a semi-local isomorphism in $\FBS_{\PBbase}$.
\end{proposition}

\begin{proof}
  Let $\DC := F_{LC}(F_{UL}(\Phi))$, so that
  $\Phi = F_{LU}(F_{CL}(\DC))$. Then, unraveling the definitions
  using~\eqref{eq:F_C'L'-def}
  and~\eqref{eq:DAS_sequence_Phi_from_s_1-FBS-def} we see that, for
  all $\pi^{\PBtotal\times \PBtotal,\SG}(q_0,q_1) \in \jcalU/\SG$,
  \begin{equation*}
    \Phi(\pi^{\PBtotal\times \PBtotal,\SG}(q_0,q_1)) = (\pi^{\PBtotal\times
      \SG,\SG}(q_0,\DC(q_0,q_1)),\PBmap(q_1)).
  \end{equation*}
  Comparison of this identity
  with~\eqref{eq:explicit_Phi_DC_and_Psi_DC} and
  Proposition~\ref{prop:Phi_DC_and_Psi_DC_are_isomorphism} lead to
  the result.
\end{proof}

\begin{remark}
  The primary use of discrete connections in the reduction and
  reconstruction of the dynamics of discrete mechanical systems is,
  precisely, because a discrete connection $\DC$ can be used to
  produce an isomorphism $\Phi_{\DC} := (F_{LU}\circ F_{CL})(\DC)$
  between $(Q\times Q)/\SG$ ---the natural space for the reduced
  system--- and $\conjB\times (Q/\SG)$ ---an easier to work with model
  (\cite{ar:leok_marsden_weinstein-a_discrete_theory_of_connections_on_principal_bundles}
  and
  \cite{ar:fernandez_tori_zuccalli-lagrangian_reduction_of_discrete_mechanical_systems}).
\end{remark}

%%%%%%%%%%%%%%%%%%%%%%%%%%%

\subsection{Right splittings and fiber product decomposition}
\label{sec:right_splittings_and_fiber_product_decomposition}

In this section we explore the relation between having a semi-local
morphism $\Psi$ from the fiber product sequence into the DAS and
having a semi-local right splitting of the DAS.

Let $\jcalU\subset \PBtotal\times \PBtotal$ be of $D$-type; we define
the open subset
$\mathcal{W}:=\hat{F}_2^{-1}(\jcalU'') \subset \conjB\times(\PBbase)$
(recall~\eqref{eq:prime_and_2prime_subsets-def}). Let
$\Sigma_D'(\jcalU)$ be the set of all semi-local morphisms
$\Psi:\mathcal{W}\rightarrow (\PBtotal\times \PBtotal)/\SG$ such that
the following diagram in $\FBS_{\PBbase}$ is commutative.
\begin{equation}\label{eq:semi_local_morphism_Psi-diagram}
  \xymatrix{
    {\conjB} \ar[r]^-{\hat{F}_1}
    \ar@{<->}[d]_{id_{\conjB}} &
    {\conjB\times(\PBbase)} \ar[r]^-{\hat{F}_2}
    \ar[d]_{\Psi} & {(\PBbase)\times (\PBbase)}
    \ar@{<->}[d]^{id_{(\PBbase)\times (\PBbase)}}\\
    {\conjB} \ar[r]_-{F_1} & {(\PBtotal\times \PBtotal)/\SG} \ar[r]_-{F_2} &
    {(\PBbase)\times (\PBbase)}
  }
\end{equation}

For $s_R\in \Sigma_R(\jcalU)$ we define
\begin{equation}\label{eq:Psi_from_s_r-def}
  \Psi:\mathcal{W}\rightarrow (\PBtotal\times \PBtotal)/\SG \stext{ by }
  \Psi(\pi^{\PBtotal\times \SG,\SG}(q,g),r) :=
  \pi^{\PBtotal\times \PBtotal,\SG}(l^{\PBtotal\times \PBtotal_2}_g(\lambda(q,s_R(\PBmap(q),r))))
\end{equation}
where $\lambda$ comes from~\eqref{eq:lambda-def}. A lengthy but
straightforward verification shows that
$\Psi\in \Sigma_D'(\jcalU)$. Then, we define
\begin{equation}
  \label{eq:F_RD'-def}
  F_{RD'}:\Sigma_R(\jcalU)\rightarrow \Sigma_D'(\jcalU)
  \stext{ so that } F_{RD'}(s_R) := \Psi,
\end{equation}
where $\Psi$ is defined by~\eqref{eq:Psi_from_s_r-def}.

Conversely, given $\Psi\in\Sigma_D'(\jcalU)$, we define
\begin{equation}
  \label{eq:s_R_from_Psi-def}
  s_R:\jcalU''\rightarrow (\PBtotal\times \PBtotal)/\SG \stext{ by }
  s_R:=\Psi\circ \hat{s}_2|_{\jcalU''},
\end{equation}
with $\hat{s}_2$ from~\eqref{eq:product_sequence_in_FBS-def}.  A
direct computation proves that $s_R\in \Sigma_R(\jcalU)$.  Thus,
\begin{equation}
  \label{eq:F_D'R-def}
  F_{D'R}:\Sigma_D'(\jcalU)\rightarrow \Sigma_R(\jcalU)
  \stext{ such that } F_{D'R}(\Psi):=s_R,
\end{equation}
where $s_R$ is given by~\eqref{eq:s_R_from_Psi-def}, is well defined.

\begin{proposition}\label{prop:F_RD''_is_injective}
  For $\jcalU\subset \PBtotal\times \PBtotal$ of $D$-type, the map $F_{D'R}$
  defined by~\eqref{eq:F_D'R-def} is a left inverse of $F_{RD'}$
  defined in~\eqref{eq:F_RD'-def}. Consequently, $F_{RD'}$ is
  one-to-one.
\end{proposition}

\begin{proof}
  Follows by evaluation and unraveling the definitions.
\end{proof}

\begin{remark}\label{rem:F_RD'_not_onto}
  With the previous definitions,
  $F_{RD'}:\Sigma_R(\jcalU)\rightarrow \Sigma_D'(\jcalU)$
  may not be onto. Indeed, it is possible to find counterexamples when
  $\PBmap:\PBtotal\rightarrow \PBbase$ is a trivial principal $\SG$-bundle with
  abelian $\SG$.
\end{remark}

Let
$\Sigma_D(\jcalU) := F_{RD'}(\Sigma_R(\jcalU))\subset
\Sigma_D'(\jcalU)$. As, by
Proposition~\ref{prop:F_RD''_is_injective}, $F_{RD'}$ is one-to-one,
its co-restriction to $\Sigma_D(\jcalU)$,
$F_{RD}:\Sigma_R(\jcalU)\rightarrow \Sigma_D(\jcalU)$, is a
bijection.

Let $\Psi\in \Sigma_D'(\jcalU)$ and define
$\hat{\mathcal{W}} := (\pi^{\PBtotal\times
  \SG,\SG}\times id_{\PBbase})^{-1}(\mathcal{W})\subset
(\PBtotal\times\SG)\times (\PBbase)$ for, as before,
$\mathcal{W}:=\hat{F}_2^{-1}(\jcalU'')$. In addition, define
\begin{equation}\label{eq:hatPsi-def}
  \hat{\Psi}:\hat{\mathcal{W}}\rightarrow \PBtotal\times \PBtotal \stext{ by }
  \hat{\Psi}(q,g,r) :=
  \lambda(q,\Psi(\pi^{\PBtotal\times\SG,\SG}(q,g),r))
\end{equation}
where $\lambda$ comes from~\eqref{eq:lambda-def}.

\begin{lemma}\label{le:hatPsi_lifts_Psi}
  For $\Psi\in \Sigma_D'(\jcalU)$ and $\hat{\Psi}$ as
  in~\eqref{eq:hatPsi-def}, $\hat{\Psi}$ is well defined and the
  following diagram in the category of smooth manifolds is
  commutative.
  \begin{equation}\label{eq:hatPsi_lifts_Psi-diagram}
    \xymatrix{
      {\PBtotal\times \SG\times (\PBbase)} \ar[d]_{\pi^{\PBtotal\times \SG,\SG}\times id_{\PBbase}}
      \ar[r]^-{\hat{\Psi}} & {\PBtotal\times \PBtotal} \ar[d]^{\pi^{\PBtotal\times \PBtotal,\SG}} \\
      {\conjB\times (\PBbase)} \ar[r]_-{\Psi} & {(\PBtotal\times \PBtotal)/\SG}
    }
  \end{equation}
\end{lemma}

\begin{proof}
  Simple exercise using the definitions and following the arrows.
\end{proof}

\begin{remark}
  It is not hard to see that, for $\Psi\in \Sigma_D'(\jcalU)$, we
  have $\Psi\in \Sigma_D(\jcalU)$ if and only if $\hat{\Psi}$
  defined in~\eqref{eq:hatPsi-def} is $\SG$-equivariant for the
  actions $l^{\PBtotal\times\SG_{ml}\times(\PBbase)}_g(q,h,r):=(q,gh,r)$ and
  $l^{\PBtotal\times \PBtotal_2}$.
\end{remark}

\begin{remark}
  All elements $\Psi\in \Sigma_D(\jcalU)$ are semi-local
  isomorphisms. Indeed, if $h_d:=(F_{RH}\circ F_{DR})(\Psi)$, so that
  $\Psi = (F_{RD}\circ F_{HR})(h_d)$, unraveling the definitions, we
  have
  \begin{equation*}
    \Psi(\pi^{\PBtotal\times\SG}(q_0,g),\PBmap(q_1)) =
    \pi^{\PBtotal\times \PBtotal,\SG}(l^{\PBtotal\times \PBtotal_2}_g(h_d(q_0,\PBmap(q_1)))).
  \end{equation*}
  Comparison of this last expression
  with~\eqref{eq:explicit_Phi_DC_and_Psi_DC} and the use of
  Proposition~\ref{prop:Phi_DC_and_Psi_DC_are_isomorphism} proves that
  $\Psi$ is a semi-local isomorphism. In addition, it proves that
  $\Psi^{-1} = \Phi\in \Sigma_U(\jcalU)$, for
  $\Phi := (F_{LU}\circ F_{CL}\circ F_{HC})(h_d)$.
\end{remark}

The next statement summarizes the results of this section concerning
the discrete Atiyah sequence (in the $\FBS_{\PBbase}$ category) of a
principal bundle and the discrete connections on the same space.

\begin{theorem}\label{thm:bijections_DAS_in_FBS}
  Let $\PBmap:\PBtotal\rightarrow \PBbase$ be a principal $\SG$-bundle and
  $\jcalU\subset \PBtotal\times \PBtotal$ be of $D$-type. Then we have the
  commutative diagram
  
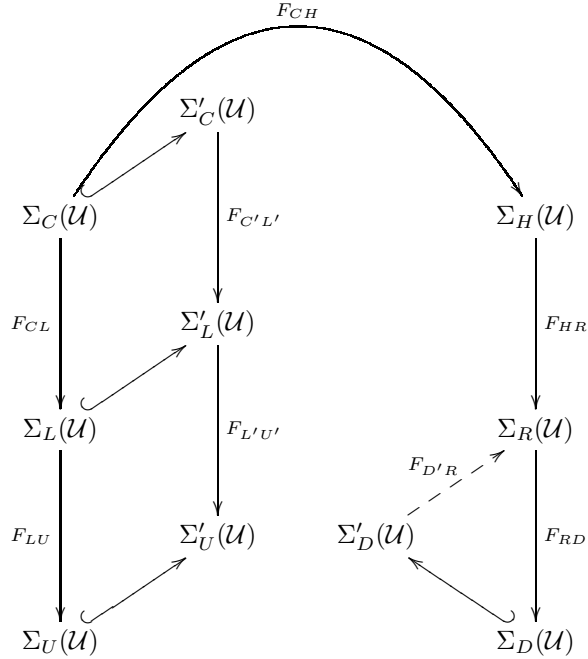
\begin{equation*}\label{eq:bijections_DAS_in_FBS-diagram}
    \xymatrix{
      {} & {\Sigma_C'(\jcalU)} \ar[dd]^{F_{C'L'}} & {} & {} \\
      {\Sigma_C(\jcalU)} \ar@/^6pc/[rrr]^{F_{CH}} \ar@{^{(}->}[ru] \ar[dd]_{F_{CL}}& {} & {} & {\Sigma_H(\jcalU)} \ar[dd]^{F_{HR}}\\
      {} & {\Sigma_L'(\jcalU)} \ar[dd]^{F_{L'U'}} & {} & {}\\
      {\Sigma_L(\jcalU)} \ar@{^{(}->}[ru] \ar[dd]_{F_{LU}} & {} & {} & {\Sigma_R(\jcalU)} \ar[dd]^{F_{RD}} \\
      {} & {\Sigma_U'(\jcalU)}  & {\Sigma_D'(\jcalU)} \ar@{-->}[ru]^{F_{D'R}} & {}\\
      {\Sigma_U(\jcalU)} \ar@{^{(}->}[ru] & {} & {} & {\Sigma_D(\jcalU)} \ar@{_{(}->}[lu]
    }
  \end{equation*}
  in the category of sets, where full arrows with straight tails are
  bijections, full arrows with curled tails are one-to-one maps while
  dashed arrows denote surjective maps. The inverse maps are obtained
  by transposing the labels of the map: for example,
  $F_{LU}^{-1} = F_{UL}$.
\end{theorem}

%%%%%%%%%%%%%%%%%%%%%%%%%%%%%%%%%%%%%%%%%%%%%%%%%%%%%%%%%% 

\section{The curvature of a discrete connection}
\label{sec:curvature_of_DC}

The curvature of a connection on a bundle has many important
properties that relate both to the geometry as well as the topology of
the underlying space. In this section we introduce a notion of
curvature for a discrete connection on a principal bundle. We motivate
our definition in the fact that continuous curvatures can be seen as
obstructions.

As it was mentioned in
Section~\ref{sec:the_atiyah_sequece-continuous_case}, a connection
$\CC$ on the principal $\SG$-bundle
$\PBmap:\PBtotal\rightarrow \PBbase$ can be seen as a certain vector
bundle map that is a morphism of Lie algebroids if and only if its
curvature $\BC$ is trivial (see Appendix A of
~\cite{bo:mackenzie-lie_groupoids_and_algebroids_in_differential_geometry}
and Remark~\ref{rem:CC_and_CB_as_obstruction} below). The curvature of
a discrete connection $\DC$ is constructed as the obstruction to
$\DC:\PBtotal\times \PBtotal\rightarrow \SG$ being a morphism of Lie
groupoids, a notion that we review next. Here we follow the notation
that is customary in Geometric Mechanics (see, for
instance,~\cite{ar:marrero_martin_martinez-discrete_lagrangian_and_hamiltonian_mechanics_on_lie_groupoids}),
that is, in some sense, opposed to the one used in
Mackenzie's~\cite{bo:mackenzie-general_theory_of_lie_groupoids_and_algebroids}.

\begin{definition}\label{def:groupoid}
  A \jdef{groupoid} over a set $M$ is a set $G$ together with the
  following maps.
  \begin{enumerate}
  \item The \jdef{source} and the \jdef{target} maps
    $\alpha,\beta:G\rightarrow M$. They define the fiber product
    $G_2 := G\FP{\beta}{\alpha} G = \{(g_1,g_2)\in G\times G:
    \beta(g_1) = \alpha(g_2)\}$.
  \item A \jdef{multiplication map} $m:G_2\rightarrow G$ (usually
    denoted by juxtaposition) such that $\alpha(g_1g_2) = \alpha(g_1)$
    and $\beta(g_1g_2) = \beta(g_2)$. It is also assumed that
    $g_1(g_2 g_3) = (g_1 g_2) g_3$ (that is, $m$ is associative,
    whenever it is defined).
  \item An \jdef{identity section} $\epsilon:M\rightarrow G$ such
    that, $\alpha \circ \epsilon = id_M = \beta\circ \epsilon$ and so
    that, for all $g\in G$, $\epsilon(\alpha(g))g = g$ and
    $g\epsilon(\beta(g)) = g$.
  \item An \jdef{inversion map} $i:G\rightarrow G$ (denoted simply by
    $i(g):=g^{-1}$) such that $g^{-1}g =\epsilon(\beta(g))$ and
    $gg^{-1} = \epsilon(\alpha(g))$.
  \end{enumerate}
  Such a groupoid is usually denoted by $G\gpoidarrows M$. A groupoid
  $G\gpoidarrows M$ is a \jdef{Lie groupoid} when $G$ and $M$ are
  smooth manifolds, $\alpha$ and $\beta$ are smooth submersions, and
  the other structure maps, $m$, $\epsilon$ and $i$, are smooth. A Lie
  groupoid is \jdef{totally intransitive} if $\alpha=\beta$, and it is
  \jdef{locally trivial} if $(\alpha,\beta):G\rightarrow M\times M$
  ---sometimes called the \jdef{anchor of $G$}--- is a surjective
  submersion.
\end{definition}

\begin{example}\label{ex:pair_groupoid-definition}
  Let $M$ be a manifold. The product $M\times M$ carries a natural
  structure of Lie groupoid over $M$ with $\alpha:=p_1$, $\beta:=p_2$,
  so that
  $(M\times M)_2 = \{((m_0,m_1),(m_1,m_2)):m_0,m_1,m_2\in M\}$ and
  \begin{equation*}
    m((m_0,m_1),(m_1,m_2)) := (m_0,m_2), \quad \epsilon(m_0):= (m_0,m_0),
    \quad i(m_0,m_1) := (m_1,m_0).
  \end{equation*}
  This groupoid is known as the \jdef{pair groupoid} over $M$.
\end{example}

\begin{example}\label{ex:lie_group_as_lie_groupoid}
  Let $\SG$ be a Lie group. Then, $\SG$ can be seen as a Lie groupoid
  over a point $\SG\gpoidarrows \{0\}$ with the multiplication and
  inversion given by the corresponding group operations. The identity
  section is the assignment $0\mapsto e$, the identity of $\SG$. The
  notion of Lie groupoid is a generalization of this example.
\end{example}

\begin{example}\label{ex:base_groupoid}
  Let $X$ be a manifold and define $\alpha=\beta=id_X$. Then
  $X\gpoidarrows X$ is a Lie groupoid with
  $X_2 = \Delta_X\subset X\times X$, so that $m(x,x):=x$,
  $\epsilon(x):=x$ and $i(x):=x$, for all $x\in X$. This groupoid is
  called the \jdef{base groupoid} and will be denoted by $X^\dagger$.
\end{example}

\begin{definition}\label{def:morphism_of_lie_groupoid}
  Given two Lie groupoids $G\gpoidarrows M$ and
  $G'\gpoidarrows M'$ a \jdef{morphism of Lie groupoids} is a
  smooth map $F:G\rightarrow G'$ such that, for all
  $(g_1,g_2)\in G_2$, $(F(g_1),F(g_2)) \in (G')_2$ and
  $F(g_1 g_2) = F(g_1) F(g_2)$. Such a morphism $F$ induces a smooth
  map $F_0:M\rightarrow M'$ in such a way that
  \begin{equation}\label{eq:morphism_of_lie_groupoid-formulas}
    \alpha'\circ F = F_0\circ \alpha,\quad \beta'\circ F = F_0\circ \beta 
    \stext{ and } F\circ \epsilon = \epsilon'\circ F_0
  \end{equation}
  (see point~\ref{it:morphism_of_loc_lie_groupoids-F_0} of
  Lemma~\ref{le:morphism_of_loc_lie_groupoid-props}).
\end{definition}

Lie groupoids and their morphisms form a category that we denote by
$\LgpdC$. If $M$ is a manifold, the Lie groupoids over $M$ and the
morphisms that induce the identity over $M$ form the subcategory
$\LgpdC_M$ of $\LgpdC$.

In order to motivate our definition of curvature, we consider first
the case of a discrete connection $\DC$ on the principal $\SG$-bundle
$\PBmap:\PBtotal\rightarrow \PBbase$ that is globally defined, that
is, $\jcalU=\PBtotal\times \PBtotal$. In this case, the discrete
connection form $\DC:\PBtotal\times \PBtotal\rightarrow \SG$ is a
smooth function between manifolds that, by
Examples~\ref{ex:pair_groupoid-definition}
and~\ref{ex:lie_group_as_lie_groupoid} happen to be Lie groupoids. It
is natural, then, to ask if $\DC$ is a morphism of groupoids.

In fact, there is a twist: rather than considering the Lie group
$(\SG,\cdot)$, we consider the \jdef{opposite group} $(\SG^{op},*)$,
that is a Lie group over the same manifold $\SG$ but with the product
$g*g':= g'\cdot g$. Then, we want to see if
$\DC\in \hom_{\LgpdC}(\PBtotal\times \PBtotal,\SG^{op})$.

As $(\SG^{op})_2=\SG\times\SG$, the condition that $\DC$ maps
multiplicable elements into multiplicable elements is always
satisfied. Then, for
$((q_0,q_1),(q_1,q_2))\in (\PBtotal\times \PBtotal)_2$ we have
\begin{gather*}
  \DC((q_0,q_1)(q_1,q_2)) = \DC(q_0,q_1) * \DC(q_1,q_2) \iff\\
  \DC(q_0,q_2) = \DC(q_1,q_2) \DC(q_0,q_1) \iff
  e = \DC(q_0,q_2)^{-1} \DC(q_1,q_2) \DC(q_0,q_1).
\end{gather*}
In general, this last condition may fail and, in order to measure its
failure, we introduce the \jdef{discrete curvature} of $\DC$. The
following definition, inspired by this analysis, is valid even when
$\DC$ is not globally defined.

\begin{definition}\label{def:discrete_curvature}
  Let $\DC$ be a discrete connection with domain $\jcalU$ on the
  principal $\SG$-bundle $\PBmap:\PBtotal\rightarrow \PBbase$. Let
  \begin{equation}\label{eq:U^{(3)}-def}
    \jcalU^{(3)}
    := \{(q_0,q_1,q_2) \in \PBtotal^3 : (q_i,q_j)\in \jcalU 
    \text{ for all } 0\leq i<j\leq 2\}.
  \end{equation}
  We define the \jdef{discrete curvature} of $\DC$ as
  $\BD: \jcalU^{(3)} \rightarrow \SG$ so that
  \begin{equation}\label{eq:BD-def}
    \BD(q_0,q_1,q_2) := \DC(q_0,q_2)^{-1}\DC(q_1,q_2)\DC(q_0,q_1).
  \end{equation}
  We say that $\DC$ is \jdef{flat} if $\BD=e$ on $\jcalU^{(3)}$.
\end{definition}

In principle, our Definition~\ref{def:discrete_curvature} is not an
exact parallel to the one used in the continuous setting. Still, a
complete parallel will be achieved later, in
Proposition~\ref{prop:s_R_is_morphism_iff_Bd=e-lLgpd} (see
Remark~\ref{rem:CC_and_CB_as_obstruction}).

\begin{remark}
  This notion of discrete curvature is, among other things, the
  obstruction to being able to trivialize a principal bundle together
  with a discrete connection, as will be discussed
  elsewhere~\cite{ar:fernandez_zuccalli-discrete_connections_on_principal_bundles_holonomy_and_curvature}. In
  the context of (discrete-time) mechanics, this notion of discrete
  curvature can be used, for example, to express the so called
  geometric phase (essentially the holonomy) around some loops when
  the structure group $\SG$ of the principal bundle is abelian
  (see~\cite{ar:fernandez_juchani_zuccalli-discrete_connections_on_principal_bundles_abelian_group_case-monteiro}
  and~\cite{ar:marsden_montgomery_ratiu-reduction_symmetry_and_phases_in_mechanics}),
  which has applications to the reconstruction of trajectories of
  systems with symmetry. Interestingly, it also appears in a seminal
  discrete version of the Ambrose--Singer Theorem (Thm. 8.1
  in~\cite{bo:kobayashi_nomizu-foundations-v1}) where the discrete
  curvature values are generators of some (discrete) holonomy groups
  (\cite{ar:fernandez_zuccalli-discrete_connections_on_principal_bundles_holonomy_and_curvature}
  and~\cite{ar:fernandez_juchani_zuccalli-discrete_connections_on_principal_bundles_abelian_group_case-monteiro}).
\end{remark}

%%%%%%%%%%%%%%%%%%%%%%%%%%%%%%%%%%%%%%%%%%%%%%%%%%%%%%%%%%

\section{The discrete Atiyah sequence in the category of local Lie groupoids}
\label{sec:DAS_in_lLgpd}

In Sections~\ref{sec:the_atiyah_sequece} and~\ref{sec:DAS_in_FBS} we
introduced and studied the discrete Atiyah sequence associated to a
principal $\SG$-bundle $\PBmap:\PBtotal\rightarrow \PBbase$ in the category $\FBS$
of fiber bundles with a section. Motivated by the properties of the
Atiyah sequence in the category of Lie algebroids, in this section we
study the diagram~\eqref{eq:discrete_atiyah_sequence-def} in the
category of Lie groupoids, focusing on its relationship with the
existence of discrete connections on $\PBmap$. In fact, as it will become
evident shortly, we will need to extend the analysis to the category
of \emph{local} Lie groupoids, because discrete connections are
usually not globally defined.

%%%%%%%%%%%%%%%%%%%%%%%%%%%%%%%

\subsection{The discrete Atiyah sequence in the category of Lie groupoids}
\label{sec:DAS_in_Lgpd}

As seen in Example~\ref{ex:pair_groupoid-definition}, the rightmost
object in~\eqref{eq:discrete_atiyah_sequence-def} is a Lie
groupoid. The next example introduces such a structure for the
object in the center.

\begin{example}\label{ex:atiyah_groupoid-definition}
  Let $\PBmap:\PBtotal\rightarrow \PBbase$ be a principal $\SG$-bundle with the
  $\SG$-action $l^\PBtotal$. It is easy to check that $(\PBtotal\times \PBtotal)/\SG$ is a
  smooth manifold and that
  $\ti{\alpha}(\pi^{\PBtotal\times \PBtotal,\SG}(q_0,q_1)):=\PBmap(q_0)$ and
  $\ti{\beta}(\pi^{\PBtotal\times \PBtotal,\SG}(q_0,q_1)):=\PBmap(q_1)$ are
  smooth surjective submersions. Then,
  \begin{equation}\label{eq:atiyah_groupoid-definition-composable}
    ((\PBtotal\times \PBtotal)/\SG)_2 =
    \{(\pi^{\PBtotal\times \PBtotal,\SG}(q_0,q_1),\pi^{\PBtotal\times
      \PBtotal,\SG}(q_1,q_2)):q_0,q_1,q_2\in \PBtotal\},
  \end{equation}
  and we define
  \begin{gather*}
    \ti{m}((\pi^{\PBtotal\times \PBtotal,\SG}(q_0,q_1),\pi^{\PBtotal\times \PBtotal,\SG}(q_1,q_2))) := 
    \pi^{\PBtotal\times \PBtotal,\SG}(q_0,q_2),\\
    \ti{\epsilon}(\PBmap(q_0)) := \pi^{\PBtotal\times \PBtotal,\SG}(q_0,q_0),\stext{ and }
    \ti{i}(\pi^{\PBtotal\times \PBtotal,\SG}(q_0,q_1)) := \pi^{\PBtotal\times \PBtotal,\SG}(q_1,q_0).
  \end{gather*}
  Standard computations show that these maps define a Lie groupoid
  $(\PBtotal\times \PBtotal)/\SG \gpoidarrows \PBbase$, that is called
  the \jdef{Atiyah} or \jdef{gauge} groupoid associated to $\PBmap$.
\end{example}

\begin{example}\label{ex:map_from_atiyah_groupoid_to_pair_as_morphism}
  It is straightforward to check that the map $F_2$ introduced
  in~\eqref{eq:discrete_atiyah_sequence-def} is a morphism of Lie
  groupoids, with the structures defined in
  Examples~\ref{ex:atiyah_groupoid-definition}
  and~\ref{ex:pair_groupoid-definition} for $(\PBtotal\times \PBtotal)/\SG$ and
  $(\PBbase)\times (\PBbase)$.  We observe that $(F_2)_0 = id_{\PBbase}$, so
  that
  $F_2\in\hom_{\LgpdC_{\PBbase}}((\PBtotal\times \PBtotal)/\SG, (\PBbase)\times
  (\PBbase))$.
\end{example}

Under certain conditions, when $G'\gpoidarrows M$ is a Lie groupoid
and $F:G\rightarrow G'$ is an embedding, there is a unique Lie
groupoid structure $G\gpoidarrows M$ for which
$F\in\hom_{\LgpdC_M}(G,G')$. An application of this construction gives
a Lie groupoid structure to the object on the left
of~\eqref{eq:discrete_atiyah_sequence-def} for which $F_1$ is a
morphism of Lie groupoids.

\begin{example}\label{ex:gpoid_struct_on_QxG/G}
  The maps
  $\alpha_{\conjB}, \beta_{\conjB}:\conjB\rightarrow \PBbase$ defined by
  \begin{equation*}
    \alpha_{\conjB}(\pi^{\PBtotal\times \SG,\SG}(q,g)) := \PBmap(q)
    \stext{ and } \beta_{\conjB} := \alpha_{\conjB}
  \end{equation*}
  together with
  \begin{gather*}
    \conjB_2 = \{(\pi^{\PBtotal\times\SG,\SG}(q,g_0),
    \pi^{\PBtotal\times\SG,\SG}(q,g_1)):q\in \PBtotal,\ g_0,g_1\in\SG\},\\
    m(\pi^{\PBtotal\times \SG,\SG}(q_0,g_0), \pi^{\PBtotal\times
      \SG,\SG}(q_0,g_1)) := \pi^{\PBtotal\times
      \SG,\SG}(q_0,g_0g_1)\\
    \epsilon(\PBmap(q)) := \pi^{\PBtotal\times\SG,\SG}(q,e), \stext{ and }
    i(\pi^{\PBtotal\times \SG,\SG}(q,g)) := \pi^{\PBtotal\times \SG,\SG}(q,g^{-1})
  \end{gather*}
  define a Lie groupoid $\conjB\gpoidarrows \PBbase$. In addition, the
  map $F_1$ defined in~\eqref{eq:discrete_atiyah_sequence-def} is a
  morphism of Lie groupoids for $\conjB$ with this structure and
  $(\PBtotal\times \PBtotal)/\SG$ with the structure given in
  Example~\ref{ex:atiyah_groupoid-definition}.
\end{example}

\begin{definition}\label{def:extensionM-LGpd}
  Given $G_1,G_3$ in $\LgpdC_M$ such that $G_1$ is totally
  intransitive and $G_3$ is locally trivial, a sequence
  $G_1\xrightarrow{\eta_1} G_2 \xrightarrow{\eta_2} G_3$ in $\LgpdC_M$
  is said to be an \jdef{extension} of $G_3$ by $G_1$ if $\eta_1$ is
  an embedding, $\eta_2$ is a surjective submersion, and, for
  $\im(\eta_1):=\eta_1(G_1)$ and
  $\ker(\eta_2):=\eta_2^{-1}(\epsilon_{G_3}(M))$,
  $\im(\eta_1)=\ker(\eta_2)$ (as sets).
\end{definition}

\begin{remark}
  There is no completely satisfactory notion of exact sequence in the
  category of Lie groupoids. One version of exactness was considered
  by Mackenzie
  in~\cite[Prop. 3.15]{bo:mackenzie-lie_groupoids_and_algebroids_in_differential_geometry}:
  it states that the sequence
  $G_1\xrightarrow{\eta_1} G_2 \xrightarrow{\eta_2} G_3$ in $\LgpdC_M$
  is \jdef{exact} if $\eta_1$ is an embedding, $\eta_2$ is a
  surjective submersion, and $\im(\eta_1) = \ker(\eta_2)$, where
  $\ker(\eta_2)=\eta_2^{-1}(\epsilon_2(M))$. Yet, in~\cite[Defniition
  1.7.15]{bo:mackenzie-general_theory_of_lie_groupoids_and_algebroids},
  he introduces essentially the same notion under the name
  \jdef{extension}. This is our
  Definition~\ref{def:extensionM-LGpd}. Of course, in abelian
  categories, exact sequences and extensions are basically the same
  notions, only the perspective changes.
\end{remark}

\begin{proposition}\label{prop:DAS_is_extensionM_in_Lgpd}
  The discrete Atiyah sequence~\eqref{eq:discrete_atiyah_sequence-def}
  is an extension of $(\PBbase)\times (\PBbase)$ by $\conjB$
  (Definition~\ref{def:extensionM-LGpd}).
\end{proposition}

\begin{proof}
  It is the special case of
  Proposition~\ref{prop:DAS_is_extensionM-lLGpd} where
  $\jcalU=\PBtotal\times \PBtotal$.
\end{proof}

The natural next step is to ask if, given a discrete connection $\DC$
with domain $\jcalU$ on the principal $\SG$-bundle $\PBmap$, the maps
$s_L$ and $s_R$ defined by~\eqref{eq:s_L_from_A-def}
and~\eqref{eq:s_R_from_h-def} are morphisms of Lie groupoids. It is
immediately clear that, unless $\DC$ is globally defined, neither
$s_L$ nor $s_R$ can be morphisms of groupoids (which are globally
defined). Thus, we work in a more general setting, that of \jdef{local
  Lie groupoids}, that we introduce next.

%%%%%%%%%%%%%%%%%%%%

\subsection{Local Lie groupoids}
\label{sec:lLgpd}

There are several slightly different notions of \emph{local} Lie
groupoid. We essentially follow the one introduced in p.40
of~\cite{ar:coste_dazord_weinstein-grupoides_symplectiques}.

\begin{definition}\label{def:local_lie_groupoid}
  A \jdef{local Lie groupoid} consists of smooth manifolds $G$ and $M$
  together with submersions $\alpha,\beta:G\rightarrow M$ as well as a
  diffeomorphism $i:G\rightarrow G$, and smooth maps
  $\epsilon:M\rightarrow G$ and $m:G_m\rightarrow G$, where
  $G_m\subset G_2:=G\FP{\beta}{\alpha} G$ is an open subset, all
  subject to the conditions stated below. As for groupoids,
  $g_1g_2:=m(g_1,g_2)$ and $g^{-1}:=i(g)$.
  \begin{enumerate}
  \item \label{it:lLgpd-epsilon}
    $\alpha\circ \epsilon = id_M = \beta\circ \epsilon$.
  \item \label{it:lLgpd-inverse} For all $(g_1,g_2)\in G_m$, we have
    that $(g_2^{-1}, g_1^{-1}) \in G_m$ and
    $g_2^{-1}g_1^{-1} = (g_1g_2)^{-1}$.
  \item \label{it:lLgpd-identity} For all $g\in G$,
    $(\epsilon(\alpha(g)),g), (g,\epsilon(\beta(g))) \in G_m$ and
    $\epsilon(\alpha(g)) g = g = g \epsilon(\beta(g))$.
  \item \label{it:lLgpd-inverses} For all $g\in G$,
    $(g,g^{-1}),(g^{-1},g)\in G_m$ and $gg^{-1} = \epsilon(\alpha(g))$
    and $g^{-1}g=\epsilon(\beta(g))$.
  \item \label{it:lLgpd-associativity} If
    $(g_1,g_2), (g_2,g_3), (g_1,g_2g_3)\in G_m$, then
    $(g_1g_2,g_3) \in G_m$ and $(g_1g_2)g_3 = g_1(g_2g_3)$.
  \end{enumerate}
  A local Lie groupoid is \jdef{totally intransitive} if
  $\alpha=\beta$, and it is \jdef{locally trivial} if
  $(\alpha,\beta):G\rightarrow M\times M$ is a surjective submersion.
\end{definition}

\begin{remark}\label{rem:totally_intransitive_as_local_Lie_group_bundles}
  It is easy to see that a totally intransitive local Lie groupoid is
  a local groupoid where multiplication is only possible within each
  fiber and, so, each fiber is a local Lie group. This is the local
  version of what, in the $\LgpdC$ category, are called \jdef{Lie group
    bundles}
  in~\cite{bo:mackenzie-general_theory_of_lie_groupoids_and_algebroids}.
\end{remark}

Many results that are valid for Lie groupoids have natural
analogues for local Lie groupoids. We include some of them
in~\ref{sec:some_basic_and_auxiliary_properties}, and use them freely
in what follows.

\begin{remark}\label{rem:epsilon(M)_2_and_G_m_are_nonempty}
  When $G$ is a local Lie groupoid over $M$, it is easy to see using
  point~\ref{it:lLgpd_basic_properties-epsilon(M)_2} of
  Lemma~\ref{le:lLgpd_basic_properties} that $G_m\neq \emptyset$ and,
  consequently, that $G_2\neq \emptyset$.
\end{remark}

Local Lie groupoids generalize the idea of Lie groupoid. The next two
statements make it more precise.

\begin{lemma}\label{le:global_lLgpd_imp_Lgpd}
  If $G$ is a local Lie groupoid over $M$ with $G_m=G_2$, then it is a
  Lie groupoid over $M$.
\end{lemma}

\begin{proof}
  Easy verification, using
  point~\ref{it:lLgpd_basic_properties-alpha_beta_m} of
  Lemma~\ref{le:lLgpd_basic_properties}.
\end{proof}

\begin{lemma}\label{le:open_in_Lgpd_is_lLgpd}
  Let $G\gpoidarrows M$ be a Lie groupoid and $U\subset G$ be an open
  subset such that $U^{-1}\subset U$ and $\epsilon(M)\subset U$. Then
  $U$ is a local Lie groupoid over $M$ with the structure maps
  $\alpha_U$, $\beta_U$, $i_U$, and $\epsilon_U$ induced by those of
  $G$ by restriction or co-restriction and where the multiplication
  $m_U$ is defined as the restriction and co-restriction of the
  multiplication $m$ on $G$ to $U_m:=U_2\cap m^{-1}(U)$ and $U$
  respectively. In particular, every Lie groupoid $G$ is a local Lie
  groupoid with the multiplication globally defined, that is, with
  $G_m=G_2$.
\end{lemma}

\begin{proof}
  Lengthy but mechanical verification following
  Definition~\ref{def:local_lie_groupoid}.
\end{proof}

\begin{example}\label{ex:gpoid_struct_on_QxG/G-lLgpd}
  As seen in Example~\ref{ex:gpoid_struct_on_QxG/G},
  $\conjB\gpoidarrows \PBbase$ is a Lie groupoid so that, by
  Lemma~\ref{le:open_in_Lgpd_is_lLgpd}, it is also a local Lie
  groupoid. In addition, as the structure maps are the same for both
  structures and $\alpha_{\conjB} = \beta_{\conjB}$, it is totally
  intransitive.
\end{example}

\begin{example}\label{ex:lLgpd_associated_with_DC}
  Let $\PBmap:\PBtotal\rightarrow \PBbase$ be a principal $\SG$-bundle; recall
  that, as seen in Examples~\ref{ex:atiyah_groupoid-definition}
  and~\ref{ex:pair_groupoid-definition}, $(\PBtotal\times \PBtotal)/\SG$ and
  $(\PBbase)\times (\PBbase)$ are Lie groupoids over $\PBbase$. Let
  $\jcalU\subset \PBtotal\times \PBtotal$ be a symmetric open subset such that
  $\Delta_\PBtotal\subset\jcalU$. Then
  $\jcalU'':=(\PBmap\times \PBmap)(\jcalU) \subset \PBbase \times \PBbase$
  and $\jcalU/\SG \subset (\PBtotal\times \PBtotal)/\SG$ are open subsets satisfying
  the conditions of Lemma~\ref{le:open_in_Lgpd_is_lLgpd}. Hence,
  $\jcalU''$ and $\jcalU/\SG$ are local Lie groupoids over $\PBbase$
  that are, usually, not Lie groupoids. In addition, when $\jcalU$ is
  of $pD$-type,
  \begin{equation}\label{eq:lLgpd_associated_with_DC-pD_type-multi_domain-U''}
    \begin{split}
      (\jcalU'')_m =
      \{((\PBmap(q_0),\PBmap(q_1)),&(\PBmap(q_1),\PBmap(q_2)))\in
      (\PBbase\times \PBbase)^2 :\\& (q_0,q_1),(q_1,q_2)\in\jcalU \text{ and }
      (q_0,q_2) \in l^{\PBtotal\times \PBtotal_2}_\SG(\jcalU)\},
    \end{split}
  \end{equation}
  and
  \begin{equation}\label{eq:lLgpd_associated_with_DC-pD_type-multi_domain-U/G}
    \begin{split}
      (\jcalU/\SG)_m = \{ (\pi^{\PBtotal\times \PBtotal,\SG}(q_0,q_1),\pi^{\PBtotal\times
        \PBtotal,\SG}(q_1,q_2)) \in ((\PBtotal\times \PBtotal)/\SG)^2 : (q_0,q_1,q_2)\in
      \jcalU^{(3)}\},
    \end{split}
  \end{equation}
  where $\jcalU^{(3)}$ is as in~\eqref{eq:U^{(3)}-def}. If $\jcalU$ is
  of
  $D$-type,~\eqref{eq:lLgpd_associated_with_DC-pD_type-multi_domain-U''}
  turns into
  \begin{equation}\label{eq:lLgpd_associated_with_DC-D_type-multi_domain-U''}
    \begin{split}
      (\jcalU'')_m =
      \{((\PBmap(q_0),\PBmap(q_1)),(\PBmap(q_1),\PBmap(q_2)))\in (\PBbase\times \PBbase)^2
      : (q_0,q_1,q_2)\in \jcalU^{(3)}\}.
    \end{split}
  \end{equation}
\end{example}

\begin{remark}\label{rem:spaces_in_DAS_as_lLgpds}
  When $\PBmap:\PBtotal\rightarrow \PBbase$ is a principal $\SG$-bundle, $\conjB$
  is a local Lie groupoid over $\PBbase$
  (Example~\ref{ex:gpoid_struct_on_QxG/G-lLgpd}). Also, if
  $\jcalU\subset \PBtotal\times \PBtotal$ is symmetric of $pD$-type, by
  Example~\ref{ex:lLgpd_associated_with_DC}, both $\jcalU''$ and
  $\jcalU/\SG$ are local Lie groupoids over $\PBbase$. In what follows,
  whenever we consider any of these three spaces as local Lie
  groupoids, it will be with the structures just mentioned.
\end{remark}

\begin{lemma}\label{le:epsilon(M)_i_invariant}
  Let $G\gpoidarrows M$ be a local Lie groupoid. Then
  $\epsilon(m)^{-1} = \epsilon(m)$ for each $m\in M$.
\end{lemma}

\begin{proof}
  For each $m\in M$ we have $\alpha(\epsilon(m))=m$, so that by
  property~\ref{it:lLgpd-identity} of
  Definition~\ref{def:local_lie_groupoid},
  $\epsilon(m)\epsilon(m) = \epsilon(m)$. The result follows by
  point~\ref{it:uniqueness_in_loc_lie_groupoid-right_identity} of
  Proposition~\ref{prop:uniqueness_in_loc_lie_groupoid}.
\end{proof}
  
%%%%%%%%%%%%%%%%%%%%%%

\subsection{Morphisms of local Lie groupoids}
\label{sec:morphisms_of_local_lie_groupoids}

\begin{definition}
  Let $G\gpoidarrows M$ and $G'\gpoidarrows M'$ be local Lie
  groupoids. A smooth map $F:G\rightarrow G'$ is a \jdef{morphism of
    local Lie groupoids} if
  \begin{enumerate}
  \item $F(\epsilon_G(M))\subset \epsilon_{G'}(M')$, and
  \item for all $(g_1,g_2)\in G_m$, we have that
    $(F(g_1),F(g_2))\in G'_m$ and $F(g_1g_2)=F(g_1)F(g_2)$.
  \end{enumerate}
  A morphism of local Lie groupoids is said to be a \jdef{submersion}
  if it is a submersion as a map between smooth manifolds.
\end{definition}

\begin{example}\label{ex:F_1_and_F_2_are_morphisms_of_lLgpd}
  Let $\PBmap:\PBtotal\rightarrow \PBbase$ be a principal $\SG$-bundle. Given
  $\jcalU \subset \PBtotal\times \PBtotal$, a symmetric $pD$-type subset, the
  maps $F_1:\conjB\rightarrow \jcalU/\SG$ and
  $F_2:\jcalU/\SG\rightarrow \jcalU''$ are defined as the
  restriction and co-restriction of the homonymous maps
  in~\eqref{eq:discrete_atiyah_sequence-def}. Viewing both domains and
  codomains as local Lie groupoids (see
  Remark~\ref{rem:spaces_in_DAS_as_lLgpds}) it is easy to check by
  direct computation that $F_1$ and $F_2$ are morphisms of local Lie
  groupoids. This result also follows from
  Lemma~\ref{le:restriction_of_morphisms_to_open_subsets}.
\end{example}

Next we analyze if the left and right semi-local splittings of the
discrete Atiyah sequence (in the category $\FBS_{\PBbase}$) can be
morphisms of local Lie groupoids.

\begin{example}\label{ex:s_L_almost_never_morphism-lLgpd}
  Let $\PBmap:\PBtotal\rightarrow \PBbase$ be a connected principal $\SG$-bundle,
  $\jcalU\subset \PBtotal\times \PBtotal$ be a symmetric $pD$-type subset and
  $s_L:\jcalU/\SG\rightarrow \conjB$ a semi-local
  left splitting of~\eqref{eq:discrete_atiyah_sequence-def}, that is,
  $s_L\in \Sigma_L'(\jcalU)$. As seen in
  Examples~\ref{ex:lLgpd_associated_with_DC}
  and~\ref{ex:gpoid_struct_on_QxG/G} (together with
  Lemma~\ref{le:global_lLgpd_imp_Lgpd}) $\jcalU/\SG$ and
  $\conjB$ are local Lie groupoids over $\PBbase$. We ask if
  $s_L$ is a morphism of local Lie groupoids. As $s_L$ is a semi-local
  morphism,
  $\phi_{\conjB} \circ s_L = \phi_{\jcalU/\SG}$ and we
  see that, for $(q_0,q_1)\in \jcalU$,
  \begin{equation}
    \label{eq:s_L_almost_never_morphism-lLgpd-projections}
    \PBmap(q_0) = \phi_{\jcalU/\SG}(\pi^{\PBtotal\times \PBtotal,\SG}(q_0,q_1))
    =\phi_{\conjB}(s_L(\pi^{\PBtotal\times \PBtotal,\SG}(q_0,q_1))).
  \end{equation}
  Using the description of the set $(\jcalU/\SG)_m$ provided
  by~\eqref{eq:lLgpd_associated_with_DC-pD_type-multi_domain-U/G}, for
  $(\pi^{\PBtotal\times \PBtotal,\SG}(q_0,q_1),\\\pi^{\PBtotal\times
    \PBtotal,\SG}(q_1,q_2)) \in (\jcalU/\SG)_m$, that is, satisfying
  $(q_0,q_1,q_2)\in \jcalU^{(3)}$,
  \begin{equation*}
    \begin{split}
      (s_L(\pi^{\PBtotal\times \PBtotal,\SG}(q_0,q_1)), s_L(\pi^{\PBtotal\times
        \PBtotal,\SG}(q_1,q_2))) \in \conjB_m = \conjB_2
    \end{split}
  \end{equation*}
  if and only if
  \begin{equation*}
    \beta_{\conjB}(s_L(\pi^{\PBtotal\times \PBtotal,\SG}(q_0,q_1))) =
    \alpha_{\conjB}(s_L(\pi^{\PBtotal\times
      \PBtotal,\SG}(q_1,q_2))).
  \end{equation*}
  As $\alpha_{\conjB}=\beta_{\conjB}=\phi_{\conjB}$,
  using~\eqref{eq:s_L_almost_never_morphism-lLgpd-projections}, this
  last condition becomes
  \begin{equation*}
    \PBmap(q_0) = \phi_{\conjB}(s_L(\pi^{\PBtotal\times \PBtotal,\SG}(q_0,q_1))) =
    \phi_{\conjB}(s_L(\pi^{\PBtotal\times
      \PBtotal,\SG}(q_1,q_2))) = \PBmap(q_1).
  \end{equation*}
  Thus, equivalently, if and only if $(q_0,q_1)\in {\mathcal
    V}_d$. Hence, if $s_L$ is a morphism of local Lie groupoids, it
  should be ${\mathcal V}_d \supset \jcalU$, thus,
  ${\mathcal V}_d = \jcalU$. But then, ${\mathcal V}_d$ is both
  open and closed in the connected space $\PBtotal\times \PBtotal$, so that
  ${\mathcal V}_d=\PBtotal\times \PBtotal$. From here it is not hard to see that
  $\PBmap:\PBtotal\rightarrow \PBbase$ is isomorphic to
  $\pi^{\SG,\SG}:\SG\rightarrow \{[e]\}$ for the left multiplication
  in $\SG$ action; in other words, $\PBtotal$ is a (left) torsor for
  $\SG$. Thus, in almost all cases, $s_L$ cannot be a morphism of
  local Lie groupoids.
\end{example}

The question of whether a semi-local right splitting $s_R$
of~\eqref{eq:discrete_atiyah_sequence-def} is a morphism of local Lie
groupoids is more interesting and, while we will see that it is not
always the case, the obstruction is, precisely, the discrete curvature
of the discrete connection $\DC := F_{HC}(F_{RH}(s_R))$. This is
explored in the next result.

\begin{proposition}\label{prop:s_R_is_morphism_iff_Bd=e-lLgpd}
  Let $\jcalU\subset \PBtotal\times \PBtotal$ be a symmetric $D$-type subset
  and let $s_R:\jcalU''\rightarrow \jcalU/\SG$ be a
  semi-local right splitting
  of~\eqref{eq:discrete_atiyah_sequence-def} (in the category
  $\FBS_{\PBbase}$). Define the discrete connection
  $\DC:= F_{HC}(F_{RH}(s_R)) \in \Sigma_C(\jcalU)$. Then $s_R$
  is a morphism of local Lie groupoids if and only if $\DC$ is
  flat.
\end{proposition}

\begin{proof}
  As $s_R$ is a semi-local morphism in $\FBS_{\PBbase}$, we have
  $s_R \circ \sigma_{(\PBbase)\times (\PBbase)} = \sigma_{(\PBtotal\times
    \PBtotal)/\SG}$ and, as $\sigma_\cdot = \epsilon_\cdot$ (the Lie groupoid
  structure map) in the current case, we have that
  $s_R(\epsilon_{\jcalU''}(\PBbase)) = \epsilon_{{\mathcal
      U}/\SG}(\PBbase)$.

  We have $s_R = F_{HR}(F_{CH}(\DC))$ for $\DC$ as in the
  statement. Then, by~\eqref{eq:HLd_from_DC}
  and~\eqref{eq:s_R_from_h-def},
  \begin{equation}\label{eq:s_R_is_morphism_iff_Bd=e-lLgpd-s_R_from_DC}
    s_R(\PBmap(q_0),\PBmap(q_1)) =
    \pi^{\PBtotal\times \PBtotal,\SG}(q_0,l^\PBtotal_{\DC(q_0,q_1)^{-1}}(q_1))
    \stext{ for all } (q_0,q_1)\in\jcalU.
  \end{equation}

  As $\jcalU$ is of $D$-type,
  recalling~\eqref{eq:lLgpd_associated_with_DC-D_type-multi_domain-U''},
  all elements of $(\jcalU'')_m$ are of the form
  \begin{equation*}
    ((\PBmap(q_0),\PBmap(q_1)),(\PBmap(q_1),\PBmap(q_2))) \stext{ with }
    (q_0,q_1,q_2)\in \jcalU^{(3)}.
  \end{equation*}
  For any of them,
  using~\eqref{eq:s_R_is_morphism_iff_Bd=e-lLgpd-s_R_from_DC}, we have
  that
  \begin{equation}
    \label{eq:s_R_is_morphism_iff_Bd=e-lLgpd-s_R_calc}
    \begin{split}
      (s_R(&\PBmap(q_0),\PBmap(q_1)),s_R(\PBmap(q_1),\PBmap(q_2))) \\=&
      (\pi^{\PBtotal\times \PBtotal,\SG}(q_0,l^\PBtotal_{\DC(q_0,q_1)^{-1}}(q_1)) ,
      \pi^{\PBtotal\times \PBtotal,\SG}(q_1,l^\PBtotal_{\DC(q_1,q_2)^{-1}}(q_2))) \\=&
      (\pi^{\PBtotal\times \PBtotal,\SG}(q_0,l^\PBtotal_{\DC(q_0,q_1)^{-1}}(q_1)) ,
      \pi^{\PBtotal\times \PBtotal,\SG}(l^\PBtotal_{\DC(q_0,q_1)^{-1}}(q_1),
      l^\PBtotal_{\DC(q_0,q_1)^{-1}\DC(q_1,q_2)^{-1}}(q_2))).
    \end{split}
  \end{equation}
  It is easy to check,
  using~\eqref{eq:lLgpd_associated_with_DC-pD_type-multi_domain-U/G},
  that this last expression is in $(\jcalU/\SG)_m$. In this
  case,
  \begin{equation*}
    \begin{split}
      s_R((\PBmap(q_0),\PBmap(q_1))(\PBmap(q_1),\PBmap(q_2))) =
      s_R(\PBmap(q_0),\PBmap(q_2)) = \pi^{\PBtotal\times
        \PBtotal,\SG}(q_0,l^\PBtotal_{\DC(q_0,q_2)^{-1}}(q_2)),
    \end{split}
  \end{equation*}
  while, from~\eqref{eq:s_R_is_morphism_iff_Bd=e-lLgpd-s_R_calc},
  \begin{equation*}
    \begin{split}
      s_R(\PBmap(q_0),&\PBmap(q_1)) s_R(\PBmap(q_1),\PBmap(q_2)) =
      \pi^{\PBtotal\times
        \PBtotal,\SG}(q_0,l^\PBtotal_{\DC(q_0,q_1)^{-1}\DC(q_1,q_2)^{-1}}(q_2)).
    \end{split}
  \end{equation*}
  Comparing the last two expressions we have that
  \begin{equation*}
    s_R((\PBmap(q_0),\PBmap(q_1))(\PBmap(q_1),\PBmap(q_2))) =
    s_R(\PBmap(q_0),\PBmap(q_1))s_R(\PBmap(q_1),\PBmap(q_2))
  \end{equation*}
  if and only if
  $\DC(q_0,q_2)^{-1} = \DC(q_0,q_1)^{-1} \DC(q_1,q_2)^{-1}$ for all
  $(q_0,q_1,q_2) \in \jcalU^{(3)}$ or, recalling the definition
  of discrete curvature~\eqref{eq:BD-def}, if and only if
  $\BD(q_0,q_1,q_2) = e$. Thus, $s_R$ is a morphism of local Lie
  groupoids if and only if $\DC$ is flat.
\end{proof}

\begin{definition}
  All local Lie groupoids and their morphisms make a category that we
  denote by $\lLgpdC$. For a manifold $M$, if we consider only
  groupoids over $M$ and morphisms between those groupoids
  $F\in\hom_{\lLgpdC}(G_1,G_2)$ such that $F_0=id_M$, we obtain the
  subcategory $\lLgpdC_M$ of $\lLgpdC$.
\end{definition}

\begin{remark}\label{rem:CC_and_CB_as_obstruction}
  A connection $\CC$ on the principal $\SG$-bundle
  $\PBmap:\PBtotal\rightarrow \PBbase$ can be seen as a right splitting (in the
  category of vector bundles) of the Atiyah
  sequence~\eqref{eq:atiyah_sequence-def} (Definition on p. 188
  of~\cite{ar:atiyah-complex_analytic_connections_in_fibre_bundles} as
  well as Definition 4.1 in Appendix A
  of~\cite{bo:mackenzie-lie_groupoids_and_algebroids_in_differential_geometry}). Following
  Mackenzie (Definition 4.10 in Appendix A
  of~\cite{bo:mackenzie-lie_groupoids_and_algebroids_in_differential_geometry}),
  the curvature $\BC$ of that connection $\CC$ is the morphism of
  vector bundles $\BC:T(\PBbase)\oplus T(\PBbase)\rightarrow \ti{\jgsg}$
  defined on sections by
  \begin{equation*}
    \phi_1(\BC(\xi_1,\xi_2)) := \CC([\xi_1,\xi_2]) - [\CC(\xi_1),\CC(\xi_2)]
  \end{equation*}
  for $\xi_1,\xi_2\in \VF(\PBbase)$ and where $\phi_1$ is the injective
  morphism appearing in~\eqref{eq:atiyah_sequence-def}. Thus, $\BC$ is
  the obstruction to $\CC([\xi_1,\xi_2]) = [\CC(\xi_1),\CC(\xi_2)]$
  or, in other words, to $\CC$ being a morphism of Lie algebroids,
  because $\CC$ always preserves the corresponding anchor maps. The
  equivalence of this definition of curvature to the more traditional
  one is given by Proposition 4.16 in Appendix A
  of~\cite{bo:mackenzie-lie_groupoids_and_algebroids_in_differential_geometry}. In
  this respect, Proposition~\ref{prop:s_R_is_morphism_iff_Bd=e-lLgpd},
  shows that the curvature of a discrete connection and the curvature
  of a connection are both the obstructions to analogous problems in
  the categories of $\lLgpdC$ and of Lie algebroids.
\end{remark}

\begin{proposition}\label{prop:DAS_for_pDtype_U_is_sequence_in_lLgpd}
  Let $\PBmap:\PBtotal\rightarrow \PBbase$ be a principal $\SG$-bundle and
  $\jcalU\subset \PBtotal\times \PBtotal$ be of symmetric $pD$-type. Then,
  \begin{enumerate}
  \item \label{it:DAS_for_pDtype_U_is_sequence_in_lLgpd-sequence} for
    the maps $F_1$ and $F_2$ defined in
    Example~\ref{ex:F_1_and_F_2_are_morphisms_of_lLgpd},
    \begin{equation}
      \label{eq:DAS_in_lLgpd}
      \xymatrix{{\conjB} \ar[r]^-{F_1} & {\jcalU/\SG} \ar[r]^-{F_2} &
        {\jcalU''}}
    \end{equation}
    is a sequence in $\lLgpdC_{\PBbase}$ and
  \item \label{it:DAS_for_pDtype_U_is_sequence_in_lLgpd-correspondence}
    if $\jcalU$ is of $D$-type and
    $s_R\in \Sigma_R(\jcalU)$ is such that
    $\DC:=F_{HC}(F_{RH}(s_R))\in \Sigma_C(\jcalU)$ is flat, then
    $s_R$ is a right splitting of~\eqref{eq:DAS_in_lLgpd} in
    $\lLgpdC_{\PBbase}$.
  \end{enumerate}
\end{proposition}

\begin{proof}
  That $\conjB$, $\jcalU/\SG$ and $\jcalU''$ are objects in
  $\lLgpdC_{\PBbase}$ was established by
  Lemma~\ref{le:open_in_Lgpd_is_lLgpd} and
  Example~\ref{ex:lLgpd_associated_with_DC}. That $F_1$ and $F_2$ are
  morphisms in $\lLgpdC_{\PBbase}$ is in
  Example~\ref{ex:F_1_and_F_2_are_morphisms_of_lLgpd}. On the other
  hand, if $s_R$ corresponds to a flat $\DC$,
  $s_R\in \hom_{\lLgpdC_{\PBbase}}(\jcalU'',\jcalU/\SG)$ by
  Proposition~\ref{prop:s_R_is_morphism_iff_Bd=e-lLgpd}. Last,
  $F_2\circ s_R = id_{\jcalU''}$ because $s_R$ is a right inverse of
  $F_2$ in $\FBS_{\PBbase}$.
\end{proof}

%%%%%

Next we introduce in $\lLgpdC_M$ the analogue of the notion of
extension in the category of Lie groupoids
(Definition~\ref{def:extensionM-LGpd}).

\begin{definition}\label{def:extensionM-lLGpd}
  Given $G_1,G_3\in\ob_{\lLgpdC_M}$ such that $G_1$ is totally
  intransitive and $G_3$ is locally trivial, a sequence
  $G_1\xrightarrow{\eta_1} G_2 \xrightarrow{\eta_2} G_3$ in
  $\lLgpdC_M$ is said to be an \jdef{extension} of $G_3$ by $G_1$ if
  $\eta_1$ is an embedding such that
  $(\eta_1\times \eta_1)^{-1}((G_2)_m)\subset (G_1)_m$, $\eta_2$ is a
  surjective submersion, and $\im(\eta_1)=\ker(\eta_2)$ (as sets). An
  extension as above is said to be \jdef{right split} if $\eta_2$ has
  a right inverse and \jdef{left split} if $\eta_1$ has a left
  inverse.
\end{definition}

\begin{remark}\label{rem:images_and_kernels_are_usually_lLgpd}
  It can be shown that when $\eta\in\hom_{\lLgpdC_M}(G,G')$ is such
  that $\eta:G\rightarrow G'$ is an embedding, the set
  $\im(\eta):=\eta(G)\subset G'$ is a local Lie groupoid over $M$ with
  the structure maps induced by those of $G'$, thus $\im(\eta)$
  together with the inclusion $j_{\im(\eta)}:\im(\eta)\rightarrow G'$
  is an embedded wide local Lie subgroupoid of $G'$. If, in addition,
  $(\eta\times\eta)^{-1}((G')_m)\subset G_m$, then
  $\eta|^{\im(\eta)}\in\hom_{\lLgpdC_M}(G,\im(\eta))$ is an
  isomorphism. Also, when $\eta'\in\hom_{\lLgpdC_M}(G,G')$ is a
  submersion, the set $K:=\ker(\eta')$ is a totally intransitive local
  Lie groupoid over $M$ such that the inclusion $j_K:K\rightarrow G$
  is in $\hom_{\lLgpdC_M}(K,G)$. Having these local Lie groupoid
  structures in mind, if the sequence
  $G_1\xrightarrow{\eta_1} G_2 \xrightarrow{\eta_2} G_3$ in
  $\lLgpdC_M$ is an extension, it can be seen that
  $\im(\eta_1)=\ker(\eta_2)$ in $\lLgpdC_M$.
\end{remark}

\begin{remark}\label{rem:preimage_of_multiplicable_elements_in_Lgpd}
  When $\eta\in\hom_{\lLgpdC_M}(G,G')$ and $G$ is a Lie groupoid, the
  condition $(\eta\times\eta)^{-1}((G')_m)\subset G_m$ is always
  satisfied.
\end{remark}

\begin{proposition}\label{prop:DAS_is_extensionM-lLGpd}
  Let $\PBmap:\PBtotal\rightarrow \PBbase$ be a principal $\SG$-bundle
  and $\jcalU\subset \PBtotal\times \PBtotal$ be of symmetric
  $pD$-type. Then, the DAS on $\jcalU$~\eqref{eq:DAS_in_lLgpd} is an
  extension in the $\lLgpdC_{\PBbase}$ category.
\end{proposition}

\begin{proof}
  Proposition~\ref{prop:DAS_for_pDtype_U_is_sequence_in_lLgpd} proves
  that~\eqref{eq:DAS_in_lLgpd} is a sequence in $\lLgpdC_{\PBbase}$ and
  we know that $\conjB$ is totally intransitive
  (Example~\ref{ex:gpoid_struct_on_QxG/G-lLgpd}). As
  $(\alpha_{\jcalU''},\beta_{\jcalU''}) = id_{\jcalU''}$, we see that
  $\jcalU''$ is locally trivial. Also, as
  $F_2:(\PBtotal\times \PBtotal)/\SG\rightarrow (\PBbase)\times (\PBbase)$ is a
  surjective submersion
  (Lemma~\ref{le:F_2_is_a_surjective_submersion}) and
  $F_2(\jcalU/\SG)=\jcalU''$ with both $\jcalU/\SG$ and $\jcalU''$
  open subsets, we have that $F_2:\jcalU/\SG\rightarrow \jcalU''$ is a
  surjective submersion. Similarly,
  $F_1:\conjB\rightarrow (\PBtotal\times \PBtotal)/\SG$ is an embedding
  (Example~\ref{ex:DAS_is_extensionM_in_FBS}) and
  $\jcalU/\SG\subset (\PBtotal\times \PBtotal)/\SG$ is open and contains
  $\VD/\SG = F_1(\conjB)$, we see that
  $F_1:\conjB\rightarrow \jcalU/\SG$ is an embedding. Condition
  $(F_1\times F_1)^{-1}((\jcalU/\SG)_m)\subset (\conjB)_m$ holds, as
  observed in
  Remark~\ref{rem:preimage_of_multiplicable_elements_in_Lgpd}, because
  $\conjB$ is a Lie groupoid.

  Last, the result follows from
  \begin{equation*}
    \begin{split}
      \ker(F_2) =& F_2^{-1}(\epsilon_{\jcalU''}(\PBbase)) =
      F_2^{-1}(\Delta_{\PBbase}) = \{\pi^{\PBtotal\times \PBtotal,\SG}(q_0,q_1) :
      \PBmap(q_0) = \PBmap(q_1)\} \\=& \{\pi^{\PBtotal\times
        \PBtotal,\SG}(q,l^\PBtotal_g(q)):q\in \PBtotal \text{ and } g\in \SG\} =\VD/\SG =
      F_1(\conj{B}) = \im(F_1).
    \end{split}
  \end{equation*}
\end{proof}

\begin{lemma}\label{le:morph_lLgpd_imp_semi-local_morph_FBS}
  Let $\PBmap:\PBtotal\rightarrow \PBbase$ be a principal $\SG$-bundle and
  $\jcalU\subset \PBtotal\times \PBtotal$ be a symmetric $pD$-type subset. If
  $s\in \hom_{\lLgpdC_{\PBbase}}(\jcalU'',\jcalU/\SG)$, then
  $s$ is a semi-local morphism between the fiber bundles with a
  section $((\PBbase)\times (\PBbase), \sigma_{(\PBbase)\times (\PBbase)})$
  and $((\PBtotal\times \PBtotal)/\SG, \sigma_{(\PBtotal\times \PBtotal)/\SG)})$ (see
  Example~\ref{ex:elements_of_FBS_Q/G}).
\end{lemma}

\begin{proof}
  It follows readily by noticing that the $\alpha$ and $\epsilon$
  local Lie groupoid structure maps are restrictions and
  co-restrictions of the projection and section maps of the
  corresponding element of $\FBS_{\PBbase}$.
\end{proof}

\begin{proposition}\label{prop:bijection_right_splittings_DAS_flat_DC-lLgpd}
  Let $\PBmap:\PBtotal\rightarrow \PBbase$ be a principal $\SG$-bundle and
  $\jcalU\subset \PBtotal\times \PBtotal$ be a symmetric $D$-type subset. Then, the
  bijection
  $F_{HR}\circ F_{CH}:\Sigma_C(\jcalU)\rightarrow \Sigma_R(\jcalU)$
  determines a bijection between the subset
  $\Sigma_C^e(\jcalU)\subset \Sigma_C(\jcalU)$ of flat discrete
  connections and the subset
  $\ti{\Sigma_R}(\jcalU)\subset \Sigma_R(\jcalU)$ of right splittings
  of extension~\eqref{eq:DAS_in_lLgpd} in the $\lLgpdC_{\PBbase}$
  category.
\end{proposition}

\begin{proof}
  That $\ti{\Sigma_R}(\jcalU)\subset \Sigma_R(\jcalU)$ follows from
  Lemma~\ref{le:morph_lLgpd_imp_semi-local_morph_FBS}. As seen in
  Section~\ref{sec:DAS_in_FBS},
  $\Sigma_R(\jcalU) = F_{HR}(F_{CH}(\Sigma_C(\jcalU)))$; then, it
  follows from
  part~\ref{it:DAS_for_pDtype_U_is_sequence_in_lLgpd-correspondence}
  of Proposition~\ref{prop:DAS_for_pDtype_U_is_sequence_in_lLgpd} that
  $\ti{\Sigma_R}(\jcalU) = F_{HR}(F_{CH}(\Sigma_C^e({\mathcal
    U})))$. Last, being $F_{HR}\circ F_{CH}$ a bijection, its
  restriction and co-restriction is a bijection between
  $\Sigma_C^e(\jcalU)$ and $\ti{\Sigma_R}(\jcalU)$.
\end{proof}

%%%

The following two results, whose proofs we omit, give an alternative
characterization of total intransitivity for local Lie groupoids.

\begin{lemma}\label{le:base_gpd_M_is_initial_in_lLgpd_M}
  Let $M$ be a manifold. Then the base groupoid $M^\dagger$ defined in
  Example~\ref{ex:base_groupoid} is an initial object in $\lLgpdC_M$.
\end{lemma}

\begin{prop}\label{le:tot_disconnected_lLgpd_M_and_hom_to_M_dagger}
  Let $M$ be a manifold and $G\in \ob_{\lLgpdC_M}$. Then
  $\hom_{\lLgpdC_M}(G,M^\dagger)$ is either empty (when $G$ is not
  totally intransitive) or equals $\{\alpha_G\}$ (when $G$ is totally
  intransitive).
\end{prop}

\begin{remark}\label{rem:M_dagger_is_not_final_in_lLgpd}
  Using
  Proposition~\ref{le:tot_disconnected_lLgpd_M_and_hom_to_M_dagger} it
  is not hard to see that, if $M$ is a manifold with more than one
  point, then $M^\dagger$ is not a terminal object in
  $\lLgpdC_M$. Thus, because of
  Lemma~\ref{le:base_gpd_M_is_initial_in_lLgpd_M}, there are no zero
  objects in that category.
\end{remark}

%%%

\begin{remark}\label{rem:categorical_stuff_in_lLgpd}
  Even though the notions of kernel and extension introduced above are
  natural adaptations to the $\lLgpdC_M$ category of standard notions
  in the $\LgpdC_M$ category, it is important to realize that they
  also have natural categorical properties in $\lLgpdC_M$. We begin by
  noticing that, as there is no zero object in $\lLgpdC_M$ (see
  Remark~\ref{rem:M_dagger_is_not_final_in_lLgpd}), the notion of
  categorical kernel discussed in
  Remark~\ref{rem:categorical_stuff_in_FBS} has to be extended. In a
  category ${\mathcal C}$ with an initial object $0$, given
  $f\in\hom_{\mathcal C}(A_1,A_2)$, a triple $(K,j,h)$ is a
  \jdef{categorical kernel} in ${\mathcal C}$ if
  $0\xleftarrow{h} K \xrightarrow{j} A_1$ is a pullback sequence in
  ${\mathcal C}$ of $0\xrightarrow{0^{A_2}} A_2 \xleftarrow{f} A_1$,
  where $0^{A_2}$ is the (unique) initial morphism. In the $\lLgpdC_M$
  category, $M^\dagger$ is an initial object
  (Lemma~\ref{le:base_gpd_M_is_initial_in_lLgpd_M}), so we can apply
  this notion of categorical kernel. Moreover, when
  $K\in\ob_{\lLgpdC_M}$ is totally intransitive, by
  Lemma~\ref{le:tot_disconnected_lLgpd_M_and_hom_to_M_dagger}, there
  is a unique $0_K\in\hom_{\lLgpdC_M}(K,M^\dagger)$, so that we may
  say that $(K,j)$ is a categorical kernel of a morphism $\eta$,
  meaning that $(K,j,0_K)$ is a categorical kernel of $\eta$. Then,
  given $G_1,G_3\in\ob_{\lLgpdC_M}$ with $G_1$ totally intransitive
  and $G_3$ locally trivial, a sequence
  $G_1\xrightarrow{\eta_1} G_2 \xrightarrow{\eta_2} G_3$ in
  $\lLgpdC_M$ is said to be a \jdef{categorical extension} if $\eta_2$
  is a surjective submersion and $(G_1,\eta_1)$ is a categorical
  kernel of $\eta_2$. It can be seen that this notion of categorical
  extension is equivalent to that of extension given in
  Definition~\ref{def:extensionM-lLGpd}. In particular, because of
  Proposition~\ref{prop:DAS_is_extensionM-lLGpd}, the DAS over
  $\jcalU$~\eqref{eq:DAS_in_lLgpd} is a categorical extension in
  $\lLgpdC_{\PBbase}$.
\end{remark}

\begin{remark}\label{rem:CC_and_DC}
  An interesting application of the characterization of flat discrete
  connections as right splittings of the DAS in the \lLgpdC category
  developed in this section is to study the relationship between
  flat discrete connections and flat (continuous) connections on a
  principal bundle. The Lie functor from the \lLgpdC category into the
  category of Lie algebroids
  (see~\cite{bo:mackenzie-general_theory_of_lie_groupoids_and_algebroids})
  maps the DAS of a principal bundle $\PBmap$ to the AS of $\PBmap$;
  in particular, right splittings of the former map to right
  splittings of the latter. In terms of connections this says that
  flat discrete connections $\DC$ on $\PBmap$ are mapped to flat
  connections $\CC$ on $\PBmap$. This is a well known elementary
  process that can be described explicitly as follows: for each
  $v_q\in T_q\PBtotal$, let $q(t)$ be a smooth curve in $\PBtotal$
  integrating $v_q$; then,
  $\CC(v_q)=\frac{d}{dt}\big|_{t=0}\DC(q,q(t))$ (see Section 5.2
  in~\cite{ar:leok_marsden_weinstein-a_discrete_theory_of_connections_on_principal_bundles}). A
  natural question is, given a flat connection $\CC$ on $\PBmap$, find
  all its ``integrals'', that is, all the flat discrete connections
  $\DC$ on $\PBmap$ that are mapped by the previous procedure to
  $\CC$. It can be proved that, essentially, there is a unique such
  $\DC$, as we will discuss elsewhere. On the other hand, for not
  necessarily flat connections (which can be reinterpreted as sections
  of the DAS in the \FBS category using the results of
  Section~\ref{sec:DAS_in_FBS}) we could approach the problem in the
  categories $\FBS$ and that of vector bundles but, so far, our
  analysis is work in progress.
\end{remark}

%%%%%%%%%%%%%%%%%%%%%

\subsection{Semidirect product of local Lie groupoids}
\label{sec:semidirect_product-lLGpd}

We saw in
Section~\ref{sec:right_splittings_and_fiber_product_decomposition} how
right splittings of~\eqref{eq:discrete_atiyah_sequence-def} in the
$\FBS_{\PBbase}$ category are related to isomorphisms in the same
category between~\eqref{eq:discrete_atiyah_sequence-def} and the fiber
product sequence~\eqref{eq:product_sequence_in_FBS-def}. In this
section we explore how right splittings of
sequence~\eqref{eq:DAS_in_lLgpd} in the $\lLgpdC_{\PBbase}$ category
are related to a semidirect product sequence in the same category.

\begin{definition}\label{def:smooth_external_action-lLgpd}
  Let $G$ be a local Lie groupoid and $H$ be a totally intransitive
  Lie groupoid, both over $M$. A \jdef{smooth external action} of $G$
  on $H$ is a smooth map
  $\bullet:G\FP{\beta_G}{\alpha_H} H\rightarrow H$ with the following
  properties (see Definition 2
  in~\cite{ar:metere_montoli-semidirect_products_of_internal_groupoids}). As
  usual, we denote $\bullet(g,h)$ by $g\bullet h$.
  \begin{enumerate}
  \item \label{it:smooth_external_action-lLgpd-alpha}
    $\alpha_H(g\bullet h) = \alpha_G(g)$ for all
    $(g,h)\in G\FP{\beta_G}{\alpha_H} H$,
  \item \label{it:smooth_external_action-lLgpd-assoc}
    $(g_1g_2)\bullet h = g_1\bullet(g_2\bullet h)$, for all
    $(g_1,g_2)\in G_m$ such that
    $(g_2,h)\in G\FP{\beta_G}{\alpha_H}H$,
  \item \label{it:smooth_external_action-lLgpd-dist}
    $g\bullet (h_1h_2) = (g\bullet h_1)(g \bullet h_2)$, for all
    $(h_1,h_2)\in H_2$ such that
    $(g,h_1)\in G\FP{\beta_G}{\alpha_H}H$,
  \item \label{it:smooth_external_action-lLgpd-identity}
    $\epsilon_G(\alpha_H(h)) \bullet h = h$, for all $h\in H$.
  \end{enumerate}
\end{definition}

\begin{lemma}\label{le::simple_properties_of_actions}
  Let $\bullet:G\FP{\beta_G}{\alpha_H} H\rightarrow H$ be a smooth
  external action of $G$ on $H$. Then,
  \begin{gather}
    g\bullet \epsilon_H(\beta_G(g)) = \epsilon_H(\alpha_G(g))
    \stext{ for all } g\in G, \stext{ and }
    \label{eq:simple_properties_of_actions-identity}\\
    (g\bullet h)^{-1} = g\bullet h^{-1} \stext{ for all } (g,h)\in
    G\FP{\beta_G}{\alpha_H}
    H \label{eq:simple_properties_of_actions-inverses}.
  \end{gather}
\end{lemma}

\begin{proof}
  \eqref{eq:simple_properties_of_actions-identity}
  and~\eqref{eq:simple_properties_of_actions-inverses} are proved
  using points~\ref{it:uniqueness_in_loc_lie_groupoid-left_unit}
  and~\ref{it:uniqueness_in_loc_lie_groupoid-left_identity} of
  Proposition~\ref{prop:uniqueness_in_loc_lie_groupoid} respectively.
\end{proof}

\begin{lemma}\label{le:construction_of_semidirect_product_lLgpoid}
  Let $\bullet$ be a smooth external action of the local Lie
  groupoid $G\gpoidarrows M$ on the totally intransitive Lie groupoid
  $H\gpoidarrows M$. Define $P:= H\FP{\beta_H}{\alpha_G} G$ and the
  maps
  \begin{gather*}
    \alpha_P,\beta_P:P\rightarrow M,\stext{ by }
    \alpha_P(h,g) := \alpha_G(g) \stext{ and } \beta_P(h,g) := \beta_G(g)\\
    m_P:P_m\rightarrow P \stext{ by }
    m_P((h_1,g_1),(h_2,g_2)):=(h_1(g_1\bullet h_2),g_1g_2)\\
    \text{ where } P_m:=\{((h_1,g_1),(h_2,g_2))\in P^2:(g_1,g_2)\in G_m\},\\
    \epsilon_P:M\rightarrow P\stext{ by } \epsilon_P(m) :=
    (\epsilon_H(m), \epsilon_G(m)),\\
    i_P:P\rightarrow P \stext{ by } i_P(h,g):= ((g^{-1}\bullet
    h)^{-1},g^{-1}).
  \end{gather*}
  Then, $(P,\alpha_P,\beta_P,m_P,\epsilon_P,i_P)$ is a local Lie
  groupoid over $M$.
\end{lemma}

\begin{proof}
  It is straightforward but lengthy, checking that the structure maps
  are well defined and that all conditions appearing in
  Definition~\ref{def:local_lie_groupoid} are met.
\end{proof}

Notice that the smooth external action $\bullet$ can be fully
recovered from the local multiplication $m_P$ of $P$.

\begin{definition}
  Given a smooth external action $\bullet$ of the local Lie groupoid
  $G\gpoidarrows M$ on the totally intransitive Lie groupoid
  $H\gpoidarrows M$ the local Lie groupoid $P$ defined in
  Lemma~\ref{le:construction_of_semidirect_product_lLgpoid} is called
  the \jdef{semidirect product local Lie groupoid} defined by
  $\bullet$ and is denoted by $H \rtimes G \gpoidarrows M$.
\end{definition}

\begin{proposition}\label{prop:semidirect_product_imp_extension}
  Let $\bullet$ be a smooth external action of the local Lie
  groupoid $G\gpoidarrows M$ on the totally intransitive Lie groupoid
  $H\gpoidarrows M$. Then, the diagram
  \begin{equation}
    \label{eq:semidirect_product_imp_extension-semidirect_extension}
    \xymatrix{{H} \ar[r]^-{j_{\rtimes}} & {H\rtimes G}
    \ar[r]^-{\rho_{\rtimes}} & {G}}
  \end{equation}
  is an extension of $G$ by $H$ where $H\rtimes G$ is the semidirect
  product of $H$ by $G$ associated to $\bullet$,
  $j_{\rtimes}(h) := (h,\epsilon_G(\beta_H(h)))$ and
  $\rho_{\rtimes}(h,g):= g$. In addition, $\rho_\rtimes$ is a
  submersion and
  $(H\rtimes G)_m = (\rho_\rtimes \times \rho_\rtimes)^{-1}(G_m)$.
  Furthermore, $s_{\rtimes}:G\rightarrow H\rtimes G$ defined by
  $s_{\rtimes}(g):=(\epsilon_H(\alpha_G(g)),g)$ is a (right) splitting
  of the semidirect product extension.
\end{proposition}

The proof of Proposition~\ref{prop:semidirect_product_imp_extension}
will follow after some preliminary work.

\begin{lemma}\label{le:semidirect_product_imp_extension-morphisms}
  With the same hypotheses of
  Proposition~\ref{prop:semidirect_product_imp_extension}, the maps
  $j_\rtimes$, $\rho_\rtimes$ and $s_\rtimes$ are morphisms in the
  $\lLgpdC_{M}$ category, $j_\rtimes$ is an embedding and $s_\rtimes$
  is a right inverse of $\rho_\rtimes$.
\end{lemma}

\begin{proof}
  We only check that $j_\rtimes\in\hom_{\lLgpdC_M}(H,H\rtimes G)$
  next; the cases for $\rho_\rtimes$ and $s_\rtimes$ being similar.
  As $\beta_H(h) = \alpha_G(\epsilon_G(\beta_H(h)))$, we have that
  $j_\rtimes$ is a well defined map which, being the co-restriction to
  $P:=H\rtimes G$ of a smooth map into $H\times G$ and whose image is
  contained in the embedded submanifold $P\subset H\times G$, is
  smooth.  For $m\in M$ we have
  $j_\rtimes(\epsilon_H(m)) =
  (\epsilon_H(m),\epsilon_G(\beta_H(\epsilon_H(m)))) = (\epsilon_H(m),
  \epsilon_G(m)) = \epsilon_P(m)$, so that
  $j_\rtimes(\epsilon_H(M)) \subset \epsilon_P(M)$.  If
  $(h_1,h_2)\in H_m=H_2$, we have that
  $\beta_H(h_1) = \alpha_H(h_2) =\beta_H(h_2)$. Then
  $(j_\rtimes(h_1), j_\rtimes(h_2)) = ((h_1,\epsilon_G(\beta_H(h_1))),
  (h_2,\epsilon_G(\beta_H(h_2)))) \in P_m$ because
  \begin{equation*}
    (\epsilon_G(\beta_H(h_1)), \epsilon_G(\beta_H(h_2))) =
    (\epsilon_G(\beta_H(h_1)), \epsilon_G(\beta_H(h_1))) \in G_m.
  \end{equation*}
  In this case
  \begin{equation*}
    \begin{split}
      j_\rtimes(h_1h_2) =& (h_1h_2, \epsilon_G(\beta_H(h_1h_2))) =
      (h_1h_2, \epsilon_G(\beta_H(h_2))) \\=&
      (h_1(\epsilon_G(\alpha_H(h_2))\bullet h_2),
      \epsilon_G(\beta_H(h_2))\epsilon_G(\beta_H(h_2))) \\=&
      (h_1(\epsilon_G(\beta_H(h_1)) \bullet
      h_2),\epsilon_G(\beta_H(h_1)) \epsilon_G(\beta_H(h_2)))\\=&
      (h_1,\epsilon_G(\beta_H(h_1))) (h_2,\epsilon_G(\beta_H(h_2))) =
      j_\rtimes(h_1) j_\rtimes(h_2),
    \end{split}
  \end{equation*}
  and, so, $j_\rtimes\in\hom_{\lLgpdC}(H,P)$. A direct computation
  shows that $(j_\rtimes)_0 = id_M$, thus
  $j_\rtimes \in \hom_{\lLgpdC_M}(H,P)$.

  Notice that if $p_1^r:H\rtimes G\rightarrow H$ is the (restriction
  of the) projection onto the first factor, it is a smooth map and
  $p_1^r\circ j_\rtimes = id_H$. Then, by
  Lemma~\ref{le:section_of_smooth_map_is_embedding}, $j_\rtimes$ is an
  embedding.

  That $s_\rtimes$ is a right inverse of $\rho_\rtimes$ follows
  immediately by direct computation.
\end{proof}

\begin{proof}[Proof (of
  Proposition~\ref{prop:semidirect_product_imp_extension})]
  By Lemmas~\ref{le:construction_of_semidirect_product_lLgpoid}
  and~\ref{le:semidirect_product_imp_extension-morphisms}
  diagram~\eqref{eq:semidirect_product_imp_extension-semidirect_extension}
  is a sequence in $\lLgpdC_M$ and
  $\rho_\rtimes\circ s_\rtimes = id_G$. Because of the last identity
  $\rho_\rtimes$ is a surjective submersion and, also, once we have
  proved
  that~\eqref{eq:semidirect_product_imp_extension-semidirect_extension}
  is an extension,
  that~\eqref{eq:semidirect_product_imp_extension-semidirect_extension}
  is split.

  As seen in
  Lemma~\ref{le:semidirect_product_imp_extension-morphisms},
  $j_\rtimes$ is an embedding and $\rho_\rtimes$ is a
  submersion. Also, as $H$ is a Lie groupoid,
  $(j_\rtimes,j_\rtimes)^{-1}((G_2)_m) \subset (G_1)_m$ as observed in
  Remark~\ref{rem:preimage_of_multiplicable_elements_in_Lgpd}. On the
  other hand, as
  $\rho_\rtimes \circ j_\rtimes = \epsilon_G\circ \beta_H$, we have
  that $\im(j_\rtimes)\subset \ker(\rho_\rtimes)$ while, if
  $(h,g)\in \ker(\rho_\rtimes)$, we see that
  $g=\rho_\rtimes(h,g) = \epsilon_G(\alpha_G(g)) =
  \epsilon_G(\beta_H(h))$ and, then,
  $(h,g)=j_\rtimes(h)\in\im(j_\rtimes)$, proving that
  $\im(j_\rtimes)=\ker(\rho_\rtimes)$.  Hence, the
  sequence~\eqref{eq:semidirect_product_imp_extension-semidirect_extension}
  is an extension in $\lLgpdC_M$.

  The only thing left to prove is that
  $(H\rtimes G)_m = (\rho_\rtimes \times \rho_\rtimes)^{-1}(G_m)$,
  which follows by simple inspection of the two sets..
\end{proof}

Sequence~\eqref{eq:semidirect_product_imp_extension-semidirect_extension}
will be called the \jdef{semidirect product sequence} in the
$\lLgpdC_M$ category.

\begin{proposition}\label{prop:split_seq=>semidirec_seq}
  Let $H \xrightarrow{j} E \xrightarrow{\rho} G$ be an extension in
  $\lLgpdC_M$ that is right split by $s$. If $H$ is a Lie groupoid and
  $E_m = (\rho\times \rho)^{-1}(G_m)$, then
  \begin{equation}\label{eq:split_seq=>semidirec_seq-bullet_def}
    \bullet:G\FP{\beta_G}{\alpha_H} H\rightarrow H \stext{ given by }
    g\bullet h := j^{-1}(s(g) j(h) s(g)^{-1}) 
  \end{equation}
  is a smooth external action of $G$ on $H$ and there is an
  isomorphism $\Phi\in\hom_{\lLgpdC_M}(E,H\rtimes G)$ such that the
  following diagram in $\lLgpdC_M$ is commutative.
  \begin{equation}\label{eq:split_seq=>semidirec_seq-seq_diag}
    \xymatrix{
      {H} \ar[r]^-{j_\rtimes} & {H\rtimes G} \ar[r]^-{\rho_\rtimes} &
      {G} \\
      {H} \ar[u]^{id_H} \ar[r]_j & {E} \ar[u]_{\Phi}
      \ar[r]_{\rho} &  {G} \ar[u]_{id_G}
    }
  \end{equation}
\end{proposition}

\begin{lemma}\label{le:split_seq=>semidirec_seq-dot_well_defined}
  With the hypotheses as in
  Proposition~\ref{prop:split_seq=>semidirec_seq},
  $\bullet:G\FP{\beta_G}{\alpha_H} H\rightarrow H$ is a smooth
  external action.
\end{lemma}

\begin{proof}
  At the level of sets, $\bullet$ is well defined. Indeed, for
  $(g,h)\in G\FP{\beta_G}{\alpha_H} H$, we have
  $\beta_G(g) = \alpha_H(h)$. Then, as
  \begin{equation*}
    \begin{split}
      (\rho\times \rho)(s(g),j(h)) =& (\rho(s(g)),\rho(j(h))) =
      (g,\epsilon_G(\alpha_H(h))) = (g,\epsilon_G(\beta_G(g))) \in G_m
    \end{split}
  \end{equation*}
  and as, by hypothesis, $E_m = (\rho\times \rho)^{-1}(G_m)$, we
  see that $(s(g),j(h))\in E_m$. Similarly, as
  \begin{equation*}
    \begin{split}
      (\rho\times \rho)(s(g)j(h), s(g^{-1})) =&
      (\rho(s(g)j(h)),\rho(s(g^{-1}))) = (g\epsilon_G(\alpha_H(h)),g^{-1}) \\=&
      (g\epsilon_G(\beta_G(g)),g^{-1}) = (g,g^{-1}) \in G_m,
    \end{split}
  \end{equation*}
  we have that $(s(g)j(h), s(g^{-1}))\in E_m$. Furthermore,
  \begin{equation*}
    \begin{split}
      \rho((s(g)j(h))s(g^{-1})) =& \rho(s(g)j(h))\rho(s(g^{-1})) = (g
      \rho(j(h))) g^{-1} = (g \epsilon_G(\alpha_H(h)))g^{-1} \\=& (g
      \epsilon_G(\beta_G(g)))g^{-1} =g g^{-1} =
      \epsilon_G(\alpha_G(g)),
    \end{split}
  \end{equation*}
  proves that
  $s(g)j(h)s(g^{-1}) \in \rho^{-1}(\epsilon_G(M)) = \ker(\rho) =
  \im(j)$. As $j|^{\im(j)}:H\rightarrow \im(j)$ is a diffeomorphism, we
  see that $g\bullet h = (j|^{\im(j)})^{-1}(s(g)j(h) s(g^{-1}))$ is a
  well defined map. In addition, as
  $(g,h)\rightarrow s(g)j(h)s(g^{-1})$ is a smooth map from
  $G\FP{\beta_G}{\alpha_H} H$ into $E$ whose image is contained in the
  embedded submanifold $\im(j)$ and $(j|^{\im(j)})^{-1}$ is a
  diffeomorphism, we conclude that $\bullet$ is a smooth map.

  Checking that $\bullet$ satisfies the remaining conditions to be a
  smooth external product follows easily by direct computation,
  recalling that $j$ is one-to-one.
\end{proof}

\begin{proof}[Proof (of
  Proposition~\ref{prop:split_seq=>semidirec_seq})]
  We know from
  Lemma~\ref{le:split_seq=>semidirec_seq-dot_well_defined} that
  $\bullet$ is a smooth external action of $G$ on $H$. Next, we want
  to define the map $\Phi:E\rightarrow H\rtimes G$.
  
  For any $e\in E$, we have that
  \begin{equation*}
    (e,s(\rho(e^{-1}))) \in E_m \iff (\rho(e),\rho(s(\rho(e^{-1})))) \in G_m
  \end{equation*}
  but, as
  $(\rho(e),\rho(s(\rho(e^{-1})))) = (\rho(e),\rho(e^{-1})) =
  (\rho(e), \rho(e)^{-1}) \in G_m$, we conclude that
  $(e,s(\rho(e^{-1}))) \in E_m$. Then, the map
  $\conj{\Phi}_1:E\rightarrow E$ defined by
  $\conj{\Phi}_1(e) := e\, s(\rho(e^{-1}))$ is well defined and
  smooth. In addition, as
  $\rho(\conj{\Phi}_1(e)) = \rho(e\, s(\rho(e^{-1}))) =
  \rho(e)\rho(s(\rho(e^{-1}))) = \rho(e) \rho(e)^{-1} =
  \epsilon_G(\alpha_G(\rho(e))) = \epsilon_G(\alpha_E(e))$, we see
  that the image of $\conj{\Phi}_1$ is contained in the embedded
  submanifold $\ker(\rho)=\im(j)\subset E$, so that
  $\Phi_1:E\rightarrow H$ defined by
  $\Phi_1(e) := {j_c}^{-1}(\conj{\Phi}_1(e))$, where
  $j_c:=j|^{\im(j)}$, is a smooth map. We define
  $\conj{\Phi}:E\rightarrow H\times G$ by
  \begin{equation*}
    \conj{\Phi}(e) := (\Phi_1(e), \rho(e)) =
    ({j_c}^{-1}(e\, s(\rho(e^{-1}))), \rho(e)),
  \end{equation*}
  that clearly is a smooth map. As
  $\beta_H({j_c}^{-1}(e\, s(\rho(e^{-1})))) = \beta_E(s(\rho(e^{-1}))) =
  \beta_G(\rho(e^{-1})) = \alpha_G(\rho(e))$, we see that the image of
  $\conj{\Phi}$ is contained in the embedded submanifold
  $H\rtimes G\subset H\times G$. Hence, its co-restriction
  $\Phi:E\rightarrow H\rtimes G$ defined by
  $\Phi:=\conj{\Phi}|^{H\rtimes G}$ is a smooth map.

  The rest of the proof follows from
  Lemma~\ref{le:split_seq=>semidirec_seq-commutative}, where
  diagram~\eqref{eq:split_seq=>semidirec_seq-seq_diag} is seen to be
  commutative in the $\lLgpdC_M$ category, and
  Lemma~\ref{le:split_seq=>semidirec_seq-isomorphism}, where $\Phi$ is
  seen to be an isomorphism in the same category.
\end{proof}

\begin{lemma}\label{le:split_seq=>semidirec_seq-commutative}
  With the hypotheses as in
  Proposition~\ref{prop:split_seq=>semidirec_seq} and $\Phi$ defined
  in the proof of that result,
  diagram~\eqref{eq:split_seq=>semidirec_seq-seq_diag} is commutative
  in the $\lLgpdC_M$ category.
\end{lemma}

\begin{proof}
  Let $P:=H\rtimes G$. A direct computation shows that
  $\Phi\circ \epsilon_E = \epsilon_P$, so that
  $\Phi(\epsilon_E(M))\subset \epsilon_P(M)$.

  For any $(e_1,e_2)\in E_m = (\rho\times \rho)^{-1}(G_m)$ we have
  that $(\rho(e_1),\rho(e_2)) \in G_m$. Then, as
  \begin{equation*}
    (\rho_{\rtimes}(\Phi(e_1)), \rho_\rtimes(\Phi(e_2))) =
    (\rho(e_1),\rho(e_2)) \in G_m, 
  \end{equation*}
  we have that
  $(\Phi(e_1),\Phi(e_2)) \in (\rho_\rtimes \times
  \rho_\rtimes)^{-1}(G_m) = P_m$ (by
  Proposition~\ref{prop:semidirect_product_imp_extension}). In this
  case, using that $j_c$ is one-to-one and that
  $\ker(\rho)=\im(\rho)$, it is straightforward to check that
  $\Phi(e_1e_2) = \Phi(e_1)\Phi(e_2)$, so that
  $\Phi\in\hom_{\lLgpdC}(E,P)$. Direct computations show that
  \begin{equation}
    \label{eq:split_seq=>semidirec_seq-comm_diag}
    \Phi_0 = id_M,\quad 
    \Phi\circ j = j_\rtimes \stext{ and }
    \rho_\rtimes \circ \Phi =\rho,
  \end{equation}
  so that $\Phi \in\hom_{\lLgpdC_M}(E,P)$ and
  diagram~\eqref{eq:split_seq=>semidirec_seq-seq_diag} is commutative.
\end{proof}

\begin{lemma}\label{le:split_seq=>semidirec_seq-isomorphism}
  With the hypotheses of
  Proposition~\ref{prop:split_seq=>semidirec_seq} and $\Phi$ defined
  in the proof of that result, $\Phi$ is an isomorphism in the
  $\lLgpdC_M$ category.
\end{lemma}

\begin{proof}
  Let $P:=H\rtimes G$. Consider now the smooth maps
  $p_1^r:P\rightarrow H$ and $p_2^r:P\rightarrow G$ obtained as
  restrictions to $P$ of the canonical projections onto the
  corresponding factor. Then, the maps $f_1,f_2:P\rightarrow E$
  defined by $f_1(p):=j_c(p_1^r(p))$ and $f_2(p):=s(p_2^r(p))$ are
  smooth. Furthermore, for $(h,g)\in P$, that is, satisfying
  $\beta_H(h) = \alpha_G(g)$, we have that
  \begin{equation*}
    \begin{split}
      (\rho(f_1(h,g)),\rho(f_2(h,g))) =& (\rho(j_c(h)),\rho(s(g))) =
      (\epsilon_G(\alpha_H(h)), g) \\=& (\epsilon_G(\beta_H(h)), g) =
      (\epsilon_G(\alpha_G(g)), g) \in G_m,
    \end{split}
  \end{equation*}
  so that $(f_1(h,g),f_2(h,g)) \in (\rho\times\rho)^{-1}(G_m) =
  E_m$. In this case, the map $\Psi:P\rightarrow E$ defined by
  $\Psi(p):= f_1(p)f_2(p)$ is well defined and, as
  that $m_E$ is smooth, $\Psi$ is smooth.

  Evaluation and a short computation show that
  $\Psi\circ \Phi = id_E$, while a lengthier one shows that
  $\Phi\circ \Psi = id_P$.  Thus, $\Psi = \Phi^{-1}$ and we have
  proved that $\Phi$ is a diffeomorphism.

  Let $(e_1,e_2) \in (\Phi\times \Phi)^{-1}(P_m)$, so that
  $(\Phi(e_1),\Phi(e_2)) \in P_m$. As $\rho_\rtimes$ is a morphism in
  $\lLgpdC_M$, we have that
  $(\rho_\rtimes(\Phi(e_1)), \rho_\rtimes(\Phi(e_2))) \in G_m$. But,
  then, by~\eqref{eq:split_seq=>semidirec_seq-comm_diag},
  $(\rho(e_1),\rho(e_2)) \in G_m$ and, as
  $E_m = (\rho\times \rho)^{-1}(G_m)$, we see that $(e_1,e_2)\in
  E_m$. All together, we proved that
  $(\Phi\rtimes \Phi)^{-1}(P_m)\subset E_m$.  Finally, by
  Lemma~\ref{le:morph_and_diffeo_imp_iso_in_lLgpd}, $\Phi$ is an
  isomorphism in the $\lLgpdC_M$ category.
\end{proof}

\begin{remark}
  All smooth external actions $\bullet$ of $G$ on $H$ come from split
  extensions via the construction described in
  Proposition~\ref{prop:split_seq=>semidirec_seq}. Indeed, we have the
  easily verified identity
  \begin{equation*}
    j_\rtimes(g\bullet h) = s_\rtimes(g) j_{\rtimes}(h) s_\rtimes(g)^{-1}
    \stext{ for all } (g,h) \in G\FP{\beta_G}{\alpha_H} H
  \end{equation*}
  and, then, we use the construction of
  Proposition~\ref{prop:split_seq=>semidirec_seq} applied to the
  semidirect product
  extension~\eqref{eq:semidirect_product_imp_extension-semidirect_extension}
  with $\Phi=id_{H\rtimes G}$ (recall
  that~\eqref{eq:semidirect_product_imp_extension-semidirect_extension}
  is split by
  Proposition~\ref{prop:semidirect_product_imp_extension}).
\end{remark}

Let $H$ be a totally intransitive Lie groupoid and
\begin{equation}
  \label{eq:an_extension_in_lLgpdC}
  H\xrightarrow{j} E\xrightarrow{\rho} G
\end{equation}
be an extension in the
$\lLgpdC_M$ category such that
\begin{equation}
  \label{eq:an_extension_in_lLgpdC-controm_of_E_m_by_G_m}
  E_m=(\rho\times\rho)^{-1}(G_m).
\end{equation}
We define $\ti{\Sigma_I}$ as the set of all pairs $(\bullet,\Phi)$
where $\bullet$ is a smooth external action of $G$ on $H$ and
$\Phi\in\hom_{\lLgpdC_M}(E,H\rtimes G)$ is an isomorphism where
$H\rtimes G$ is constructed using $\bullet$, such that
diagram~\eqref{eq:split_seq=>semidirec_seq-seq_diag} is
commutative. We also define $\ti{\Sigma_R}$ as the set of right
splittings of $H\xrightarrow{j} E\xrightarrow{\rho} G$ in the
$\lLgpdC_M$ category. There is a map
$\ti{F}_{RI}:\ti{\Sigma_R}\rightarrow \ti{\Sigma_I}$ defined by
$\ti{F}_{RI}(s) := (\bullet,\Phi)$, both provided by
Proposition~\ref{prop:split_seq=>semidirec_seq}. Conversely, define
$\ti{F}_{IR}:\ti{\Sigma_I}\rightarrow \ti{\Sigma_R}$ by
$\ti{F}_{IR}(\bullet,\Phi):=s$, for $s:=\Phi^{-1}\circ s_\rtimes$,
where $s_\rtimes$ is defined in
Proposition~\ref{prop:semidirect_product_imp_extension}, using
$\bullet$.

\begin{theorem}\label{prop:equiv_split_extensions_and_semi_direct_products}
  With the definitions as above, both $\ti{F}_{RI}$ and $\ti{F}_{IR}$
  are well defined maps that are inverses of each other. In other
  words, in the $\lLgpdC_M$ category, there is a bijective
  correspondence between the right splittings of an
  extension~\eqref{eq:an_extension_in_lLgpdC}
  satisfying~\eqref{eq:an_extension_in_lLgpdC-controm_of_E_m_by_G_m}
  and the isomorphisms of this extension with semidirect product
  sequences making~\eqref{eq:split_seq=>semidirec_seq-seq_diag}
  commutative.
\end{theorem}

\begin{proof}
  That $\ti{F}_{RI}$ is well defined is the content of
  Proposition~\ref{prop:split_seq=>semidirec_seq}. If
  $(\bullet,\Phi)\in\ti{\Sigma_I}$, then
  $\Phi\in\hom_{\lLgpdC_M}(E,H\rtimes G)$ is an isomorphism; thus,
  $\Phi^{-1} \in \hom_{\lLgpdC_M}(H\rtimes G,E)$ and, as
  $s_\rtimes\in\hom_{\lLgpdC_M}(G,H\rtimes G)$ by
  Proposition~\ref{prop:semidirect_product_imp_extension}, then
  $s := \Phi^{-1} \circ s_\rtimes \in \hom_{\lLgpdC_M}(G,E)$ and
  $\rho \circ s = \rho \circ \Phi^{-1} \circ s_\rtimes =
  \rho_\rtimes\circ \Phi \circ \Phi^{-1} \circ s_\rtimes = id_G$, so
  that $s\in \ti{\Sigma_R}$, and $\ti{F}_{IR}$ is well defined.

  For $s\in\ti{\Sigma_R}$, let $(\bullet,\Phi):=\ti{F}_{RI}(s)$ and
  $\conj{s}:=\ti{F}_{IR}(\bullet,\Phi)$. Then, for $g\in G$, recalling the
  expression for $\Phi^{-1} = \Psi$ given in the proof of
  Lemma~\ref{le:split_seq=>semidirec_seq-isomorphism},
  \begin{equation*}
    \begin{split}
      \conj{s}(g) =& \Phi^{-1}(s_\rtimes(g)) =
      \Phi^{-1}(\epsilon_H(\alpha_G(g)),g) =
      j_c(\epsilon_H(\alpha_G(g))) s(g) = \epsilon_K(\alpha_G(g))
      s(g) \\=& \epsilon_E(\alpha_E(s(g))) s(g) = s(g),
    \end{split}
  \end{equation*}
  and we see that
  $\ti{F}_{IR} \circ \ti{F}_{RI} = id_{\ti{\Sigma_R}}$.

  For $(\bullet,\Phi)\in\ti{\Sigma_I}$, let
  $s:=\ti{F}_{IR}(\bullet,\Phi)$ and
  $(\conj{\bullet},\conj{\Phi}):=\ti{F}_{RI}(s)$. By construction,
  $\conj{\bullet}$
  satisfies~\eqref{eq:split_seq=>semidirec_seq-bullet_def} and, then,
  a direct computation shows that
  $(\Phi \circ j)(g\ \conj{\bullet}\ h) = (\Phi \circ j)(g\bullet h)$
  for all $(g,h)\in G\FP{\beta_G}{\alpha_H} H$, implying that
  $\conj{\bullet} = \bullet$.

  If $p_j^r$ is the restriction of the projection onto the
  corresponding component, as both $\conj{\Phi}$ and $\Phi$ make
  diagram~\eqref{eq:split_seq=>semidirec_seq-seq_diag} commutative, we
  have that $p_2^r\circ \Phi= \rho= p_2^r\circ \conj{\Phi}$.

  For $e\in E$, as seen in the proof of
  Proposition~\ref{prop:split_seq=>semidirec_seq}, we have
  \begin{equation}\label{eq:equiv_split_extensions_and_semi_direct_products}
    \begin{split}
      (p_1^r\circ \conj{\Phi})(e) =& {j_c}^{-1}(e\, s(\rho(e^{-1}))) =
      {j_c}^{-1}(e\, \Phi^{-1}(s_\rtimes(\rho(e^{-1})))).
    \end{split}
  \end{equation}
  A straightforward computation now shows that
  $e\, \Phi^{-1}(s_\rtimes(\rho(e^{-1}))) =
  \Phi^{-1}(j_\rtimes(p_1^r(\Phi(e))))$.  Back
  in~\eqref{eq:equiv_split_extensions_and_semi_direct_products}, and
  noticing that $\Phi\circ j_c = j_\rtimes$, we have
  \begin{equation*}
    \begin{split}
      p_1^r\circ \conj{\Phi} =& {j_c}^{-1}\circ \Phi^{-1}\circ
      j_\rtimes\circ (p_1^r \circ \Phi) = p_1^r \circ \Phi.
    \end{split}
  \end{equation*}
  All together, we conclude that $\conj{\Phi}=\Phi$ and, then,
  $\ti{F}_{RI}\circ \ti{F}_{IR} = id_{\ti{\Sigma_I}}$, so that
  $\ti{F}_{RI}$ and $\ti{F}_{IR}$ are mutually inverse maps.
\end{proof}

Last, we specialize
Theorem~\ref{prop:equiv_split_extensions_and_semi_direct_products} to
the case of discrete Atiyah sequences.

\begin{corollary}\label{cor:flat_DC_right)splittings_and_semidirect_product}
  Let $\PBmap:\PBtotal\rightarrow \PBbase$ be a principal $\SG$-bundle and
  $\jcalU\subset \PBtotal\times \PBtotal$ be a symmetric $D$-type
  subset. Then, there are bijective correspondences between
  \begin{enumerate}
  \item\label{it:flat_DC_right)splittings_and_semidirect_product-DC}
    $\Sigma_C^e(\jcalU)$, the set of flat discrete connections
    on $\PBmap$ with domain $\jcalU$,
  \item\label{it:flat_DC_right)splittings_and_semidirect_product-split}
    $\ti{\Sigma_R}(\jcalU)$, the set of right splittings of the
    discrete Atiyah sequence over
    $\jcalU$~\eqref{eq:DAS_in_lLgpd} in the $\lLgpdC_{\PBbase}$
    category, and
  \item\label{it:flat_DC_right)splittings_and_semidirect_product-SDP}
    $\ti{\Sigma_I}$, the set of isomorphisms from the discrete Atiyah
    sequence over $\jcalU$~\eqref{eq:DAS_in_lLgpd} to any
    semidirect product extension of $\jcalU''$ by $\conjB$.
  \end{enumerate}
\end{corollary}

\begin{proof}
  The equivalence
  between~\ref{it:flat_DC_right)splittings_and_semidirect_product-DC}
  and~\ref{it:flat_DC_right)splittings_and_semidirect_product-split}
  is given by
  Proposition~\ref{prop:bijection_right_splittings_DAS_flat_DC-lLgpd},
  while the equivalence
  between~\ref{it:flat_DC_right)splittings_and_semidirect_product-split}
  and~\ref{it:flat_DC_right)splittings_and_semidirect_product-SDP}
  if given by
  Theorem~\ref{prop:equiv_split_extensions_and_semi_direct_products},
  applied to the extension~\eqref{eq:DAS_in_lLgpd}.
\end{proof}

%%%%%%%%%%%%%%%%%%%%%%%%%%%%%%%%%%%%%%%%%%%%%%%%%%%%%%%%%%

\appendix

\section{Some basic and convenient properties}
\label{sec:some_basic_and_auxiliary_properties}

Here we include some basic results that we need for the main
presentation. In most cases, they are analogues of well known results
that we need in the \emph{local} groupoid context or that we state
explicitly because we are using the conventions
of~\cite{ar:marrero_martin_martinez-discrete_lagrangian_and_hamiltonian_mechanics_on_lie_groupoids}
instead of those
of~\cite{bo:mackenzie-general_theory_of_lie_groupoids_and_algebroids}. Most
of the proofs will be omitted because they are simple exercises or
adaptations of their Lie groupoid version.

\begin{lemma}\label{le:lLgpd_basic_properties}
  Let $G$ be a local Lie groupoid over $M$. Then, the following
  assertions are true.
  \begin{enumerate}
  \item \label{it:lLgpd_basic_properties-alpha_beta_onto}
    $\alpha,\beta:G\rightarrow M$ are onto.
  \item \label{it:lLgpd_basic_properties-alpha_beta_inverses} For
    $g\in G$, $\alpha(g^{-1})=\beta(g)$ and
    $\beta(g^{-1}) = \alpha(g)$.
  \item \label{it:lLgpd_basic_properties-alpha_beta_m} If
    $(g_1,g_2)\in G_m$, then $\alpha(g_1g_2)=\alpha(g_1)$ and
    $\beta(g_1g_2)=\beta(g_2)$.
  \item \label{it:lLgpd_basic_properties-epsilon_embedding}
    $\epsilon:M\rightarrow G$ is an embedding and
    $\epsilon(M)\subset G$ is an embedded submanifold.
  \item \label{it:lLgpd_basic_properties-epsilon(M)_2}
    $\epsilon(M)_2:=\{(\epsilon_1,\epsilon_2)\in
    \epsilon(M)^2:\beta(\epsilon_1) = \alpha(\epsilon_2)\}\subset
    G_m$.
  \end{enumerate}
\end{lemma}

The following result is a revised version of Proposition 1.1.2
in~\cite{bo:mackenzie-general_theory_of_lie_groupoids_and_algebroids}
adapted to the local context.

\begin{proposition}\label{prop:uniqueness_in_loc_lie_groupoid}
  Let $G\gpoidarrows M$ be a local Lie groupoid and $g\in G$. Then,
  the following assertions are true.
  \begin{enumerate}
  \item \label{it:uniqueness_in_loc_lie_groupoid-left_unit} If
    $h\in G$ is such that $(h,g)\in G_m$ and $hg=g$, then
    $h=\epsilon(\alpha(g))$.
  \item \label{it:uniqueness_in_loc_lie_groupoid-right_unit} If
    $h\in G$ is such that $(g,h)\in G_m$ and $gh=g$, then
    $h=\epsilon(\beta(g))$.
  \item \label{it:uniqueness_in_loc_lie_groupoid-left_identity} If
    $h\in G$ is such that $(h,g)\in G_m$ and $hg=\epsilon(\beta(g))$,
    then $h=g^{-1}$.
  \item \label{it:uniqueness_in_loc_lie_groupoid-right_identity} If
    $h\in G$ is such that $(g,h)\in G_m$ and $gh=\epsilon(\alpha(g))$,
    then $h=g^{-1}$.
  \end{enumerate}
\end{proposition}

\begin{lemma}\label{le:morphism_of_loc_lie_groupoid-props}
  Let $F:G\rightarrow G'$ be a morphism of the local Lie groupoids
  $G\gpoidarrows M$ and $G'\gpoidarrows M'$. Then, the following
  assertions are true
  \begin{enumerate}
  \item \label{it:morphism_of_loc_lie_groupoids-F_of_identity}
    $F\circ \epsilon\circ \alpha = \epsilon'\circ \alpha'\circ F$ and
    $F\circ \epsilon\circ \beta = \epsilon'\circ \beta'\circ F$.
  \item \label{it:morphism_of_loc_lie_groupoids-F_of_inverse}
    $F(g^{-1}) = F(g)^{-1}$ for all $g\in G$.
  \item \label{it:morphism_of_loc_lie_groupoids-F_0} The map
    \begin{equation}
      \label{eq:F_0_from_F-def}
      F_0:M\rightarrow M' \stext{ defined by }
      F_0:=\alpha'\circ F\circ \epsilon
    \end{equation}
    is the unique smooth map between those spaces satisfying
    \begin{equation}\label{eq:morphism_of_loc_lie_groupoid-formulas}
      \alpha'\circ F = F_0\circ \alpha,\quad \beta'\circ F = F_0\circ \beta 
      \stext{ and } F\circ \epsilon = \epsilon'\circ F_0.
    \end{equation}
  \end{enumerate}
\end{lemma}

\begin{proof}
  \begin{enumerate}
  \item For any $g\in G$, $(g,\epsilon(\beta(g))) \in G_m$, so that
    $(F(g),F(\epsilon(\beta(g))))\in G'_m$ and
    $F(g) = F(g\epsilon(\beta(g))) = F(g)
    F(\epsilon(\beta(g)))$. Applying
    point~\ref{it:uniqueness_in_loc_lie_groupoid-right_unit} of
    Proposition~\ref{prop:uniqueness_in_loc_lie_groupoid} to
    $h:=F(\epsilon(\beta(g)))$, we obtain that
    $F(\epsilon(\beta(g))) = \epsilon'(\beta'(F(g)))$. The other
    identity is obtained in the same way, starting from
    $g=\epsilon(\alpha(g))g$.
  \item It follows from
    point~\ref{it:uniqueness_in_loc_lie_groupoid-right_identity} of
    Proposition~\ref{prop:uniqueness_in_loc_lie_groupoid}, applied to
    $F(g)$ and $ F(g^{-1})$.
  \item If $\ti{F}_0:M\rightarrow M'$ satisfies
    $\alpha' \circ F = \ti{F}_0\circ \alpha$, then
    $\ti{F}_0 = \ti{F}_0\circ \alpha \circ \epsilon = \alpha' \circ F
    \circ \epsilon = F_0$, so that $F_0$ is unique. Being $F_0$ a
    composition of smooth functions, it is smooth. The formulas
    in~\eqref{eq:morphism_of_loc_lie_groupoid-formulas} now follow
    using point~\ref{it:morphism_of_loc_lie_groupoids-F_of_identity}
    of the current result.
  \end{enumerate}
\end{proof}

Just as certain open subsets of Lie groupoids determine local Lie
groupoids, the restriction and co-restriction of morphisms of Lie
groupoids determine morphisms of the corresponding local Lie
groupoids, as shown by the following result.

\begin{lemma}\label{le:restriction_of_morphisms_to_open_subsets}
  Let $F\in \hom_{\LgpdC_M}(G,G')$, $U\subset G$ and $U'\subset G'$ be
  open subsets such that $U^{-1}\subset U$, $(U')^{-1}\subset U'$,
  $\epsilon_G(M)\subset U$ and $\epsilon_{G'}(M)\subset U'$. Assume
  that $F(U)\subset U'$. Then $F|^{U'}_U$ ($F$ restricted to $U$ and
  co-restricted to $U'$) is a morphism of local Lie groupoids, where
  $U$ and $U'$ are local Lie groupoids with the structure given in
  Lemma~\ref{le:open_in_Lgpd_is_lLgpd}. Furthermore,
  $(F|^{U'}_U)_0=id_M$.
\end{lemma}

\begin{lemma}\label{le:morph_and_diffeo_imp_iso_in_lLgpd}
  Let $F\in \hom_{\lLgpdC_{M}}(G,G')$ such that
  $(F\times F)^{-1}((G')_m) \subset G_m$ and that $F:G\rightarrow G'$
  is also a diffeomorphism. Then $F$ is an isomorphism in the
  $\lLgpdC_{M}$ category.
\end{lemma}

%%%%%%%%%%%%%%%%%%%%%%%%%%%%%%%%%%%%%%%%%%%%%%%%%%%%%%%%%%

%\bibliography{math}

%%%%%%%%%%%%%%%%%%%%%%%%%%%%%%%%%%%%%%%%%%%%%%%%%%%%%%%%%%

\newcommand{\etalchar}[1]{$^{#1}$}
\def\cprime{$'$} \def\polhk#1{\setbox0=\hbox{#1}{\ooalign{\hidewidth
  \lower1.5ex\hbox{`}\hidewidth\crcr\unhbox0}}} \def\cprime{$'$}
  \def\cprime{$'$}
\providecommand{\bysame}{\leavevmode\hbox to3em{\hrulefill}\thinspace}
\providecommand{\MR}{\relax\ifhmode\unskip\space\fi MR }
\providecommand{\MRhref}[2]{%
  \href{http://www.ams.org/mathscinet-getitem?mr=#1}{#2}
}
\providecommand{\href}[2]{#2}

%%%%%%%%%%%%%%%%%%%%%%%%%%%%%%%%%%%%%%%%%%%%%%%%%%%%%%%%%%

\end{document}